\documentclass[opre,nonblindrev]{informs3}
\RequirePackage{fix-cm}

\OneAndAHalfSpacedXI

\usepackage{endnotes}
\let\footnote=\endnote

\usepackage{placeins}
\usepackage{natbib}
\bibpunct[, ]{(}{)}{,}{a}{}{,}%
\def\bibfont{\small}%

\TheoremsNumberedThrough     
\ECRepeatTheorems
\EquationsNumberedThrough    

\usepackage[english]{babel}
\usepackage[T1]{fontenc}
\usepackage{amsfonts,amssymb}
\usepackage[sc]{mathpazo}
\usepackage{mathtools}
\usepackage{bbm}
\usepackage{booktabs}
\usepackage{multirow}
\usepackage{rotating}
\usepackage{enumitem}
\usepackage{xcolor}
\usepackage[normalem]{ulem}
\usepackage{tikz}
\usetikzlibrary{arrows.meta}
\definecolor{ctxblue}{HTML}{2F80ED}
\definecolor{ctxcyan}{HTML}{23B5D3}
\definecolor{ctxgreen}{HTML}{31A66A}
\definecolor{ctxorange}{HTML}{F2994A}
\definecolor{ctxpurple}{HTML}{8E5BC7}
\definecolor{ctxred}{HTML}{E25555}
\definecolor{ctxink}{HTML}{263238}
\definecolor{ctxgray}{HTML}{7D8790}
\definecolor{ctxlight}{HTML}{F6F8FB}

\tikzset{
  ctx panel/.style={
    draw=ctxgray!55,
    fill=white,
    rounded corners=2pt,
    line width=0.55pt
  },
  ctx title/.style={
    font=\footnotesize\bfseries,
    text=ctxink,
    align=center
  },
  ctx label/.style={
    font=\scriptsize,
    text=ctxink,
    align=center,
    inner sep=1.5pt
  },
  ctx note/.style={
    font=\scriptsize,
    text=ctxgray!85!black,
    align=center,
    inner sep=1.5pt
  },
  ctx arrow/.style={
    -{Latex[length=2.1mm,width=1.3mm]},
    line width=1.0pt,
    rounded corners=1pt
  },
  ctx thin arrow/.style={
    -{Latex[length=1.7mm,width=1.05mm]},
    line width=0.75pt,
    rounded corners=1pt
  },
  ctx source/.style={
    circle,
    draw=ctxorange!85!black,
    fill=ctxorange!72,
    minimum size=5.2mm,
    inner sep=0pt,
    line width=0.65pt
  },
  ctx mass/.style={
    circle,
    draw=ctxink!75,
    fill=white,
    minimum size=3.2mm,
    inner sep=0pt,
    line width=0.55pt
  },
  ctx box/.style={
    draw=ctxgray!55,
    fill=ctxlight,
    rounded corners=2pt,
    line width=0.55pt,
    align=center,
    font=\scriptsize,
    inner sep=3pt
  },
  ctx equation/.style={
    draw=ctxpurple!60,
    fill=ctxpurple!7,
    rounded corners=2pt,
    line width=0.65pt,
    align=center,
    font=\scriptsize,
    inner sep=3pt
  }
}

\usepackage{comment}
\usepackage{placeins}
\usepackage{hyperref}

\hypersetup{
  colorlinks=true,
  linkcolor=purple,
  citecolor=purple,
  urlcolor=cyan,
  pdftitle={Contextual Distributionally Robust Chance-Constrained Programs: Exact Reformulations and Valid Inequalities},
  pdfauthor={Hanbin Yang, Nan Jiang, Guodong Lyu},
  pdfpagemode=UseOutlines
}

\def\b#1{\boldsymbol{#1}}
\newcommand{\R}{\mathbb R}
\newcommand{\dist}{\mathop{\rm dist}}
\newcommand{\pos}[1]{\left[#1\right]_+}

\usepackage{natbib}
\bibpunct[, ]{(}{)}{,}{a}{}{,}
\def\bibfont{\small}

\begin{document}

\RUNAUTHOR{Yang, Jiang, and Lyu}

\renewcommand{\TITLEfont}{\fontsize{18}{22}\selectfont}
\TITLE{Contextual Distributionally Robust Chance-Constrained\\Programs: Exact Reformulations and Valid Inequalities}
\RUNTITLE{Contextual Distributionally Robust Chance-Constrained Programs}

\ARTICLEAUTHORS{
\AUTHOR{Hanbin Yang, Nan Jiang, Guodong Lyu}
\AFF{The Hong Kong University of Science and Technology}
}

\ABSTRACT{%
We study a contextual distributionally robust chance-constrained programming (C-DRCCP) model that controls violation probabilities conditional on a neighborhood of a target context. The ambiguity set is defined over the joint distribution of contexts and outcomes using the Wasserstein distance. For affine safety systems with right-hand-side uncertainty, we derive an exact mixed-integer programming (MIP) reformulation. The reformulation reveals that worst-case conditional risk is determined by two coupled uses of the transportation budget: moving empirical mass into the conditioning neighborhood and moving conditioned mass into the failure region. Equivalently, robust feasibility holds when the minimum transportation cost needed to reach the conditional risk boundary is at least the available budget. Calibrating the ambiguity set gives a finite-sample guarantee for the true conditional violation probability.

\noindent This characterization leads to strong valid inequalities that exploit the contextual transport. Each contextual mass allocation determines both the conditional scenario weights and the budget available to move outcomes into failure, whereas conventional quantile cuts utilize only the induced weights. Accounting for the remaining budget generates \emph{strengthened quantile cuts} that strictly dominate standard quantile bounds and yield tighter MIP coefficients. We also derive \emph{rank inequalities} from an exchange property of contextual transportation costs, linking feasible failure patterns across contextual allocations. Extensive numerical experiments demonstrate that these inequalities close much of the linear-relaxation gap, speed up the solution process, and enable large-sample instances that are intractable for the baseline formulation. Using a real-world anesthesiologist deployment case study, we further reveal a statistical trade-off in contextual robustness. Conditioning on the target context improves relevance by focusing on observations from samples with similar contexts. When sufficient data are available, the C-DRCCP model better differentiates risk across contexts and reallocates protection accordingly, improving the balance between cost and reliability. 
}




\MANUSCRIPTNO{}

\renewcommand{\theARTICLETOP}{}
\RRHFirstLine{}
\LRHFirstLine{}
\RRHSecondLine{\itshape\theRUNTITLE}
\LRHSecondLine{\itshape\theRUNTITLE}
\ECRRHSecondLine{\itshape\theRUNTITLE}
\ECLRHSecondLine{\itshape\theRUNTITLE}

\maketitle


\section{Introduction}\label{sec:introduction}
Side information is increasingly available before operational decisions are made~\citep{bertsimas2020predictive}. Such information includes demand forecasts, weather conditions, customer attributes, market indicators, and other covariates that help characterize the uncertainty entering operational constraints. For example, \citet{mao2026predictive} illustrate how demand distributions depend on exogenous features and service attributes, and emphasize that ignoring the resulting heteroskedasticity can degrade operational decisions. 

Contextual optimization can utilize the side information for decision-making, which is especially important in chance-constrained problems, where feasibility depends on the distribution under which the constraint is evaluated.  A decision that appears feasible under the marginal distribution may still violate the desired reliability level under the conditional distribution relevant to the realized context~\citep{rahimian2023contextual}. Therefore, ignoring contextual information can lead to decisions that are feasible on average but unreliable in specific regimes. This motivates methods that control conditional violation probabilities, particularly in applications where violations have serious consequences.

Let $X\in\mathcal{X}$ denote the side information and $Y\in\mathcal{Y}$ the uncertain outcome. The decision incurs cost $c(z)$ and is safe for outcomes in $\mathcal{S}(z)\subseteq\mathcal{Y}$. We collect $N$ joint samples $\{(\hat{x}_i,\hat{y}_i)\}_{i\in[N]}$, where $\hat{x}_i$ is a realized context, $\hat{y}_i$ is the associated outcome, and $[N]\coloneqq\{1,\ldots,N\}$. Their empirical distribution is $\widehat{\mathbb{P}}_N=N^{-1}\sum_{i\in[N]}\delta_{(\hat{x}_i,\hat{y}_i)}$. To account for distributional uncertainty, we consider an ambiguity set of joint distributions around this empirical distribution. We measure the cost of moving a joint observation $(x,y)$ to $(x',y')$ using the additive ground metric $\mathrm{dist}((x,y),(x',y'))=d_{\mathcal{X}}(x,x')+d_{\mathcal{Y}}(y,y')$, where $d_{\mathcal{X}}$ and $d_{\mathcal{Y}}$ measure distances in the context and outcome spaces, respectively. For joint distributions $\mathbb{P}_1$ and $\mathbb{P}_2$ on $\mathcal{X}\times\mathcal{Y}$, let $W$ denote the associated type-1 Wasserstein distance, i.e., $W(\mathbb{P}_1,\mathbb{P}_2)\coloneqq\inf_{\pi\in\Pi(\mathbb{P}_1,\mathbb{P}_2)}\int \mathrm{dist}(\xi,\xi')\,\pi(\mathrm{d}\xi,\mathrm{d}\xi')$, where $\Pi(\mathbb{P}_1,\mathbb{P}_2)$ is the set of probability measures on $(\mathcal{X}\times\mathcal{Y})^2$ with marginals $\mathbb{P}_1$ and $\mathbb{P}_2$, and $\xi,\xi'\in\mathcal{X}\times\mathcal{Y}$ denote joint realizations. Here, we allow candidate joint distributions $\mathbb{P}$ to lie within a Wasserstein radius $\theta>0$ of $\widehat{\mathbb{P}}_N$.

Because contextual information can alter the relevant outcome distribution~\citep{sim2025analytics}, we seek to assess risk under the conditional distribution of $Y$. However, conditional on a null set, the distribution of $Y$ can be specified arbitrarily without changing the joint distribution and is not uniquely determined by it (see Proposition~\ref{prop:null-kernel}). Motivated by the use of nearby contexts to characterize the local outcome distribution in~\citet{nguyen2024robustifying}, we condition on a prescribed \emph{fiber} set $\mathcal{N}$, and impose a lower bound $\varepsilon\in(0,1]$ on the probability of this event. These requirements define the ambiguity set
\begin{equation}\label{eq:ambiguity-set}
\mathcal{P}_{\theta,\varepsilon}(\mathcal{N})
\coloneqq
\left\{\mathbb{P}\;\middle|\;
W(\mathbb{P},\widehat{\mathbb{P}}_N)\leq\theta,\ 
\mathbb{P}_X(\mathcal{N})\geq\varepsilon
\right\},
\end{equation}
where $\mathbb{P}_X$ denotes the context marginal of $\mathbb{P}$. Given a prescribed limit on the conditional violation probability $\epsilon\in(0,1)$, we minimize cost while requiring this limit to hold under every distribution in $\mathcal{P}_{\theta,\varepsilon}(\mathcal{N})$. After observing a target context $x_0\in\mathcal{X}$, we choose a decision $z$ from the feasible decision set $\mathcal{Z}\subseteq\R^{\dim(z)}$. This gives the \emph{contextual} distributionally robust chance-constrained programming (C-DRCCP) model
\begin{subequations}\label{opt:drc3p}
\begin{align}
\min_{z\in\mathcal{Z}}\quad &c(z) \label{eq:drc3p-objective}\\
\mathrm{s.t.}\quad&
\sup_{\mathbb{P}\in\mathcal{P}_{\theta,\varepsilon}(\mathcal{N})}
\mathbb{P}\left[Y\notin\mathcal{S}(z)\mid X\in\mathcal{N}\right]
\leq\epsilon, \label{eq:drc3p-chance}
\end{align}
\end{subequations}
where $\mathcal{N}=\mathcal{N}_\gamma(x_0)$ is a neighborhood of $x_0$, and the parameter $\gamma\geq0$ controls the degree of localization. The chance constraint controls the probability that at least one safety requirement fails. The ambiguity set on the joint distribution of $(X,Y)$ avoids the need to \emph{estimate} a conditional distribution \emph{before} optimization. 

\subsection{Literature Review}
The contextual optimization literature examines how observed side information can improve decisions by estimating the conditional distribution relevant to the current setting~\citep{sadana2025survey}. Nonparametric prescriptions use local weights or other prediction methods to connect observations to decisions~\citep{ban2019newsvendor,bertsimas2020predictive,rahimian2023contextual}. Residual methods first estimate the covariate effect and then construct an ambiguity set for the unexplained uncertainty~\citep{kannan2024residuals}. \citet{zhang2026managing} integrate prediction and optimization into a moment-based distributionally robust optimization (DRO) model. The ambiguity model of~\citet{nguyen2024robustifying} is closest to ours. They study worst-case conditional expectations for convex portfolio losses and exploit convexity to derive conic reformulations. Related conditional DRO models infer local distributions through trimmings and partial transport~\citep{esteban2022trimmings}, connect several convex conditional risk constructions through unions of transport balls~\citep{xie2025conditional}, or preserve the information structure of a decision rule through causal transport~\citep{yang2026causal}. They focus on conditional expected losses or convex risk measures. In contrast, model~\eqref{opt:drc3p} imposes a chance constraint through a discontinuous indicator, and the ambiguous distribution controls both the failure mass and the conditioning mass. Thus, the reformulations for convex losses cannot directly apply to the mixed-integer programming (MIP) reformulation required here.

DRO protects decisions against uncertainty in an estimated distribution. Ambiguity sets based on moments, marginal distributions, and other structural information can capture asymmetry and dependence not captured by a nominal model~\citep{natarajan2008asymmetric,natarajan2010tractable,chen2019infinitely,chen2022marginals}. Wasserstein ambiguity sets instead place distributions near an empirical reference distribution and support tractable dual reformulations and statistical guarantees~\citep{mohajerin2018data,gao2023distributionally}. For chance constraints over Wasserstein balls, exact reformulations are available for general safety sets and for important classes of right-hand and left-hand side uncertainty~\citep{xie2021drcc,chen2022data,honguyen2022rhs,honguyen2023strong}. Convex approximations can be useful, but the common approximations are generally incomparable and can be overconservative~\citep{chen2023approximations,zhou2026sfla}. Application-driven work has developed MIP formulations and constraint generation for distributionally robust chance-constrained hub location~\citep{zhao2023hub}. These results explain both the statistical appeal of DRCCP and the need for exact formulations when the chance constraint is computationally difficult.

Chance constraints under finite-support distributions have been studied extensively in integer programming. Binary variables indicating scenario violations yield an MIP formulation, whose continuous relaxation can be strengthened using mixing inequalities~\citep{luedtke2010integer,kucukyavuz2012mixing,abdi2016mixing,kilincc2022joint}. \citet{song2014packing} derive probabilistic cover inequalities, coefficient strengthening, and a projection procedure for chance-constrained binary packing problems. \citet{ahmed2017nonanticipative} use nonanticipativity relaxations to obtain stronger dual bounds, two primal MIP formulations, and exact algorithms. Branch-and-cut decomposition and quantile closures provide more general solution tools~\citep{luedtke2014branch,xie2018quantile}. In the Wasserstein setting,~\citet{jiang2024terminator} combine inner and outer approximations through variable fixing to reduce the number of active scenario binaries, while~\citet{pathy2025decomposition} derive probability cuts and mixing inequalities for finite-support DRCCP with polyhedral probability ambiguity. Our formulation builds on this integer-programming foundation, and its probability weights are derived from contextual transport. Normalizing those weights could \emph{discard} the remaining transport budget, which is precisely the information our strengthened cuts use.

\subsection{Main Results and Contributions}
We derive tractable reformulations of the C-DRCCP model and develop valid inequalities to improve computational efficiency. Our main contributions are as follows.
\begin{enumerate}[leftmargin=2em, wide]
\item \textbf{Exact Reformulation and Finite-Sample Guarantees.} We derive an exact MIP reformulation of the C-DRCCP model for affine polyhedral safety systems with right-hand-side uncertainty. The derivation uses a linear-programming (LP) representation of the adversary's mass allocation problem for each fixed decision and a compact dual certificate. We also establish a finite-sample guarantee: when the Wasserstein radius and the lower bound on neighborhood probability are calibrated at prescribed confidence levels, any decision selected from the data and feasible for the model satisfies the target risk limit conditional on the neighborhood with high confidence.

\item \textbf{Valid Inequalities and Computational Benefits.} We develop valid inequalities by exploiting the coupling between contextual mass allocation and outcome transportation. Each contextual mass allocation determines both the conditional scenario probabilities and the budget available to move outcomes into failure. Accounting for this remaining budget yields \emph{strengthened quantile cuts} that strictly dominate those based solely on conditional scenario probabilities and provide tighter coefficients for the MIP formulation. We further exploit the ordering of contextual transportation costs and the geometry of feasible mass allocations to derive \emph{rank inequalities} and compact extended formulations.

\item[] The numerical experiments show that these valid inequalities substantially strengthen the root relaxation and accelerate the solution process. In the capacitated transportation problem, our strengthened formulation achieves a mean root gap of $0.69\%$ across 108 instances, compared with $89.48\%$ for the baseline. The advantages are particularly pronounced for instances with larger training samples: with $N=1{,}000$, the strengthened formulation solves all 12 instances, with a median solution time of $39.67$ seconds, whereas the baseline reaches the one-hour time limit on every instance without attaining that gap.

\item \textbf{Case Study and Managerial Insights.} 
 Using a real-world anesthesiologist deployment case study, we illustrate a statistical trade-off: conditioning on the target context improves relevance by focusing on samples with similar contexts, but reduces the \emph{effective} sample size for estimating local tail risk. Thus, the benefit of contextual information depends on whether the available data support reliable estimation within the conditioning neighborhood. With sufficient data, C-DRCCP can better distinguish risk across contexts and allocate protection accordingly, making less conservative decisions when demand is low and providing greater protection when demand is high. This can improve the balance between cost and reliability.

\end{enumerate}

The remainder of the paper proceeds as follows. Section~\ref{sec:model} justifies the conditioning requirement in the ambiguity set and characterizes contextual transport geometry. Section~\ref{sec:exact} derives the exact MIP reformulation, characterizes robust feasibility through the minimum transportation cost of reaching the conditional risk boundary, and establishes statistical guarantees. Section~\ref{sec:strengthening} develops coefficient tightening and valid inequalities for strengthening the reformulations. Section~\ref{sec:computational} reports numerical experiments on a capacitated transportation problem and an anesthesiologist deployment case study. Finally, Section~\ref{sec:conclusion} summarizes the paper. Appendix~\ref{appendix:model-relations} characterizes the relation to the DRCCP model. Appendix~\ref{sec:additionaltables} reports additional numerical results, including computational comparisons of alternative implementations of the exact inequality reformulation.

\section{Contextual Wasserstein Ambiguity}\label{sec:model}
This section justifies the conditioning structure and the joint Wasserstein ambiguity set in model~\eqref{opt:drc3p}. We first show the motivation for conditioning on a positive-probability neighborhood. We then characterize the least contextual transportation needed to create the required conditioning mass, which explains the roles of the neighborhood radius and provides the contextual cost terms used in the exact reformulations.

\subsection{Neighborhood Conditioning}\label{sec:neighborhood-conditioning}

Let $X$ and $Y$ take values in standard Borel spaces $\mathcal{X}$ and $\mathcal{Y}$, respectively. \citet[Theorem~8.5]{Kallenberg2021} states that a regular conditional distribution of $Y$ given $X$ exists and that any two versions of this conditional distribution agree $\mathbb{P}_X$-almost everywhere. This \emph{almost-everywhere} uniqueness, however, does not determine the conditional distribution at a prescribed context having zero probability under $\mathbb{P}_X$. The following proposition shows that the conditional distribution at such a context can be changed arbitrarily without changing the joint distribution of $(X,Y)$.

\begin{proposition}\label{prop:null-kernel}
Let $\mathbb{P}$ be a probability distribution on $\mathcal{X}\times\mathcal{Y}$, and suppose that $\mathbb{P}_X(\{x_0\})=0$. If $K$ is a regular conditional distribution of $Y$ given $X$, then, for every probability distribution $Q$ on $\mathcal{Y}$, there exists another regular conditional distribution $K^Q$ that agrees with $K$ on $\mathcal{X}\setminus\{x_0\}$ and satisfies $K^Q(x_0,\cdot)=Q$.
\end{proposition}

Proposition~\ref{prop:null-kernel} shows that the value of $\mathbb{P}(Y\notin\mathcal{S}(z)\mid X=x_0)$ is not determined by the joint distribution whenever $x_0$ is a null point. Thus, we condition on the event $X\in\mathcal{N}_\gamma(x_0)$, where $\mathcal{N}_\gamma(x_0)$ is a Borel neighborhood of the target context $x_0$ and $\gamma\geq 0$ controls its size. Throughout, $x_0\in\mathcal{N}_\gamma(x_0)$, $\sup_{x\in\mathcal{N}_\gamma(x_0)}d_{\mathcal{X}}(x,x_0)\leq\gamma$, and $\mathcal{N}_0(x_0)=\{x_0\}$. The parameters $\gamma$ and $\varepsilon$ serve different purposes. The neighborhood radius $\gamma$ controls localization: reducing it focuses the risk requirement more closely on $x_0$ but leaves less local information. The mass lower bound $\varepsilon$ ensures that each distribution in $\mathcal{P}_{\theta,\varepsilon}(\mathcal{N})$ assigns positive probability to the conditioning event, so the conditional probability is well defined. 

\begin{assumption}\label{ass:transport}
We impose the following assumptions on the C-DRCCP model~\eqref{opt:drc3p}:
\begin{enumerate}[label=(\roman*),leftmargin=2em]
\item $(\mathcal{X},d_{\mathcal{X}})$ and $(\mathcal{Y},d_{\mathcal{Y}})$ are Polish metric spaces whose metrics generate the stated Borel $\sigma$-fields, and transport is allowed on the product $\mathcal{X}\times\mathcal{Y}$.
\item The set $\mathcal{N}_\gamma(x_0)\subseteq \mathcal{X}$ is measurable. The distances from each $\hat{x}_i$ to $\mathcal{N}_\gamma(x_0)$ and $\mathcal{N}_\gamma(x_0)^c$ are finite.
\item For $z\in\mathcal{Z}$, the failure set $\mathcal{S}^c(z)\subseteq\mathcal{Y}$ is a nonempty Borel set, and $\dist_{\mathcal{Y}}(\hat{y}_i,\mathcal{S}^c(z))$ is finite for $i\in[N]$.
\end{enumerate}
\end{assumption}
Assumption~\ref{ass:transport} guarantees the existence of an optimal coupling under the continuous ground metric~\citep{villani2009optimal}. However, the closest points in $\mathcal{N}$, $\mathcal{N}^c$, or $\mathcal{S}^c(z)$ may not exist, and we work with infimum distances and later impose a strict radius condition that provides the slack needed for approximate projections.

\subsection{Contextual Transport Geometry}\label{sec:geometry}

We now translate the joint ambiguity model into the contextual transportation costs that enter the reformulations. For every sample $i\in[N]$, define its inward and outward contextual distances by $\kappa_i^{\mathrm{in}}\coloneqq\inf_{x\in\mathcal{N}}d_{\mathcal{X}}(\hat{x}_i,x)$ and $\kappa_i^{\mathrm{out}}\coloneqq\inf_{x\in\mathcal{N}^c}d_{\mathcal{X}}(\hat{x}_i,x)$. If $w_i\in[0,1/N]$ units of mass from sample $i$ are assigned to the conditioning neighborhood $\mathcal{N}_\gamma(x_0)$, then the infimum contextual transportation cost attributable to that sample is
\[w_i\kappa_i^{\mathrm{in}}+\left(\frac{1}{N}-w_i\right)\kappa_i^{\mathrm{out}}=\frac{1}{N}\kappa_i^{\mathrm{out}}+\bar{\kappa}_iw_i,\]
where $\bar{\kappa}_i\coloneqq\kappa_i^{\mathrm{in}}-\kappa_i^{\mathrm{out}}$ represents the change in transportation cost when mass from sample $i$ is assigned inside. Define the minimum transportation cost required to assign at least $\varepsilon$ probability mass to the conditioning neighborhood as $\theta_{\min}(x_0,\gamma,\varepsilon)\coloneqq\inf\{W(\mathbb{P},\widehat{\mathbb{P}}_N) \;|\;\mathbb{P}_X(\mathcal{N}_\gamma(x_0))\geq\varepsilon\}$. Proposition~\ref{prop:min-radius} establishes the equivalent LP reformulation for computing $\theta_{\min}(x_0,\gamma,\varepsilon)$, as in~\citet{nguyen2024robustifying}.
\begin{proposition}\label{prop:min-radius}
Under Assumption~\ref{ass:transport}, the value $\theta_{\min}(x_0,\gamma,\varepsilon)$ is the optimal value of the linear program
\begin{equation}\label{eq:min-radius-lp}
\theta_{\min}(x_0,\gamma,\varepsilon)
=
\min_{\substack{0\leq w_i\leq1/N,\ i\in[N]\\
\sum_{i\in[N]}w_i\geq\varepsilon}}
\left\{
K_0+\sum_{i\in[N]}\bar{\kappa}_iw_i
\right\},
\end{equation}
where $K_0\coloneqq N^{-1}\sum_{i\in[N]}\kappa_i^{\mathrm{out}}$. Consequently, $\mathcal{P}_{\theta,\varepsilon}(\mathcal{N})$ is nonempty whenever $\theta>\theta_{\min}(x_0,\gamma,\varepsilon)$, and it is empty whenever $\theta<\theta_{\min}(x_0,\gamma,\varepsilon)$.
\end{proposition}

Problem~\eqref{eq:min-radius-lp} can be viewed as a fractional knapsack problem. Indeed, we can first keep all mass with $\bar{\kappa}_i<0$, then, if necessary, fill the remaining mass requirement in increasing order of $\bar{\kappa}_i$, splitting at most one sample. Thus, $\theta_{\min}$ can be computed in $O(N\log N)$ time. Throughout the remainder of the paper, we impose the strict condition
\begin{equation}\label{eq:strict-radius}
\theta>\theta_{\min}(x_0,\gamma,\varepsilon).
\end{equation}
Besides guaranteeing that the ambiguity set is nonempty, condition~\eqref{eq:strict-radius} supplies the small transportation slack needed when a closest point in $\mathcal{N}^c$ or $\mathcal{S}^c(z)$ does not exist. When $\mathcal{X}$ is the entire normed linear space and $\mathcal{N}_\gamma(x_0)=\{x \mid \|x-x_0\|_{\mathcal{X}}\leq\gamma\}$ is its closed norm ball, the contextual coefficients are available in closed form: $\kappa_i^{\rm in}=[\|\hat{x}_i-x_0\|_{\mathcal{X}}-\gamma]_+$, $\kappa_i^{\rm out}=[\gamma-\|\hat{x}_i-x_0\|_{\mathcal{X}}]_+$, and $\bar{\kappa}_i=\|\hat{x}_i-x_0\|_{\mathcal{X}}-\gamma$. Thus, all coefficients that enter the minimum-radius problem can be computed directly from the sample distances to $x_0$. 

\section{Exact Reformulations and Theoretical Properties}\label{sec:exact}
Consider an affine safety system with right-hand-side uncertainty at the observed context $x_0$. Let $[P]\coloneqq\{1,\ldots,P\}$ index the safety constraints, assume $\mathcal{Y}=\R^{\dim(y)}$, and let $\b a_p(x_0)\in\R^{\dim(z)}$, $\b b_p(x_0)\in\R^{\dim(y)}$, and $d_p(x_0)\in\R$ be known coefficients for each $p\in[P]$. The safety set is
\begin{equation}\label{eq:polyhedral-safety}
\mathcal{S}(z)
\coloneqq
\left\{
y\in\R^{\dim(y)} \;\middle|\;
y^\top\b b_p(x_0)+d_p(x_0)-z^\top\b a_p(x_0)>0
\quad\forall p\in[P]
\right\},
\end{equation}
which covers a broad class of operational reliability requirements. For example, it captures demand-coverage requirements in capacitated transportation~\citep{chen2022data}, vehicle routing~\citep{zhang2021robust}, and unit commitment~\citep{zhou2026sfla}, where shipments, routed capacity, or committed generation capacity must cover the corresponding uncertain demand. It also captures staffing requirements in healthcare management, where scheduled capacity together with a planned allowance must cover workload~\citep{rath2026staff}. 

The strict inequalities classify boundary points as failures, so the associated failure set $\mathcal{S}^c(z)=\bigcup_{p\in[P]}\{y\mid y^\top\b b_p(x_0)+d_p(x_0)-z^\top\b a_p(x_0)\leq0\}$ is closed. For a decision $z\in\mathcal{Z}$ and each sample $i\in[N]$, define $\Delta_i(z)\coloneqq\dist_{\mathcal{Y}}\left(\hat{y}_i,\mathcal{S}^c(z)\right)=\inf_{y\in\mathcal{S}^c(z)}d_{\mathcal{Y}}(\hat{y}_i,y)$. Since $\mathbb{P}_X(\mathcal{N})\geq\varepsilon>0$ for every $\mathbb{P}\in\mathcal{P}_{\theta,\varepsilon}(\mathcal{N})$, we have
\[\mathbb{P}\left(Y\in\mathcal{S}^c(z)\mid X\in\mathcal{N}\right)\leq\epsilon\quad\Longleftrightarrow\quad\mathbb{P}\left((X,Y)\in\mathcal{N}\times\mathcal{S}^c(z)\right)-\epsilon\mathbb{P}_X(\mathcal{N})\leq0,\]
which avoids direct dualization of a probability ratio. 

\subsection{LP Reformulation Based on Mass Allocation}\label{sec:mass-flow}
We first establish an exact LP reformulation from the transport plan. For each empirical sample, the formulation tracks the mass transported into the conditioning neighborhood $\mathcal{N}$ and the portion of that mass transported into the failure set $\mathcal{S}^c(z)$. This decomposition reveals the trade-off faced by the adversary: transporting mass into $\mathcal{N}$ changes the denominator of the conditional violation probability, whereas transporting part of that mass into $\mathcal{S}^c(z)$ increases the numerator and incurs an additional outcome transportation cost. This interpretation also provides the foundation for the dual and MIP reformulations that follow.

\begin{theorem}\label{thm:mass-flow}
Suppose that Assumption~\ref{ass:transport} and condition~\eqref{eq:strict-radius} hold. For every fixed $z\in\mathcal{Z}$,
\begin{align}
\sup_{\mathbb{P}\in\mathcal{P}_{\theta,\varepsilon}(\mathcal{N})}
\left\{
\mathbb{P}\left((X,Y)\in\mathcal{N}\times\mathcal{S}^c(z)\right)
-\epsilon\mathbb{P}_X(\mathcal{N})
\right\} = \max_{w,r}\quad &
\sum_{i\in[N]}r_i-\epsilon\sum_{i\in[N]}w_i
\label{eq:mass-flow}\\
\mathrm{s.t.}\quad&
K_0+\sum_{i\in[N]}\bar{\kappa}_iw_i
+\sum_{i\in[N]}\Delta_i(z)r_i\leq\theta \notag\\
&\sum_{i\in[N]}w_i\geq\varepsilon\notag\\
&0\leq r_i\leq w_i\leq\frac{1}{N} && \forall i\in[N]. \notag
\end{align}
\end{theorem}

We interpret the variables in~\eqref{eq:mass-flow} as follows: variable $w_i$ is the mass from source $i$ that reaches the conditioning neighborhood, while $r_i$ is the portion of that mass that also reaches failure. The term $K_0+\sum_{i\in[N]}\bar{\kappa}_iw_i$ is the cost of moving mass into the conditioning neighborhood, and $\sum_{i\in[N]}\Delta_i(z)r_i$ is the additional outcome cost needed to create failure. The adversary also chooses the total probability mass assigned to the conditioning neighborhood.

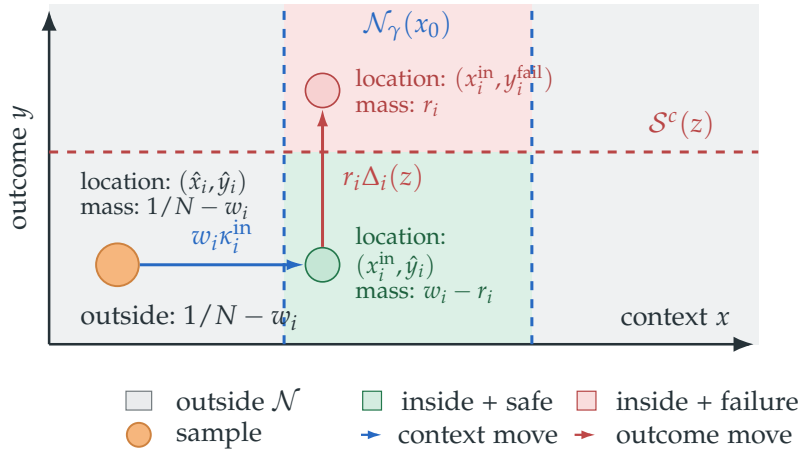
\begin{figure}[htbp]
  \centering
  \caption{Allocation of mass from an empirical distribution.}\label{fig:contextual-mass-flow}
  \resizebox{0.65\linewidth}{!}{%
  \begin{tikzpicture}[x=1cm,y=1cm,>=Latex]
    \begin{scope}
      \clip (0.75,1.45) rectangle (9.05,5.43);
      \fill[ctxgray!12] (0.75,1.45) rectangle (9.05,5.43);
      \fill[ctxgreen!17] (3.50,1.45) rectangle (6.40,3.70);
      \fill[ctxred!16] (3.50,3.70) rectangle (6.40,5.43);
    \end{scope}

    \draw[->,line width=0.85pt,ctxink]
      (0.75,1.45) -- (9.00,1.45);
    \node[anchor=south east,font=\footnotesize,text=ctxink]
      at (8.86,1.53) {context $x$};
    \draw[->,line width=0.85pt,ctxink]
      (0.75,1.45) -- (0.75,5.36);
    \node[rotate=90,anchor=south,font=\footnotesize,text=ctxink]
      at (0.7,3.52) {outcome $y$};

    \draw[ctxblue!82!black,dashed,line width=1.0pt]
      (3.50,1.45) -- (3.50,5.43);
    \draw[ctxblue!82!black,dashed,line width=1.0pt]
      (6.40,1.45) -- (6.40,5.43);
    \draw[ctxred!82!black,dashed,line width=1.0pt]
      (0.75,3.70) -- (9.05,3.70);
    \node[font=\footnotesize\bfseries,text=ctxblue!82!black]
      at (4.95,5.18) {$\mathcal{N}_\gamma(x_0)$};
    \node[anchor=west,font=\footnotesize\bfseries,
      text=ctxred!82!black] at (7.63,4.02) {$\mathcal{S}^c(z)$};

    \node[ctx source,minimum size=5.0mm] (si) at (1.55,2.38) {};
    \node[anchor=south west,align=left,font=\scriptsize,text=ctxink]
      at (0.98,2.78)
      {location: $(\hat{x}_i,\hat{y}_i)$\\[-1pt]
       mass: $1/N-w_i$};
    \node[ctx mass,minimum size=4.0mm,draw=ctxgreen!72!black,
      fill=ctxgreen!28] (isafe) at (3.95,2.38) {};
    \node[ctx mass,minimum size=4.0mm,draw=ctxred!80!black,
      fill=ctxred!27] (ifail) at (3.95,4.42) {};

    \draw[ctx arrow,ctxblue!82!black] (si) -- (isafe)
      node[midway,above=2pt,font=\footnotesize\bfseries,
        text=ctxblue!82!black,fill opacity=0.94,text opacity=1,
        inner sep=1.2pt]
      {$w_i\kappa_i^{\rm in}$};
    \draw[ctx arrow,ctxred!84!black] (isafe) -- (ifail)
      node[midway,right=5pt,font=\footnotesize\bfseries,
        text=ctxred!84!black,fill opacity=0.94,text opacity=1,
        inner sep=1.2pt]
      {$r_i\Delta_i(z)$};

    \node[anchor=north west,font=\footnotesize,text=ctxink] at (0.98,2.05) {outside: $1/N-w_i$};
    \node[anchor=west,align=left,font=\scriptsize,text=ctxgreen!58!black,
      text width=2.05cm] at (4.20,2.38)
      {location: $(x_i^{\rm in},\hat{y}_i)$\\[-1pt]
       mass: $w_i-r_i$};
    \node[anchor=west,align=left,font=\scriptsize,text=ctxred!78!black,
      text width=3.05cm] at (4.20,4.42)
      {location: $(x_i^{\rm in},y_i^{\rm fail})$\\[-1pt]
       mass: $r_i$};

    \filldraw[fill=ctxgray!16,draw=ctxink!65] (1.65,0.68) rectangle (1.95,0.90);
    \node[anchor=west,font=\footnotesize,text=ctxink] at (2.10,0.79) {outside $\mathcal{N}$};
    \filldraw[fill=ctxgreen!21,draw=ctxgreen!72!black]
      (4.42,0.68) rectangle (4.64,0.90);
    \node[anchor=west,font=\footnotesize,text=ctxink] at (4.74,0.79)
      {inside + safe};
    \filldraw[fill=ctxred!20,draw=ctxred!80!black]
      (6.93,0.68) rectangle (7.15,0.90);
    \node[anchor=west,font=\footnotesize,text=ctxink] at (7.25,0.79)
      {inside + failure};

    \node[ctx source,minimum size=2.8mm] at (1.80,0.39) {};
    \node[anchor=west,font=\footnotesize,text=ctxink] at (2.10,0.39)
      {sample};
    \draw[ctx thin arrow,ctxblue!82!black]
      (4.40,0.39) -- (4.65,0.39);
    \node[anchor=west,font=\footnotesize,text=ctxink] at (4.70,0.39)
      {context move};
    \draw[ctx thin arrow,ctxred!84!black]
      (6.90,0.39) -- (7.15,0.39);
    \node[anchor=west,font=\footnotesize,text=ctxink] at (7.25,0.39)
      {outcome move};

  \end{tikzpicture}%
  }
  \par\footnotesize Note. $x_i^{\rm in}$ and $y_i^{\rm fail}$ denote representative destination points. Mass $w_i$ enters $\mathcal{N}_\gamma(x_0)$ and its nested portion $r_i$ reaches failure.
\end{figure}

Figure~\ref{fig:contextual-mass-flow} summarizes the resulting geometry. In an unconditional chance constraint, every empirical sample contributes the fixed mass $1/N$, so the adversary only decides how much mass to move toward failure and prioritizes small outcome distances. Here, the adversary also constructs the conditioning event. It may favor a sample with a larger outcome distance when moving its context into $\mathcal{N}$ is inexpensive, or reject a sample with a smaller outcome distance when moving its context into $\mathcal{N}$ is costly. The mass allocations are coupled: $w_i$ creates conditional denominator mass, $r_i\leq w_i$ creates numerator mass, and both compete for the same transportation budget. The two movements in the figure describe this additive cost decomposition, not a sequential restriction on the joint choice of $(w,r)$.

In particular, the chance constraint~\eqref{eq:drc3p-chance} holds if and only if the optimal value of~\eqref{eq:mass-flow} is nonpositive. The next result follows from the LP duality applied to Theorem~\ref{thm:mass-flow}. In the unnormalized dual, the dual multiplier of the Wasserstein budget constraint is multiplied by each distance $\Delta_i(z)$, which is decision-dependent. In contrast, normalization removes these bilinear terms, yielding the compact certificate below.

\begin{theorem}\label{thm:dual-certificate}
Suppose that Assumption~\ref{ass:transport} and condition~\eqref{eq:strict-radius} hold. For a fixed $z\in\mathcal{Z}$, the chance constraint~\eqref{eq:drc3p-chance} holds if and only if there exist $t\geq0$, $v_i\geq0$, $s_i\in\R$, and $\lambda\in\R$ such that
\begin{subequations}\label{eq:author-certificate}
\begin{align}
& K_0+\epsilon\varepsilon t-\varepsilon\lambda-\frac{1}{N}\sum_{i\in[N]}v_i\geq\theta \label{eq:certificate-budget}\\
& v_i+\lambda \geq s_i &&\forall i\in[N] \label{eq:certificate-link}\\
& \Delta_i(z)+\bar{\kappa}_i \geq t-s_i &&\forall i\in[N] \label{eq:certificate-distance}\\
& s_i\geq-\bar{\kappa}_i &&\forall i\in[N] \label{eq:certificate-s-domain}\\
& \lambda \leq\epsilon t. \label{eq:certificate-lambda-domain}
\end{align}
\end{subequations}
\end{theorem}

While the ambiguity set follows the joint optimal transport construction of~\citet{nguyen2024robustifying}, their results rely on the conic portfolio reformulations and exploit the structure of mean-variance and mean-CVaR losses.
Starting from the allocation problem~\eqref{eq:mass-flow}, Theorem~\ref{thm:dual-certificate} reformulates a conditional chance constraint, in which the decision determines a discontinuous failure set, and the adversary chooses both the mass entering the conditioning neighborhood and the portion of that mass entering failure. The following corollary simplifies Theorem~\ref{thm:dual-certificate} to the singleton conditioning set $\mathcal{N}=\{x_0\}$.

\begin{corollary}\label{cor:singleton}
Let $\mathcal{N}=\{x_0\}$, and suppose that Assumption~\ref{ass:transport} and condition~\eqref{eq:strict-radius} hold. Suppose that $\inf_{x\in\mathcal{X}\setminus\{x_0\}}d_{\mathcal{X}}(x_0,x)=0$. Define $\kappa_i=d_{\mathcal{X}}(\hat{x}_i,x_0)$. Then $K_0=0$, $\bar{\kappa}_i=\kappa_i$, and the robust conditional chance constraint is equivalent to the existence of $t\in\R_+$, $v\in\R_+^N$, $s\in\R^N$, and $\lambda\in\R$ satisfying
\begin{subequations}\label{eq:singleton-certificate}
\begin{align}
& \epsilon\varepsilon Nt-\varepsilon N\lambda-\sum_{i\in[N]} v_i \geq\theta N\\
& v_i+\lambda\geq s_i &&\forall i\in[N]\\
& \Delta_i(z)+\kappa_i\geq t-s_i &&\forall i\in[N]\\
& s_i\geq-\kappa_i &&\forall i\in[N]\\
& \lambda\leq\epsilon t.
\end{align}
\end{subequations}
\end{corollary}

\subsection{Exact MIP Reformulation}\label{sec:mip}
We now use the affine safety system~\eqref{eq:polyhedral-safety} to express the distances $\Delta_i(z)$ explicitly. Assume that $d_{\mathcal{Y}}$ is induced by a norm $\|\cdot\|_{\mathcal{Y}}$ on $\R^{\dim(y)}$, let $\|\cdot\|_{\mathcal{Y},*}$ be its dual norm, and suppose that $\b b_p(x_0)\neq0$ for all $p\in[P]$. Define the normalized sample margin
\begin{equation}\label{eq:normalized-margin}
m_{ip}(z)
\coloneqq
\frac{\hat{y}_i^\top\b b_p(x_0)+d_p(x_0)-z^\top\b a_p(x_0)}
{\|\b b_p(x_0)\|_{\mathcal{Y},*}}.
\end{equation}

\begin{lemma}[\citealt{chen2022data}]\label{lem:distance-formula}
For every $z\in\mathcal{Z}$ and $i\in[N]$, $\Delta_i(z)=\pos{\min_{p\in[P]}m_{ip}(z)}$.
\end{lemma}

Assume from now on that $\mathcal{Z}$ is a nonempty bounded mixed-integer polyhedron and that $c(z)$ is linear. Choose finite constants satisfying
\begin{subequations}\label{eq:valid-bounds}
\begin{align}
M_{ip}
&\geq
\max_{z\in\mathcal{Z}}\pos{-m_{ip}(z)}
&&\forall i\in[N],\ p\in[P] \label{eq:valid-M}\\
\overline\Delta_i
&\geq
\max_{z\in\mathcal{Z}}
\pos{\min_{p\in[P]}m_{ip}(z)}
&&\forall i\in[N]. \label{eq:valid-Deltabar}
\end{align}
\end{subequations}
Boundedness of $\mathcal{Z}$ guarantees that these constants exist. 
For every sample, introduce a continuous variable $\delta_i\geq0$ and a binary variable $\mathfrak{z}_i$. The inequalities
\begin{subequations}\label{eq:distance-hypograph}
\begin{align}
& 0\leq\delta_i \leq\overline\Delta_i(1-\mathfrak{z}_i) &&\forall i\in[N] \label{eq:hypograph-zero}\\
& \delta_i \leq m_{ip}(z)+M_{ip}\mathfrak{z}_i &&\forall i\in[N],\ p\in[P] \label{eq:hypograph-margin}
\end{align}
\end{subequations}
describe the hypograph $0\leq\delta_i\leq\Delta_i(z)$ after projection. Binary $\mathfrak{z}_i=1$ permits a nonpositive sample margin and forces $\delta_i=0$, whereas $\mathfrak{z}_i=0$ enforces all normalized margins above $\delta_i$.

\begin{theorem}\label{thm:exact-mip}
Suppose that Assumption~\ref{ass:transport}, condition~\eqref{eq:strict-radius}, safety system~\eqref{eq:polyhedral-safety}, and bounds~\eqref{eq:valid-bounds} hold. If $c(z)$ is linear and $\mathcal{Z}$ is a bounded mixed-integer polyhedron, then the C-DRCCP model~\eqref{opt:drc3p} is equivalent to
\begin{subequations}\label{opt:exact-mip}
\begin{align}
\min\quad &c(z) \label{eq:mip-objective}\\
\mathrm{s.t.}\quad&
K_0+\epsilon\varepsilon t-\varepsilon\lambda
-\frac{1}{N}\sum_{i\in[N]}v_i
\geq\theta \label{eq:mip-budget}\\
&v_i+\lambda\geq s_i
&&\forall i\in[N] \label{eq:mip-link}\\
&\delta_i+\bar{\kappa}_i\geq t-s_i
&&\forall i\in[N] \label{eq:mip-distance}\\
&0\leq\delta_i\leq\overline\Delta_i(1-\mathfrak{z}_i)
&&\forall i\in[N] \label{eq:mip-zero}\\
&\delta_i\leq m_{ip}(z)+M_{ip}\mathfrak{z}_i
&&\forall i\in[N],\ p\in[P] \label{eq:mip-margin}\\
&s_i\geq-\bar{\kappa}_i,\quad v_i\geq0
&&\forall i\in[N] \label{eq:mip-s-domain}\\
&\lambda\leq\epsilon t,\quad t\geq0,\quad
\lambda\in\R,\quad s\in\R^N \label{eq:mip-dual-domain}\\
&z\in\mathcal{Z},\quad
\delta\in\R_+^N,\quad
\mathfrak{z}\in\{0,1\}^N. \label{eq:mip-domain}
\end{align}
\end{subequations}
\end{theorem}
Thus, formulation~\eqref{opt:exact-mip} adds one binary variable per sample, $3N+2$ continuous variables, and $O(NP)$ linear inequalities to the deterministic model. The contextual geometry enters only through the precomputed constants $K_0$ and $\bar{\kappa}_i$ without additional disjunction.

\begin{remark}\label{remark:lhs}
If the safety system~\eqref{eq:polyhedral-safety} is defined by an individual chance constraint with left-hand-side uncertainty, i.e., $\mathcal{S}(z) \coloneqq \left\{y\in\R^{\dim(y)} \;\middle|\; y^\top\b b(x_0)+d(x_0)-z^\top(\b a(x_0) + \b A(x_0)y)>0\right\}$,  the distances $\Delta_i(z)$ admit explicit expressions, and an exact MIP reformulation can be derived from the constraint system~\eqref{eq:author-certificate}. $\hfill\blacksquare$
\end{remark}

\subsection{Minimum Transportation Cost to the Conditional Risk Boundary}\label{sec:cut-generation}
The LP reformulation~\eqref{eq:mass-flow} evaluates the largest excess failure mass that an adversary can create within a fixed transportation budget. We now reverse this perspective and minimize the transportation cost required to attain the prescribed conditional risk limit. This view provides a direct interpretation of robust feasibility and yields an exact inequality reformulation in terms of the outcome distance variables. For $d\in\R_+^N$, let $\Phi(d)$ denote the optimal value of the LP problem obtained from~\eqref{eq:mass-flow} by replacing $\Delta(z)$ with $d$. Thus,
\begin{equation}\label{eq:artificial-mass-flow}
\Phi(d)=\max\left\{
\mathbf{1}^\top r-\epsilon\mathbf{1}^\top w\;\middle|\;
K_0+\bar{\kappa}^\top w+d^\top r\leq\theta,\ 
\mathbf{1}^\top w\geq\varepsilon,\ 
0\leq r\leq w\leq\frac{1}{N}\mathbf{1}
\right\}.
\end{equation}
To formalize this view, we can define the polytope of mass allocations at the conditional risk boundary as
\begin{equation}\label{eq:critical-allocation-polytope}
\mathcal{C}_{\varepsilon,\epsilon}
\coloneqq
\left\{(w,r)\in\R^{2N}\;\middle|\;
0\leq r\leq w\leq\frac{1}{N}\mathbf{1},
\mathbf{1}^\top w\geq\varepsilon,\ 
\mathbf{1}^\top r=\epsilon\mathbf{1}^\top w
\right\},
\end{equation}
and
\begin{equation}\label{eq:critical-radius-value}
\vartheta(d)
\coloneqq
\min_{(w,r)\in\mathcal{C}_{\varepsilon,\epsilon}}
\left\{K_0+\bar{\kappa}^\top w+d^\top r\right\}.
\end{equation}

For every $(w,r)\in\mathcal{C}_{\varepsilon,\epsilon}$, $\mathbf{1}^\top w$ is the mass assigned to the conditioning neighborhood and $\mathbf{1}^\top r$ is the portion assigned to failure. Because $\mathbf{1}^\top w\geq\varepsilon>0$, the equality $\mathbf{1}^\top r=\epsilon\mathbf{1}^\top w$ places the allocation exactly at the conditional risk boundary. The polytope depends only on $(N,\varepsilon,\epsilon)$ and is nonempty and compact. Thus, $\vartheta(d)$ is the minimum total transportation cost required to reach the risk boundary.

\begin{theorem}\label{thm:critical-radius}
Suppose that condition~\eqref{eq:strict-radius} holds. For every finite $d\in\R_+^N$, we have
\begin{equation}\label{eq:critical-radius-equivalence}
\Phi(d)\leq0
\quad\Longleftrightarrow\quad
\vartheta(d)\geq\theta.
\end{equation}
Thus, a decision $z$ satisfies the robust conditional chance constraint if and only if $\vartheta(\Delta(z))\geq\theta$.
\end{theorem}

Theorem~\ref{thm:critical-radius} compares the transportation budget with the cost of reaching the conditional risk boundary. If $\vartheta(d)<\theta$, the adversary can reach the boundary with budget remaining and create positive excess failure mass. If $\vartheta(d)\geq\theta$, it cannot exceed the prescribed conditional risk limit. Because $\vartheta$ is the minimum of affine functions over $\mathcal{C}_{\varepsilon,\epsilon}$, the latter condition is equivalent to the following exact inequality reformulation.

\begin{corollary}\label{cor:critical-radius-cuts}
In model~\eqref{opt:exact-mip}, the variables $(t,s,v,\lambda)$ and constraints~\eqref{eq:mip-budget}--\eqref{eq:mip-distance} and~\eqref{eq:mip-s-domain}--\eqref{eq:mip-dual-domain} can be replaced by
\begin{equation}\label{eq:critical-radius-cuts}
K_0+\bar{\kappa}^\top\bar{w}+\bar{r}^\top\delta\geq\theta
\qquad
\forall(\bar{w},\bar{r})\in\mathcal{C}_{\varepsilon,\epsilon}.
\end{equation}
The resulting MIP formulation has the same projection onto $z$ as formulation~\eqref{opt:exact-mip}. Exact separation at a candidate $\bar{\delta}$ requires solving the LP problem~\eqref{eq:critical-radius-value}. If $\vartheta(\bar{\delta})<\theta$, an optimal solution $(\bar{w},\bar{r})$ identifies a violated inequality. The LP problem has an optimal solution for which all but at most two blocks $(\bar{w}_i,\bar{r}_i)$ belong to $\{(0,0),(1/N,0),(1/N,1/N)\}$.
\end{corollary}

Each $(\bar{w},\bar{r})\in\mathcal{C}_{\varepsilon,\epsilon}$ defines a halfspace in the distance variables. The failure allocation $\bar{r}$ gives its coefficients, while $K_0+\bar\kappa^\top\bar{w}$ gives its constant term. Because $\mathcal{C}_{\varepsilon,\epsilon}$ is a nonempty compact polytope, it suffices to impose the inequalities indexed by its extreme points. Moreover, Corollary~\ref{cor:critical-radius-cuts} shows that a separating inequality can be obtained from an allocation in which all but at most two empirical sources are assigned entirely outside the conditioning neighborhood, inside the neighborhood and safe, or inside the neighborhood and failing. Thus, the representation separates the contextual cost of forming the conditioning event from the outcome cost of reaching failure. We also use Example~\ref{example:cg} to demonstrate how different context allocations can produce inequalities of different strength even when the failure vector remains unchanged.

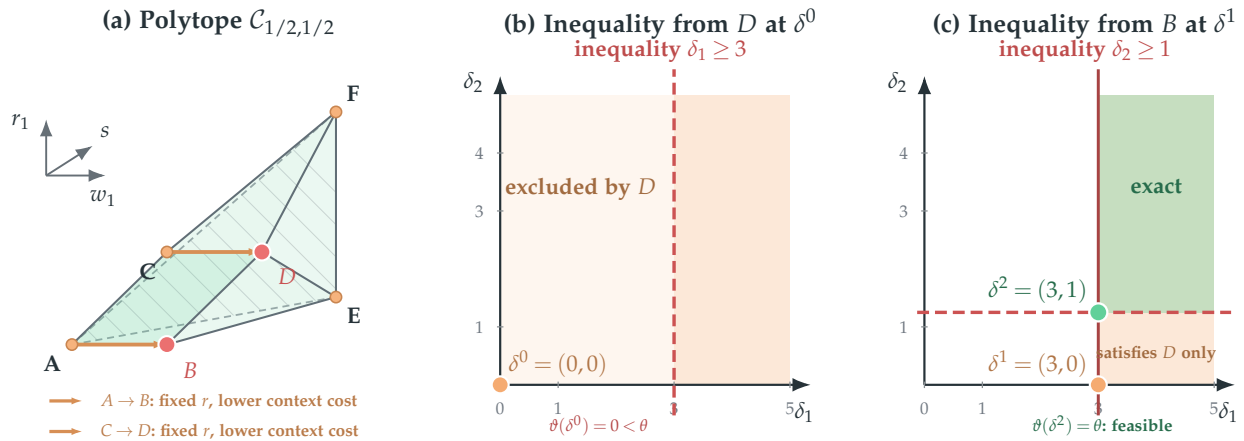
\begin{figure}[htbp]
  \centering
  \caption{Geometry of the exact inequality reformulation in the example with two samples in Example~\ref{example:cg}.}
  \label{fig:critical-radius-separation-example}
  \definecolor{cgGreen}{HTML}{47C88A}
  \definecolor{cgOrange}{HTML}{F4A261}
  \definecolor{cgRed}{HTML}{E85D5D}
  \definecolor{cgEdge}{HTML}{34454F}
  \resizebox{\linewidth}{!}{%
  \begin{tikzpicture}[
    x=1cm,y=1cm,font=\footnotesize,>=Latex,
    visible edge/.style={draw=cgEdge!78,line width=0.72pt},
    hidden edge/.style={draw=cgEdge!48,line width=0.64pt,
      dash pattern=on 2.5pt off 1.7pt}
  ]
    \node[ctx title] at (2.725,5.27)
      {(a) Polytope $\mathcal{C}_{1/2,1/2}$};

    \coordinate (O) at (0.58,3.38);
    \draw[->,line width=0.68pt,draw=ctxgray!82!black]
      (O) -- (1.32,3.38);
    \draw[->,line width=0.68pt,draw=ctxgray!82!black]
      (O) -- (0.58,4.08);
    \draw[->,line width=0.68pt,draw=ctxgray!82!black]
      (O) -- (1.16,3.75);
    \node[anchor=north,font=\scriptsize,text=ctxink] at (1.29,3.33) {$w_1$};
    \node[anchor=east,font=\scriptsize,text=ctxink] at (0.52,4.05) {$r_1$};
    \node[anchor=south west,font=\scriptsize,text=ctxink] at (1.12,3.72) {$s$};

    \coordinate (A) at (0.90,1.28);
    \coordinate (B) at (2.08,1.28);
    \coordinate (C) at (2.08,2.43);
    \coordinate (D) at (3.26,2.43);
    \coordinate (E) at (4.18,1.87);
    \coordinate (F) at (4.18,4.17);

    \fill[fill=cgGreen,fill opacity=0.035] (A) -- (B) -- (E) -- cycle;
    \fill[fill=cgGreen,fill opacity=0.035] (A) -- (C) -- (F) -- cycle;
    \fill[fill=cgGreen,fill opacity=0.045] (A) -- (E) -- (F) -- cycle;
    \begin{scope}
      \clip (A) -- (E) -- (F) -- cycle;
      \foreach \t in {0.4,0.8,...,7.2} {
        \draw[draw=cgEdge!22,line width=0.22pt]
          (\t,0.55) -- ++(-3.8,3.8);
      }
    \end{scope}

    \draw[hidden edge] (A) -- (E);
    \draw[hidden edge] (A) -- (F);

    \fill[fill=cgGreen,fill opacity=0.085] (B) -- (D) -- (E) -- cycle;
    \fill[fill=cgGreen,fill opacity=0.085] (C) -- (D) -- (F) -- cycle;
    \fill[fill=cgGreen,fill opacity=0.065] (D) -- (E) -- (F) -- cycle;
    \fill[fill=cgGreen,fill opacity=0.20]  (A) -- (B) -- (D) -- (C) -- cycle;
    \foreach \u/\v in {A/B,A/C,B/D,B/E,C/D,C/F,D/E,D/F,E/F} {
      \draw[visible edge] (\u) -- (\v);
    }

    \draw[-{Latex[length=1.65mm,width=1.0mm]},draw=cgOrange!88!black,
      line width=1.35pt] (A) -- (B);
    \draw[-{Latex[length=1.65mm,width=1.0mm]},draw=cgOrange!88!black,
      line width=1.35pt] (C) -- (D);

    \foreach \p/\lab/\pos in {A/A/below left,C/C/below left,E/E/below right,F/F/above right} {
      \filldraw[fill=cgOrange!82,draw=cgOrange!82!black,line width=0.55pt]
        (\p) circle (2.0pt);
      \node[font=\scriptsize\bfseries,text=ctxink,\pos] at (\p) {\lab};
    }
    \filldraw[fill=cgRed!88,draw=white,line width=0.75pt]
      (B) circle (3.0pt);
    \node[anchor=north west,font=\scriptsize\bfseries,text=cgRed!86!black]
      at (2.14,1.20) {$B$};
    \filldraw[fill=cgRed!88,draw=white,line width=0.75pt]
      (D) circle (3.0pt);
    \node[anchor=north west,font=\scriptsize\bfseries,text=cgRed!86!black]
      at (3.32,2.36) {$D$};


    \draw[-{Latex[length=1.45mm,width=0.9mm]},draw=cgOrange!88!black,
      line width=1.15pt] (0.62,0.58) -- (1.04,0.58);
    \node[anchor=west,font=\tiny\bfseries,text=cgOrange!76!black]
      at (1.12,0.58) {$A\to B$: fixed $r$, lower context cost};
    \draw[-{Latex[length=1.45mm,width=0.9mm]},draw=cgOrange!88!black,
      line width=1.15pt] (0.62,0.20) -- (1.04,0.20);
    \node[anchor=west,font=\tiny\bfseries,text=cgOrange!76!black]
      at (1.12,0.20) {$C\to D$: fixed $r$, lower context cost};

    \node[ctx title] at (8.21,5.27) {(b) Inequality from $D$ at $\delta^0$};
    \node[font=\scriptsize\bfseries,text=cgRed!82!black] at (8.21,4.93)
      {inequality $\delta_1\geq3$};

    \coordinate (Ob) at (6.22,0.78);
    \def\sb{0.72}
    \fill[cgOrange!10] (6.22,0.78) rectangle (9.82,4.38);
    \fill[cgOrange!25] (8.38,0.78) rectangle (9.82,4.38);

    \draw[->,line width=0.78pt,draw=ctxink] (Ob) -- (10.08,0.78);
    \draw[->,line width=0.78pt,draw=ctxink] (Ob) -- (6.22,4.62);
    \node[anchor=north,font=\scriptsize,text=ctxink] at (10.00,0.71) {$\delta_1$};
    \node[anchor=east,font=\scriptsize,text=ctxink] at (6.15,4.54) {$\delta_2$};
    \foreach \x/\lab in {0/0,1/1,3/3,5/5} {
      \draw[draw=ctxink!55,line width=0.45pt]
        ({6.22+\x*\sb},0.74) -- ({6.22+\x*\sb},0.82);
      \node[anchor=north,font=\tiny,text=ctxgray!90!black]
        at ({6.22+\x*\sb},0.70) {\lab};
    }
    \foreach \y/\lab in {1/1,3/3,4/4} {
      \draw[draw=ctxink!55,line width=0.45pt]
        (6.18,{0.78+\y*\sb}) -- (6.26,{0.78+\y*\sb});
      \node[anchor=east,font=\tiny,text=ctxgray!90!black]
        at (6.15,{0.78+\y*\sb}) {\lab};
    }

    \draw[draw=cgRed!88!black,line width=1.20pt,
      dash pattern=on 4.2pt off 2.1pt] (8.38,0.4) -- (8.38,4.7);
    \node[align=center,font=\scriptsize\bfseries,text=cgOrange!66!black]
      at (7.22,3.22) {excluded by $D$};

    \filldraw[fill=cgOrange!90,draw=white,line width=0.75pt] (Ob) circle (3.0pt);
    \node[anchor=south west,font=\scriptsize\bfseries,text=cgOrange!72!black,
      fill opacity=0.88,text opacity=1,inner sep=1.0pt]
      at (6.28,0.84) {$\delta^0=(0,0)$};
    \node[font=\tiny\bfseries,text=cgRed!78!black] at (7.44,0.28) {$\vartheta(\delta^0)=0<\theta$};

    \node[ctx title] at (13.485,5.27) {(c) Inequality from $B$ at $\delta^1$};
    \node[font=\scriptsize\bfseries,text=cgRed!82!black] at (13.485,4.93)
      {inequality $\delta_2\geq1$};

    \coordinate (Oc) at (11.49,0.78);
    \def\sc{0.72}
    \fill[cgOrange!25] (13.65,0.78) rectangle (15.09,4.38);
    \fill[fill=cgGreen,fill opacity=0.30] (13.65,1.68) rectangle (15.09,4.38);

    \draw[->,line width=0.78pt,draw=ctxink] (Oc) -- (15.35,0.78);
    \draw[->,line width=0.78pt,draw=ctxink] (Oc) -- (11.49,4.62);
    \node[anchor=north,font=\scriptsize,text=ctxink] at (15.27,0.71) {$\delta_1$};
    \node[anchor=east,font=\scriptsize,text=ctxink] at (11.42,4.54) {$\delta_2$};
    \foreach \x/\lab in {0/0,1/1,3/3,5/5} {
      \draw[draw=ctxink!55,line width=0.45pt]
        ({11.49+\x*\sc},0.74) -- ({11.49+\x*\sc},0.82);
      \node[anchor=north,font=\tiny,text=ctxgray!90!black]
        at ({11.49+\x*\sc},0.70) {\lab};
    }
    \foreach \y/\lab in {1/1,3/3,4/4} {
      \draw[draw=ctxink!55,line width=0.45pt]
        (11.45,{0.78+\y*\sc}) -- (11.53,{0.78+\y*\sc});
      \node[anchor=east,font=\tiny,text=ctxgray!90!black]
        at (11.42,{0.78+\y*\sc}) {\lab};
    }

    \draw[draw=cgRed!72!black,line width=1.05pt] (13.65,0.4) -- (13.65,4.7);
    \draw[draw=cgRed!88!black,line width=1.20pt,
      dash pattern=on 4.2pt off 2.1pt]
      (11.4,1.68) -- (15.3,1.68);

    \node[align=center,font=\tiny\bfseries,text=cgOrange!66!black]
      at (14.37,1.16) {satisfies $D$ only};
    \node[align=center,font=\scriptsize\bfseries,text=cgGreen!54!black]
      at (14.37,3.25) {exact};

    \filldraw[fill=cgOrange!90,draw=white,line width=0.75pt]
      (13.65,0.78) circle (3.0pt);
    \node[anchor=south east,font=\scriptsize\bfseries,text=cgOrange!72!black,
      fill opacity=0.88,text opacity=1,inner sep=1.0pt]
      at (13.56,0.84) {$\delta^1=(3,0)$};

    \filldraw[fill=cgGreen!86,draw=white,line width=0.80pt] (13.65,1.68) circle (3.2pt);
    \node[anchor=south east,font=\scriptsize\bfseries,text=cgGreen!54!black,
      fill opacity=0.88,text opacity=1,inner sep=1.0pt]
      at (13.56,1.76) {$\delta^2=(3,1)$};
    \node[font=\tiny\bfseries,text=cgGreen!54!black]
      at (13.70,0.28) {$\vartheta(\delta^2)=\theta$: feasible};
  \end{tikzpicture}%
  }
\end{figure}

\begin{example}\label{example:cg}
Consider an example with two samples. Let $N=2$, $\varepsilon=\epsilon=1/2$, and let the conditioning neighborhood be a norm ball with radius $\gamma=1/2$. If the distances from the two sample contexts to $x_0$ are $0$ and $2$, respectively, then the expressions in Section~\ref{sec:geometry} give $K_0=1/4$ and $\bar\kappa=(-1/2,3/2)$. Set $\theta=3/4$. The minimum contextual transportation cost is zero, so condition~\eqref{eq:strict-radius} holds. Writing $s=\mathbf{1}^\top w$, the equalities $w_2=s-w_1$ and $r_2=s/2-r_1$ show that $(w_1,r_1,s)$ provides a unique coordinate system for the complete polytope $\mathcal{C}_{1/2,1/2}$ with vertices
\[\begin{array}{c|cccccc}
 &A&B&C&D&E&F\\ \hline
(w_1,r_1,s)
 &(0,0,\tfrac12)&(\tfrac14,0,\tfrac12)
 &(\tfrac14,\tfrac14,\tfrac12)&(\tfrac12,\tfrac14,\tfrac12)
 &(\tfrac12,0,1)&(\tfrac12,\tfrac12,1).
\end{array}\]

At any candidate $\bar{\delta}\in\mathbb R_+^2$, the LP problem~\eqref{eq:critical-radius-value} computes $\vartheta(\bar{\delta})$ and thereby checks the entire inequality family. If $\vartheta(\bar{\delta})\geq\theta$, then every inequality indexed by $\mathcal{C}_{1/2,1/2}$ is satisfied. Otherwise, an optimal allocation $(\bar{w},\bar{r})$ identifies a most violated inequality, $K_0+\bar\kappa^\top\bar{w}+\bar{r}^\top\delta\geq\theta$.

Figure~\ref{fig:critical-radius-separation-example} evaluates the inequality family at two distance vectors. At $\delta^0=(0,0)$, vertex $D:(w,r)=((1/2,0),(1/4,0))$ is the unique optimal solution and gives the violated inequality $\delta_1\geq3$. At $\delta^1=(3,0)$, vertex $B:(w,r)=((1/4,1/4),(0,1/4))$ is optimal and gives the violated inequality $\delta_2\geq1$. At $\delta^2=(3,1)$, the value of $\vartheta$ equals $\theta$. Vertices $A$ and $B$ have the same failure allocation $r=(0,1/4)$, but their contextual transportation costs $K_0+\bar\kappa^\top w$ are $1$ and $1/2$, respectively. Thus, the inequality associated with $A$ is redundant over $\mathbb R_+^2$, whereas $B$ gives $\delta_2\geq1$. Similarly, $C$ and $D$ have the same failure allocation $r=(1/4,0)$, but their contextual transportation costs are $1/2$ and $0$. Hence, the inequality $\delta_1\geq3$ associated with $D$ dominates the inequality $\delta_1\geq1$ associated with $C$. The inequalities associated with $E$ and $F$ are also redundant, and the exact robust region is $\left\{\delta\in\mathbb R_+^2\;\middle|\;\delta_1\geq3,\ \delta_2\geq1\right\}$.
$\hfill\blacksquare$
\end{example}

For the coefficient calculations in Section~\ref{sec:radius-adjusted-bounds}, we use the following dual representation of~\eqref{eq:critical-radius-value}:
\begin{equation}\label{eq:critical-radius-dual}
\vartheta(d)=
\max_{\alpha,\mu,\nu,\tau\geq0}
\left\{K_0+\varepsilon\alpha-\frac{1}{N}\sum_{i \in [N]}\nu_i\;\middle|\;
\alpha+\mu_i-\nu_i-\epsilon\tau\leq\bar{\kappa}_i,\ 
\tau-\mu_i\leq d_i\quad \forall i \in [N]
\right\}.
\end{equation}
Although the dual variable $\tau$ of the equality constraint $\mathbf{1}^\top r = \epsilon\mathbf{1}^\top w$ should be free, any negative value of $\tau$ can be increased to zero without changing the objective or violating a constraint. Therefore, we impose $\tau\geq0$ without loss of generality.
\begin{remark}
Using Theorem~\ref{thm:critical-radius}, Appendix~\ref{appendix:model-relations} characterizes the relation between the C-DRCCP model and the unconditional DRCCP model of \citet{chen2022data}. In general, the C-DRCCP model cannot be obtained by changing the Wasserstein radius alone.$\hfill\blacksquare$
\end{remark}

\subsection{Statistical Guarantees}\label{sec:statistical}
The reformulations are exact for any empirical sample and any parameters $(\theta,\varepsilon,\gamma)$ that satisfy our assumptions. Theorem~\ref{thm:finite-sample} shows that, when the Wasserstein radius and the minimum probability of the conditioning event are calibrated at their respective confidence levels, any decision selected from the data and feasible for the C-DRCCP model~\eqref{opt:drc3p} controls the true conditional violation probability.

\begin{theorem}\label{thm:finite-sample}
Let $\mathbb{P}^\star$ be the true distribution, let $\widehat{\mathbb{P}}_N$ be formed from $N$ independent and identically distributed observations, and fix $\beta_W,\beta_{\mathcal{N}}\in[0,1)$ with $\beta_W+\beta_{\mathcal{N}}<1$. The Borel neighborhood $\mathcal{N}_N$, the radius $\theta_N\geq0$, the minimum probability $\varepsilon_N>0$, and decision $z_N$ may all depend on the data. Assume that the resulting random quantities and the events below are measurable, that $z_N$ is feasible for the realized C-DRCCP model almost surely, and that
\begin{subequations}
\begin{align}
(\mathbb{P}^\star)^N\!\left[
W(\mathbb{P}^\star,\widehat{\mathbb{P}}_N)\leq\theta_N
\right]&\geq1-\beta_W, \label{eq:W-coverage}\\
(\mathbb{P}^\star)^N\!\left[
\mathbb{P}_X^\star(\mathcal{N}_N)\geq\varepsilon_N
\right]&\geq1-\beta_{\mathcal{N}}. \label{eq:mass-coverage}
\end{align}
Then 
\begin{align}
(\mathbb{P}^\star)^N\Big[
&\mathbb{P}_X^\star(\mathcal{N}_N)\geq\varepsilon_N,\;\mathbb{P}^\star\left((X,Y)\in
\mathcal{N}_N\times\mathcal{S}^c(z_N)\right)
\leq\epsilon\mathbb{P}_X^\star(\mathcal{N}_N)
\Big]
\geq1-\beta_W-\beta_{\mathcal{N}},
\label{eq:finite-sample-guarantee}
\end{align}
\end{subequations}
that is, the conditional violation probability is well-defined and at most $\epsilon$.
\end{theorem}

The coverage event~\eqref{eq:W-coverage} can be calibrated through concentration inequalities under explicit tail and moment assumptions, as in~\citet{mohajerin2018data}, or through a justified resampling or validation procedure~\citep{gao2023distributionally}. For a neighborhood $\mathcal{N}$ fixed independently of the data and $\beta_{\mathcal{N}}\in(0,1)$, let $\widehat{p}_N=N^{-1}\sum_{i\in[N]}\mathbbm 1\{\hat{x}_i\in\mathcal{N}\}$. A one-sided Hoeffding inequality implies that, with probability at least $1-\beta_{\mathcal{N}}$, we have $\mathbb{P}_X^\star(\mathcal{N})\geq\widehat{p}_N-\sqrt{\frac{\log(1/\beta_{\mathcal{N}})}{2N}}$, where the right-hand side may be used as $\varepsilon_N$. 

Neighborhood conditioning avoids the nonuniqueness in Proposition~\ref{prop:null-kernel} and allows the ambiguity set to contain continuous joint distributions. When the conditional failure probability varies smoothly near $x_0$, its average over the neighborhood approximates its value at $x_0$. Proposition~\ref{prop:pointwise} quantifies this approximation and explains how the neighborhood radius $\gamma$ controls the approximation error.

\begin{proposition}\label{prop:pointwise}
For each decision $z\in\mathcal{Z}$, let $p_z:\mathcal{N}_\gamma(x_0)\to[0,1]$ be a fixed measurable version of $\mathbb{P}^\star\left(Y\notin\mathcal{S}(z)\mid X=x\right)$. Suppose that there exists $L\in[0,\infty)$ such that $|p_z(x)-p_z(x_0)|\leq Ld_{\mathcal{X}}(x,x_0)$ for $x\in\mathcal{N}_\gamma(x_0)$ and $z\in\mathcal{Z}$. If $\mathbb{P}_X^\star(\mathcal{N}_\gamma(x_0))>0$, then every $z\in\mathcal{Z}$ satisfies
\begin{equation}\label{eq:pointwise-bound}
p_z(x_0)
\leq
\mathbb{P}^\star\left(Y\notin\mathcal{S}(z)\mid
X\in\mathcal{N}_\gamma(x_0)\right)+L\gamma.
\end{equation}
Thus, if $L\gamma<\epsilon$ and the model is enforced with the tighter risk limit $\epsilon-L\gamma$, then ambiguity-set coverage implies that the pointwise violation probability at $x_0$ is at most $\epsilon$.
\end{proposition}
When $x_0$ is a null point, the joint distribution determines $p_z(x)$ only $\mathbb{P}_X^\star$-almost everywhere and leaves $p_z(x_0)$ unrestricted. The Lipschitz condition supplies the missing local structure by restricting attention to versions for which $p_z(x_0)$ differs from the conditional failure probability at a nearby context $x$ by at most $L d_{\mathcal{X}}(x,x_0)$. Consequently, control of the conditional failure probability over $\mathcal{N}_\gamma(x_0)$ can be transferred to pointwise control at $x_0$, with the explicit localization penalty $L\gamma$.

\section{Valid Inequalities}\label{sec:strengthening}

The exact reformulation~\eqref{opt:exact-mip} represents the entire conditional transportation problem, but its LP relaxation does not use all of the structure carried by the binary variables. This section develops complementary remedies. The transportation budget remaining after the context allocation yields strengthened quantile cuts and binary knapsack and rank inequalities.

Define the polytope of contextual mass allocations
\begin{equation}\label{eq:context-allocation-polytope}
\mathcal{W}_\theta
\coloneqq
\left\{w\in\R^N\;\middle|\;
0\leq w_i\leq\frac{1}{N},
\mathbf{1}^\top w\geq\varepsilon,\ 
K_0+\bar{\kappa}^\top w\leq\theta
\right\},
\end{equation}
and let $B(w)\coloneqq\theta-K_0-\bar{\kappa}^\top w$ denote the transportation budget left after mass is moved according to $w$. Condition~\eqref{eq:strict-radius} implies that $\mathcal{W}_\theta$ is nonempty and contains a point with $B(w)>0$.

For every sample $i$ and safety constraint $p$, decompose the normalized margin as
\begin{equation}\label{eq:margin-decomposition}
c_{ip}
\coloneqq
\frac{\hat y_i^\top\b b_p(x_0)}{\|\b b_p(x_0)\|_{\mathcal{Y},*}},
\quad
\beta_p(z)
\coloneqq
\frac{d_p(x_0)-z^\top\b a_p(x_0)}{\|\b b_p(x_0)\|_{\mathcal{Y},*}},
\quad \text{and}\quad
m_{ip}(z)=c_{ip}+\beta_p(z).
\end{equation}
Here, $c_{ip}$ is the outcome component of sample $i$, whereas $\beta_p(z)$ is the component shared across all samples in safety constraint $p$. Thus, a larger value of $\beta_p(z)$ moves all samples toward the safe side of the $p$-th boundary and increases their distance to its failure region.

\subsection{Probability-Cut Closure}\label{sec:conditional-probability-set}

Every nonzero $w\in\mathcal{W}_\theta$ induces conditional scenario probabilities by normalization, $p=w/(\mathbf{1}^\top w)$. The set of all such probabilities is $\mathcal{D}_\theta^{\rm ctx}\coloneqq\left\{p\in\R_+^N\;\middle|\;p=\frac{w}{\mathbf{1}^\top w}\text{ for some }w\in\mathcal{W}_\theta\right\}$.
The following lift makes its polyhedral nature explicit and places the probability cuts in the existing finite-support theory.

\begin{proposition}\label{prop:induced-probability-set}
Set $\mathcal{D}_\theta^{\rm ctx}$ is the projection onto $p$ of
\begin{equation}\label{eq:induced-probability-lift}
\mathbf{1}^\top p=1,\qquad
0\leq p_i\leq\frac qN\ \ \forall i,\qquad
1\leq q\leq\frac1\varepsilon,\qquad
\bar{\kappa}^\top p\leq(\theta-K_0)q.
\end{equation}
Moreover, every feasible $\mathfrak{z}$ of problem~\eqref{opt:exact-mip} satisfies
\begin{equation}\label{eq:probability-closure}
p^\top\mathfrak{z}\leq\epsilon
\quad\forall p\in\mathcal{D}_\theta^{\rm ctx},
\end{equation}
or, equivalently, $w^\top\mathfrak{z}\leq\epsilon\mathbf{1}^\top w$ for every $w\in\mathcal{W}_\theta$. The full family~\eqref{eq:probability-closure} is equivalent to the existence of $\tau^{\rm P},\eta^{\rm P}\geq0$ and $\nu^{\rm P}\in\R_+^N$ satisfying
\begin{subequations}\label{eq:compact-probability-closure}
\begin{align}
& (\theta-K_0)\tau^{\rm P}-\varepsilon\eta^{\rm P}+\frac{1}{N}\sum_{i\in[N]}\nu_i^{\rm P}\leq0 \label{eq:compact-probability-budget}\\
& \bar{\kappa}_i\tau^{\rm P}-\eta^{\rm P}+\nu_i^{\rm P} \geq\mathfrak{z}_i-\epsilon &&\forall i\in[N].\label{eq:compact-probability-link}
\end{align}
\end{subequations}
\end{proposition}

The probability cut~\eqref{eq:probability-closure} extends the probabilistic cover inequality~\citep{song2014packing}. The contextual content is the set~\eqref{eq:induced-probability-lift} and, as shown below, the additional information carried by $B(w)$ that disappears after normalization. For comparison with the strengthened quantile cut derived below, define the quantile
\begin{equation}\label{eq:probability-quantile}
q_p
\coloneqq
\min\left\{q\in\left\{c_{ip}\;\middle|\;i\in[N]\right\}\;\middle|\;
\rho_\theta\!\left(\left\{i\;\middle|\;c_{ip}\leq q\right\}\right)>0
\right\},
\end{equation}
where $w(A)\coloneqq\sum_{i\in A}w_i$ and $\rho_\theta(A)\coloneqq\max_{w\in\mathcal{W}_\theta}\left\{w(A)-\epsilon\mathbf{1}^\top w\right\}$. It is well-defined because the set in~\eqref{eq:probability-quantile} contains $\max_i c_{ip}$.
Indeed, $\rho_\theta([N])=(1-\epsilon)\max_{w\in\mathcal{W}_\theta}\mathbf{1}^\top w>0$.
\begin{subequations}
\begin{proposition}\label{prop:probability-quantile-row-bound}
For every $p\in[P]$, every feasible solution to model~\eqref{opt:exact-mip} satisfies the quantile cut
\begin{equation}\label{eq:probability-quantile-row-bound}
\beta_p(z)\geq-q_p.
\end{equation}
Thus, after adding~\eqref{eq:probability-quantile-row-bound}, every coefficient $M_{ip}$ in~\eqref{eq:mip-margin} can be replaced, without changing the feasible region of the resulting mixed-integer formulation~\eqref{opt:exact-mip}, by
\begin{equation}\label{eq:probability-quantile-M}
M_{ip}^{\rm Q}
\coloneqq
\min\left\{M_{ip},\pos{q_p-c_{ip}}\right\}.
\end{equation}
\end{proposition}
\end{subequations}
\subsection{Strengthened Quantile Cuts and Coefficient Tightening}\label{sec:radius-adjusted-bounds}

The probability-quantile bound~\eqref{eq:probability-quantile-row-bound} uses only whether selected samples can have zero distance. A sharper threshold also incurs the transportation costs needed to create additional failure mass. For $\beta\in\R$, define
\begin{equation}\label{eq:radius-row-value}
d_i^p(\beta)\coloneqq\pos{c_{ip}+\beta},
\quad
\Psi_p(\beta)\coloneqq\vartheta(d^p(\beta)),
\quad \text{and} \quad
\underline{\beta}_p\coloneqq
\min\left\{\beta\in\R\;\middle|\;\Psi_p(\beta)\geq\theta\right\}.
\end{equation}
Interpreting $\beta$ as a candidate value of the common decision component $\beta_p(z)$, $d_i^p(\beta)$ is the unit transportation cost of moving mass from sample $i$ into failure through safety constraint $p$. The value $\Psi_p(\beta)$ is the least total transportation cost required for the adversary to reach the conditional risk boundary $\epsilon$, including both contextual and outcome transportation. Accordingly, $\underline{\beta}_p$ is the critical decision margin at which the least cost reaches the budget $\theta$.
\begin{subequations}
\begin{theorem}\label{thm:radius-adjusted-threshold}
For every $p\in[P]$, the threshold $\underline{\beta}_p$ in~\eqref{eq:radius-row-value} exists and is finite. Every feasible decision for the robust model satisfies the strengthened quantile cut
\begin{equation}\label{eq:radius-row-bound}
\beta_p(z)\geq\underline{\beta}_p.
\end{equation}
If the safety system contains only safety constraint $p$, inequality~\eqref{eq:radius-row-bound} is also sufficient for the robust chance constraint and therefore gives its exact scalar threshold.
Moreover,
\begin{equation}\label{eq:radius-dominates-probability}
\underline{\beta}_p>-q_p.
\end{equation}
Consequently, after adding~\eqref{eq:radius-row-bound}, every coefficient $M_{ip}$ in~\eqref{eq:mip-margin} can be replaced without changing the feasible region of model~\eqref{opt:exact-mip} by
\begin{equation}\label{eq:radius-adjusted-M}
\widetilde M_{ip}
\coloneqq
\min\left\{M_{ip},\pos{-c_{ip}-\underline{\beta}_p}\right\}.
\end{equation}
These coefficients are componentwise no larger than the probability-quantile coefficients~\eqref{eq:probability-quantile-M}.
\end{theorem}
\end{subequations}
For $p\in[P]$ and $w\in\mathcal{W}_\theta$ with $B(w)>0$, define
\begin{subequations}\label{eq:fixed-allocation-threshold-definition}
\begin{align}
g_{w,p}(\beta)&\coloneqq\min\left\{\sum_{i \in [N]} r_i d_i^p(\beta)\;\middle|\;0\leq r\leq w,\ \mathbf{1}^\top r=\epsilon\mathbf{1}^\top w\right\},\\
\beta_p^\star(w)&\coloneqq\min\left\{\beta\in\R\;\middle|\;g_{w,p}(\beta)\geq B(w)\right\}.
\end{align}
\end{subequations}
The function $g_{w,p}(\beta)$ is the least cost of transporting outcomes to reach the conditional risk boundary under allocation $w$. Indeed, the conditional failure probability is $(\mathbf{1}^\top r)/(\mathbf{1}^\top w)$, so reaching the boundary $\epsilon$ requires $\mathbf{1}^\top r=\epsilon\mathbf{1}^\top w$. Because the costs of transporting outcomes are nonnegative, it suffices to allocate exactly this amount of failure mass. If $g_{w,p}(\beta)<B(w)$, the adversary reaches the boundary with budget left. Moreover, because $\epsilon<1$, the boundary allocation leaves positive unused capacity, and the adversary can use any strict budget slack to add a sufficiently small amount of failure mass and cross the boundary. Thus, $\beta_p^\star(w)$ is the common decision margin at which the cheapest plan that reaches the boundary for the allocation $w$ \emph{exhausts} the residual budget.

For a fixed $\beta$, assigning one unit of mass from sample $i$ to failure costs $d_i^p(\beta)$, and sample $i$ has capacity $w_i$. The common shift and positive-part operation preserve the ordering of the $c_{ip}$. Hence, if a cheaper sample has unused capacity while a more expensive sample carries positive failure mass, shifting mass from the latter to the former preserves $\mathbf{1}^\top r=\epsilon\mathbf{1}^\top w$ and weakly decreases cost. Repeating this exchange fills the cheapest capacities first: an initial set of samples is full, at most one subsequent sample is partial, and all remaining samples are unused. Taking the supremum of the resulting fixed-$w$ thresholds gives $\underline{\beta}_p$. By contrast, the probability quantile records only which distances are zero. At $\beta=-q_p$, the minimum critical cost remains strictly below $\theta$, so the quantile threshold leaves transport budget unused. Moving to the first crossing $\underline{\beta}_p$ explains both the strict dominance in~\eqref{eq:radius-dominates-probability} and the coefficients in~\eqref{eq:radius-adjusted-M}. This improvement can be substantial even with one sample, as shown in Example~\ref{example:quantilecut}. 
\begin{example}\label{example:quantilecut}
Let $N=1$, $K_0=\bar{\kappa}_1=c_{1p}=0$, $\varepsilon=\epsilon=1/2$, and $\theta=1$. The quantile gives the bound $-q_p=0$. In contrast, $\Psi_p(\beta)=\tfrac14\pos{\beta}$, so $\underline{\beta}_p=4$. The stronger bound captures the transportation cost required to place the failure mass at the risk boundary. $\hfill\blacksquare$
\end{example}

The threshold $\underline{\beta}_p$ is computable exactly by linear programming. Sort the at most $N$ breakpoints $-c_{ip}$. On each of the resulting $N+1$ intervals, $d_i^p(\beta)$ is either zero or $c_{ip}+\beta$. By the reformulation in~\eqref{eq:critical-radius-dual}, the smallest feasible $\beta$ on that interval is the value of a linear program with variables $(\beta,\alpha,\mu,\nu,\tau)$, the interval bounds, and
\begin{equation}\label{eq:radius-threshold-preprocessing}
K_0+\varepsilon\alpha-\frac{1}{N}\sum_{i \in [N]}\nu_i\geq\theta,
\qquad
\alpha+\mu_i-\nu_i-\epsilon\tau\leq\bar{\kappa}_i,
\qquad
\tau-\mu_i\leq
\begin{cases}
0&c_{ip}+\beta\leq0\\
c_{ip}+\beta&c_{ip}+\beta\geq0
\end{cases},
\end{equation}
with $\alpha,\mu,\nu,\tau\geq0$. The minimum over the feasible interval linear programs is $\underline{\beta}_p$. Alternatively, monotonicity permits bisection with warm-started solves of~\eqref{eq:critical-radius-value}. One can also compute margin cuts from fixed context allocations. The next proposition gives their exact relation to $\underline\beta_p$.

\begin{proposition}\label{prop:fixed-allocation-threshold}
Fix $p\in[P]$. For every $w\in\mathcal{W}_\theta$ with $B(w)>0$, the threshold $\beta_p^\star(w)$ is finite, and every feasible decision for the robust model satisfies $\beta_p(z)\geq\beta_p^\star(w)$. Moreover,
\begin{equation}\label{eq:fixed-allocation-supremum}
\underline\beta_p = \sup\left\{\beta_p^\star(w)\;\middle|\; w\in\mathcal{W}_\theta,\ B(w)>0
\right\}.
\end{equation}
\end{proposition}

The strengthened quantile cuts also supply the construction of \emph{mixing inequalities}. Define $u_p(z)=\beta_p(z)-\underline{\beta}_p$ and $h_{ip}=\pos{-c_{ip}-\underline{\beta}_p}$. After adding $u_p(z)\geq0$ and using~\eqref{eq:radius-adjusted-M}, the LP relaxation implies the base inequalities $u_p(z)\geq h_{ip}(1-\mathfrak{z}_i)$. Indeed, when $h_{ip}>0$, constraints $\delta_i\geq0$ and~\eqref{eq:mip-margin} give $u_p(z)\geq h_{ip}-\widetilde M_{ip}\mathfrak z_i\geq h_{ip}(1-\mathfrak z_i)$ because $\widetilde M_{ip}\leq h_{ip}$. When $h_{ip}=0$, the inequality follows from $u_p(z)\geq0$. Corollary~\ref{cor:radius-adjusted-mixing} records the resulting mixing inequalities.

\begin{corollary}\label{cor:radius-adjusted-mixing}
Fix $p\in[P]$ and a sequence of distinct indices $j_1,\ldots,j_\ell$ such that $h_{j_1p}\geq\cdots\geq h_{j_\ell p}>0$, and set $h_{j_{\ell+1}p}=0$. Every feasible solution to model~\eqref{opt:exact-mip} satisfies
\begin{equation}\label{eq:radius-adjusted-mixing}
u_p(z)
+\sum_{s=1}^{\ell}
\left(h_{j_sp}-h_{j_{s+1}p}\right)
\mathfrak{z}_{j_s}
\geq h_{j_1p}.
\end{equation}
\end{corollary}

Inequalities~\eqref{eq:radius-adjusted-mixing} extend the classical mixing inequalities~\citep{gunluk2001mixing,honguyen2022rhs}, which can be separated in $O(N\log N)$ time for each safety constraint.

\subsection{Strict Probability Cuts and Rank Inequalities}\label{sec:residual-probability-cuts}

\paragraph{Strict probability cuts and the convex hull.}
Normalization into $\mathcal{D}_\theta^{\rm ctx}$ discards the scale of $w$ and the remaining transportation budget $B(w)$. However, the probability cut~\eqref{eq:probability-closure} can be strengthened at every integer point if $B(w)>0$.
\begin{subequations}
\begin{theorem}\label{thm:strict-probability-cuts}
Every feasible solution to model~\eqref{opt:exact-mip} satisfies
\begin{equation}\label{eq:strict-probability-cut}
w^\top\mathfrak{z}<\epsilon\mathbf{1}^\top w
\qquad
\forall w\in\mathcal{W}_\theta\text{ with }B(w)>0.
\end{equation}
Equivalently, for any fixed such $w$,
\begin{equation}\label{eq:predecessor-probability-cut}
w^\top\mathfrak{z}\leq\gamma(w),
\qquad
\gamma(w)
\coloneqq
\max\left\{w^\top a\;\middle|\;
a\in\{0,1\}^N,\ w^\top a<\epsilon\mathbf{1}^\top w
\right\}.
\end{equation}
If $\varepsilon=1$, then $w=\mathbf{1}/N$ and condition~\eqref{eq:strict-radius} gives $B(w)>0$, so~\eqref{eq:predecessor-probability-cut} becomes
\begin{equation}\label{eq:strict-uniform-cardinality}
\sum_{i \in [N]}\mathfrak{z}_i\leq\lceil\epsilon N\rceil-1.
\end{equation}
In particular, this improves the usual non-strict cardinality cut by one whenever $\epsilon N$ is an integer.
\end{theorem}
\end{subequations}
For an arbitrary $w$, computing $\gamma(w)$ is a binary knapsack problem. The case with at most two fractional components is nevertheless easy, and Theorem~\ref{thm:vertex-predecessor-hull} shows that every vertex of $\mathcal{W}_\theta$ has this property. For such a $w$, scale $\omega=Nw$, and define $F(w)=\{i\mid\omega_i=1\}$, $J(w)=\{i\mid0<\omega_i<1\}$, and $Z(w)=\{i\mid\omega_i=0\}$. Let $b(w)=\epsilon\sum_{i \in [N]}\omega_i$ and, for $a\in\{0,1\}^{J(w)}$, let $\alpha(a)=\sum_{j\in J(w)}\omega_ja_j$. Define the admissible patterns and the corresponding cardinality bounds by
\begin{equation}\label{eq:strict-vertex-patterns}
\mathcal{T}_\theta(w)=
\begin{cases}
\left\{a\;\middle|\;\alpha(a)\leq b(w)\right\}&B(w)=0\\
\left\{a\;\middle|\;\alpha(a)<b(w)\right\}&B(w)>0
\end{cases},
\qquad
k_a=\min\{|F(w)|,\;\widehat{k}_a\},
\end{equation}
where
\begin{equation}\label{eq:strict-vertex-capacity}
\widehat{k}_a=
\begin{cases}
\lfloor b(w)-\alpha(a)\rfloor&B(w)=0\\
\lceil b(w)-\alpha(a)\rceil-1&B(w)>0
\end{cases}.
\end{equation}
Thus, after pattern $a$ is fixed, $k_a$ is the largest admissible number of indices $i\in F(w)$ for which $\mathfrak z_i=1$.

\begin{theorem}\label{thm:vertex-predecessor-hull}
Every vertex of $\mathcal{W}_\theta$ has at most two fractional components. More generally, fix any $w\in\mathcal{W}_\theta$ with $|J(w)|\leq2$, and let $\mathcal{K}_\theta(w)$ be the binary set defined by $\omega^\top\mathfrak{z}\leq b(w)$ when $B(w)=0$ and by $\omega^\top\mathfrak{z}<b(w)$ when $B(w)>0$. Then $\operatorname{conv}(\mathcal{K}_\theta(w))$ is the projection onto $\mathfrak{z}$ of
\begin{subequations}\label{eq:vertex-predecessor-hull}
\begin{align}
&\sum_{a\in\mathcal{T}_\theta(w)}\psi_a=1\\
&\psi_a\geq0 &&\forall a\in\mathcal{T}_\theta(w) \label{eq:vertex-predecessor-pattern}\\
&\mathfrak{z}_j=\sum_{a\in\mathcal{T}_\theta(w)}a_j\psi_a
&&\forall j\in J(w) \label{eq:vertex-predecessor-special}\\
&\mathfrak{z}_i=\sum_{a\in\mathcal{T}_\theta(w)}y_i^a
&&\forall i\in F(w) \label{eq:vertex-predecessor-full}\\
&\sum_{i\in F(w)}y_i^a\leq k_a\psi_a
&&\forall a\in\mathcal{T}_\theta(w) \label{eq:vertex-predecessor-rank}\\
&0\leq y_i^a\leq\psi_a
&&\forall i\in F(w),\ a\in\mathcal{T}_\theta(w) \label{eq:vertex-predecessor-disagg}\\
&0\leq\mathfrak{z}_i\leq1
&&\forall i\in Z(w).\label{eq:vertex-predecessor-zero}
\end{align}
\end{subequations}
The formulation uses at most four pattern variables and $O(N)$ continuous variables and constraints. Equivalently, eliminating $y$ yields, in addition to nonnegativity,
\begin{equation}\label{eq:vertex-predecessor-rank-cuts}
\sum_{i\in S}\mathfrak{z}_i
\leq
\sum_{a\in\mathcal{T}_\theta(w)}
\psi_a\min\{|S|,k_a\}
\qquad
\forall S\subseteq F(w),
\end{equation}
which can be separated after sorting the components of $F(w)$.
\end{theorem}

The formulation is inexpensive for a minimum-radius allocation. An optimal solution $w^0$ of problem~\eqref{eq:min-radius-lp} can be chosen as a vertex and has at most one fractional component. Because $B(w^0)=\theta-\theta_{\min}>0$, it remains a vertex of $\mathcal{W}_\theta$, and~\eqref{eq:vertex-predecessor-hull} then needs at most two patterns. Further allocations can be collected from the probability-cut separation problem or the LP separation problem. Corollary~\ref{cor:critical-radius-cuts} ensures that the latter also returns allocations with at most two fractional $w$-components. Intersecting the formulations for several fixed allocations is valid, but may be weaker than the convex hull of their full intersection.

\paragraph{Rank inequalities across allocations.}
An important observation is that different allocations can exclude different binary patterns. To combine these restrictions, we first ask whether a set $S\subseteq[N]$ can supply all failure mass at the risk boundary without moving outcomes. Restricting the allocations in~\eqref{eq:critical-allocation-polytope} to failures supported on $S$ gives the minimum contextual cost
\begin{equation}\label{eq:contextual-rank-cost}
T(S)\coloneqq\min\left\{K_0+\bar\kappa^\top w\;\middle|\;
(w,r)\in\mathcal C_{\varepsilon,\epsilon},\quad r_i=0\quad\forall i\notin S\right\},
\end{equation}
with value $+\infty$ if infeasible. If $T(S)<\theta$, the adversary can reach the boundary using $S$ without exhausting the transportation budget. Thus, all labels in $S$ cannot equal one simultaneously. Larger contextual costs make a failure set more expensive to exploit, which identifies the most favorable set of any given size.

\begin{theorem}\label{thm:contextual-rank}
For $A\subseteq[N]$, let $S_k(A)$ consist of the $k$ indices in $A$ with the largest contextual costs, breaking ties by a fixed order, and set $S_0(A)=\varnothing$. Then $T(S)\leq T(S_k(A))$ for every $S\subseteq A$ with $|S|=k$. Define $R(A)=\max\{k\in\{0,\ldots,|A|\}\mid T(S_k(A))\geq\theta\}$. Every feasible solution to model~\eqref{opt:exact-mip} satisfies
\begin{equation}\label{eq:contextual-rank-cut}
\sum_{i\in A}\mathfrak z_i\leq R(A).
\end{equation}
The bound is sharp for the binary transportation projection in which $\vartheta(d)\geq\theta$, $d_i=0$ whenever $\mathfrak z_i=1$, and the other distances are unrestricted finite nonnegative values.
\end{theorem}

In particular, $T(S_k(A))<\theta$ excludes every $k$-subset of $A$, even when different subsets require different witness allocations. The contextual order identifies the extremal subset for the entire allocation family.

\section{Numerical Experiments}\label{sec:computational}
We conduct two sets of numerical experiments. First, we use a capacitated transportation problem, adapted from the benchmark instances in~\citet{chen2022data} and~\citet{honguyen2022rhs}, to evaluate the computational strength of the proposed MIP formulations and valid inequalities. We also use this problem to examine how training sample size and distribution shifts affect the out-of-sample cost and reliability of contextual and unconditional models. Second, we consider an anesthesiologist deployment problem using real-world data from~\citet{rath2026staff}. 
This case study examines whether contextual information helps allocate protection across different risk conditions and how this allocation affects operational cost and reliability.

All experiments are implemented in Julia 1.12.5 and solved with Gurobi 10.0.1 on an Apple M3 Max with 48~GB of memory. Each solver run uses one thread, \texttt{NumericFocus=2}, and primal and integer feasibility tolerances of $10^{-9}$. Unless otherwise stated, each final MIP solve has a time limit of $3{,}600$ seconds and a relative MIP gap tolerance of $10^{-8}$.

\subsection{Capacitated Transportation Problem}\label{sec:transportation}
Let $[F]$ and $[D]$ index factories and distribution centers, respectively. Shipping one unit from factory $f$ to distribution center $d$ incurs cost $c_{fd}\geq0$, and factory $f$ has capacity $m_f>0$. The decision $q_{fd}\geq0$ specifies the shipped quantity. The random vector $Y\in\R_+^D$ represents demands, while $X\in\R^K$ contains covariates observed before the shipment plan is chosen. Given a target context $x_0$, define  the conditioning neighborhood $\mathcal{N}_\gamma(x_0)=\{x\in\R^K \mid \|x-x_0\|_{\mathcal{X}}\leq\gamma\}$. The contextual transportation problem is
\begin{align*}
\min_{q\geq0}\quad
&\sum_{f\in[F]}\sum_{d\in[D]}c_{fd}q_{fd}\\*
\mathrm{s.t.}\quad
&\mathbb{P}\left(\sum_{f\in[F]}q_{fd} > Y_d\quad \forall d\in[D]\,\middle|\,X\in\mathcal{N}_\gamma(x_0)\right)\geq1-\epsilon && \forall \mathbb{P}\in\mathcal{P}_{\theta,\varepsilon}(\mathcal{N}_\gamma(x_0)) \span\\*
&\sum_{d\in[D]}q_{fd}\leq m_f &&\forall f\in[F].
\end{align*}
For a target context $x_0$, a shipment plan minimizes transportation cost subject to factory capacities and must strictly exceed demand at every distribution center with conditional probability at least $1-\epsilon$ under every distribution in $\mathcal{P}_{\theta,\varepsilon}(\mathcal{N}_\gamma(x_0))$. 

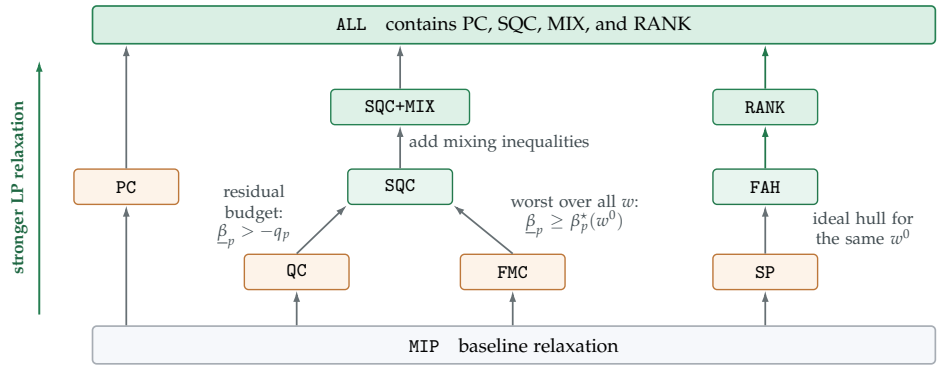
\begin{figure}[htbp]
  \centering
  \caption{Dominance relations among LP relaxations.}
  \label{fig:root-lp-dominance}
  \resizebox{0.75\linewidth}{!}{%
  \begin{tikzpicture}[x=1cm,y=1cm,>=Latex]
    \tikzset{
      dom node/.style={
        rounded corners=2pt,
        line width=0.65pt,
        minimum width=1.72cm,
        minimum height=0.58cm,
        align=center,
        font=\footnotesize\bfseries,
        inner sep=2.5pt
      },
      dom bar/.style={
        rounded corners=2pt,
        line width=0.78pt,
        minimum width=13.85cm,
        minimum height=0.62cm,
        align=center,
        font=\footnotesize\bfseries,
        inner sep=2.5pt
      },
      dom arrow/.style={
        -{Latex[length=1.9mm,width=1.15mm]},
        draw=ctxink!70,
        line width=0.82pt
      },
      dom reason/.style={
        rounded corners=1.5pt,
        fill opacity=0.94,
        text opacity=1,
        text=ctxink!88,
        align=center,
        font=\scriptsize,
        inner sep=2pt
      }
    }

    \node[dom bar,draw=ctxgray!68,fill=ctxlight] (mip) at (8.42,0.52)
      {\texttt{MIP}\quad\normalfont baseline relaxation};
    \node[dom bar,draw=ctxgreen!78!black,fill=ctxgreen!18] (all) at (8.42,5.78)
      {\texttt{ALL}\quad\normalfont contains PC, SQC, MIX, and RANK};

    \node[dom node,draw=ctxorange!78!black,fill=ctxorange!12] (pc)
      at (2.05,3.12) {\texttt{PC}};
    \node[dom node,draw=ctxorange!78!black,fill=ctxorange!12] (qc)
      at (4.85,1.72) {\texttt{QC}};
    \node[dom node,draw=ctxorange!78!black,fill=ctxorange!12] (fmc)
      at (8.40,1.72) {\texttt{FMC}};
    \node[dom node,draw=ctxorange!78!black,fill=ctxorange!12] (sp)
      at (12.55,1.72) {\texttt{SP}};

    \node[dom node,draw=ctxgreen!72!black,fill=ctxgreen!11] (sqc)
      at (6.55,3.12) {\texttt{SQC}};
    \node[dom node,minimum width=2.28cm,draw=ctxgreen!78!black,fill=ctxgreen!16]
      (sqcmix) at (6.55,4.43) {\texttt{SQC+MIX}};
    \node[dom node,draw=ctxgreen!72!black,fill=ctxgreen!11] (fah)
      at (12.55,3.12) {\texttt{FAH}};

    \node[dom node,draw=ctxgreen!78!black,fill=ctxgreen!16] (rank)
      at (12.55,4.43) {\texttt{RANK}};
    \draw[dom arrow,draw=ctxgreen!72!black] (fah.north) -- (rank.south);
    \draw[dom arrow,draw=ctxgreen!72!black] (rank.north) -- ([xshift=4.13cm]all.south);

    \draw[dom arrow] ([xshift=-6.37cm]mip.north) -- (pc.south);
    \draw[dom arrow] ([xshift=-3.57cm]mip.north) -- (qc.south);
    \draw[dom arrow] ([xshift=-0.02cm]mip.north) -- (fmc.south);
    \draw[dom arrow] ([xshift=4.13cm]mip.north) -- (sp.south);

    \draw[dom arrow] (qc.north) -- (sqc.south west);
    \draw[dom arrow] (fmc.north) -- (sqc.south east);
    \draw[dom arrow] (sqc.north) -- (sqcmix.south);

    \draw[dom arrow] (sp.north) -- (fah.south);

    \draw[dom arrow] (pc.north) -- ([xshift=-6.37cm]all.south);
    \draw[dom arrow] (sqcmix.north) -- ([xshift=-1.87cm]all.south);

    \node[dom reason,text width=2.05cm] at (4.15,2.63)
      {residual budget:\\[-1pt]
       $\underline\beta_p>-q_p$};
    \node[dom reason,text width=3.10cm] at (9.42,2.63)
      {worst over all $w$:\\[-1pt]
       $\underline\beta_p\geq\beta_p^\star(w^0)$};
    \node[dom reason] at (8.18,3.83)
      {add mixing inequalities};
    \node[dom reason,text width=2.35cm] at (14.16,2.44)
      {ideal hull for the same $w^0$};

    \draw[-{Latex[length=1.9mm,width=1.15mm]},draw=ctxgreen!76!black,
      line width=0.92pt] (0.62,1.02) -- (0.62,5.20);
    \node[rotate=90,font=\scriptsize\bfseries,text=ctxgreen!62!black]
      at (0.29,3.11) {stronger LP relaxation};
  \end{tikzpicture}%
  }
  \par\footnotesize Note. Arrows indicate established containment relations for the exact inequality families. 
\end{figure}

We compare ten variants of the MIP reformulation. The baseline MIP formulation, denoted by \texttt{MIP}, is given by~\eqref{opt:exact-mip}. \texttt{PC}, the probability-cut closure, adds~\eqref{eq:compact-probability-closure}. \texttt{QC}, the quantile cut, adds~\eqref{eq:probability-quantile-row-bound} and the tightened coefficients~\eqref{eq:probability-quantile-M}. Let $w^0$ be a minimum-radius allocation obtained from~\eqref{eq:min-radius-lp}. \texttt{FMC}, the fixed-allocation margin cut, adds $\beta_p(z)\geq\beta_p^\star(w^0)$ from~\eqref{eq:fixed-allocation-threshold-definition}. \texttt{SQC}, the strengthened quantile cut, adds~\eqref{eq:radius-row-bound}, which accounts for the residual transportation budget, and the tightened coefficients~\eqref{eq:radius-adjusted-M}. \texttt{SQC}+\texttt{MIX} augments \texttt{SQC} with the mixing inequalities~\eqref{eq:radius-adjusted-mixing}, separated at the root. \texttt{SP}, the strict probability cut, adds~\eqref{eq:predecessor-probability-cut} generated by $w^0$, whereas \texttt{FAH}, the fixed-allocation hull, installs the ideal formulation~\eqref{eq:vertex-predecessor-hull} for the same allocation. \texttt{RANK} adds the prefix inequalities~\eqref{eq:contextual-rank-cut} to \texttt{MIP}. We choose $w^0$ greedily in nondecreasing order of $\bar\kappa_i$, using the same order for \texttt{RANK}. \texttt{ALL} combines~\eqref{eq:compact-probability-closure},~\eqref{eq:radius-row-bound},~\eqref{eq:radius-adjusted-M},~\eqref{eq:radius-adjusted-mixing}, and~\eqref{eq:contextual-rank-cut}. Figure~\ref{fig:root-lp-dominance} summarizes the theoretical partial order.

We measure computational performance using relaxation strength, final optimality gaps, elapsed time, and the search-tree size. We distinguish three measures of computation time. The root time $T^{\rm root}$ includes the LP solves and any separation of mixing inequalities at the root. The MIP time $T^{\rm MIP}$ records the elapsed time of the final MIP solve. The total algorithm time $T^{\rm alg}$ additionally includes preprocessing, model construction, and any root cut generation and transfer required by the method. For method $m$, let $V_m^{\rm LP}$ be the root LP relaxation value and let $V^{\rm ref}$ be the optimal value. The root-relaxation gap is $g_m^{\rm LP}=\frac{\max\{V^{\rm ref}-V_m^{\rm LP},0\}}{|V^{\rm ref}|}\times100\%$. We also report the final relative gap $g^{\rm final}$ between the incumbent objective and the best valid bound, and the number $B$ of branch-and-bound nodes. Gaps are reported as percentages and computation times in seconds.

\subsubsection{Strength of the Valid Inequalities}\label{sec:computational-results}

Let $F$, $D$, and $K$ denote the numbers of factories, distribution centers, and contextual features, respectively. The aggregate experiment uses three network scales $(F,D,K)\in\{(5,20,3)$, $(10,50,3)$, $(15,100,3)\}$, nested training sizes $N\in\{50,100,500\}$ drawn from a master sample of size $500$, and four radius labels (the minimum, denoted by NM, and the normalized interpolation fractions $0.1$, $0.5$, and $1.0$). Each design is evaluated at target contexts corresponding to low, central, and high demand. 

\paragraph{Root-relaxation strength.} Table~\ref{tab:formal-root-strength} shows that the exact MIP reformulation has a weak relaxation, with mean and median root gaps of $89.48\%$ and $95.34\%$. The theoretical relations among the strengthening blocks are reflected in the results, but they do not form a single dominance chain. Within the margin-cut family,~\eqref{eq:radius-dominates-probability} and~\eqref{eq:fixed-allocation-supremum} imply that \texttt{SQC} dominates \texttt{QC} and \texttt{FMC}, respectively. Consistently, the mean root gap is $2.95\%$ under \texttt{SQC}, compared with $3.56\%$ under \texttt{QC} and $8.64\%$ under \texttt{FMC}, while \texttt{QC} and \texttt{FMC} themselves are not generally ordered. \texttt{SQC}+\texttt{MIX} is theoretically no weaker than \texttt{SQC}. Based on the computational results, \texttt{SQC}+\texttt{MIX} produces root gaps nearly identical to those of \texttt{SQC}. Within the probability-cut family, \texttt{FAH} provides an ideal formulation of the binary set underlying \texttt{SP} for a fixed allocation and is at least as strong as \texttt{SP}. The mean gaps of $17.41\%$ for \texttt{FAH} and $18.07\%$ for \texttt{SP} indicate a modest improvement. \texttt{PC} imposes all non-strict probability cuts, whereas \texttt{SP} and \texttt{FAH} exploit strictness for one allocation, so these formulations are not generally ordered. Although \texttt{FMC}, \texttt{SP}, and \texttt{FAH} all use $w^0$, \texttt{FMC} strengthens the safety margins, whereas \texttt{SP} and \texttt{FAH} strengthen the constraints on the sample failure indicators, so the shared allocation does not imply an ordering between the two families. \texttt{RANK} reduces the mean gap to $6.52\%$, improving on \texttt{SP} and \texttt{FAH}, while \texttt{SQC} gives a smaller mean gap of $2.95\%$. These results show that restrictions on failure patterns complement the margin inequalities. \texttt{ALL} combines their complementary strength and gives the \emph{smallest} mean and median gaps, $0.69\%$ and $0.02\%$.

\begin{table}[htbp]
\centering
\caption{Aggregate root-relaxation strength and computational performance.}
\label{tab:formal-root-strength}
\setlength{\tabcolsep}{5.5pt}
\resizebox{0.9\linewidth}{!}{%
\begin{tabular}{l|cccc|ccc|cc}\toprule[1pt]\midrule
& \multicolumn{4}{c|}{Root relaxation} & \multicolumn{3}{c|}{Final MIP solve} & \multicolumn{2}{c}{Nodes}\\\cmidrule(lr){2-5} \cmidrule(lr){6-8} \cmidrule(lr){9-10}
Method & Mean $g^{\rm LP}$ & Med. $g^{\rm LP}$ & Med. $T^{\rm root}$ & Max $T^{\rm root}$ & Max $g^{\rm final}$ & Med. $T^{\rm MIP}$ & Med. $T^{\rm alg}$ & Med. $B$ & Max $B$\\\midrule
\texttt{MIP} & 89.48 & 95.34 & 0.040 & 0.494 & 2.934 & 0.344 & 0.353 & 42.5 & 326{,}715\\
\texttt{PC} & 19.58 & 18.71 & 0.046 & 0.955 & 0.000 & 0.365 & 0.380 & 26.5 & 24{,}190\\
\texttt{QC} & 3.56 & 0.20 & 0.058 & 0.808 & 0.630 & 0.244 & 0.347 & 1.5 & 195{,}942\\
\texttt{FMC} & 8.64 & 9.69 & 0.061 & 0.819 & 0.000 & 0.332 & 0.352 & 6 & 130{,}976\\
\texttt{SQC} & 2.95 & 0.16 & 0.062 & 0.820 & 0.000 & 0.249 & 0.589 & 1 & 132{,}718\\
\texttt{SQC+MIX} & 2.83 & 0.15 & 0.092 & 6.480 & 0.000 & 0.248 & 0.666 & 1 & 123{,}629\\
\texttt{SP} & 18.07 & 18.31 & 0.029 & 0.864 & 0.000 & 0.339 & 0.348 & 20 & 17{,}001\\
\texttt{FAH} & 17.41 & 17.70 & 0.029 & 1.236 & 0.000 & 0.336 & 0.345 & 20 & 17{,}001\\
\texttt{RANK} & 6.52 & 0.40 & 0.029 & 1.656 & 0.000 & 0.133 & 0.144 & 1 & 10{,}046\\
\texttt{ALL} & 0.69 & 0.02 & 0.084 & 1.594 & 0.000 & 0.178 & 0.613 & 1 & 679\\
\midrule\bottomrule[1pt]
\end{tabular}}
\end{table}

\paragraph{Final MIP performance and cost.} Stronger roots can shrink the search tree, but preprocessing can dominate total time. \texttt{SP} and \texttt{FAH} are the lightest strengthenings: their median root time is $0.029$ seconds and their median algorithm times are $0.348$ and $0.345$ seconds, close to $0.353$ seconds for the \texttt{MIP} formulation, while the median and maximum node counts fall from $42.5$ and $326{,}715$ to $20$ and $17{,}001$. \texttt{QC}, \texttt{SQC}, and \texttt{SQC}+\texttt{MIX} lower the median optimization time from $0.344$ seconds to $0.244$--$0.249$ seconds. Yet \texttt{SQC} and \texttt{SQC}+\texttt{MIX} have median algorithm times of $0.589$ and $0.666$ seconds. However, \texttt{MIX} gives only a modest improvement in the root bound and raises the maximum time spent at the root node relaxation. Across the 108 instances, \texttt{RANK} has a median MIP solution time of $0.133$ seconds and a median total algorithm time of $0.144$ seconds. The corresponding medians for \texttt{ALL} are $0.178$ and $0.613$ seconds, respectively. Both have a median of one node, and their maximum node counts are $10{,}046$ and $679$, respectively. All final gaps are zero except the maxima for the \texttt{MIP} formulation ($2.934\%$) and the \texttt{QC} formulation ($0.630\%$).

\paragraph{Scale dependence.}
This scale extension also includes $(30,100,3)$ with training sizes up to $N=1{,}000$. Table~\ref{tab:formal-scale-performance} in Appendix~\ref{sec:additionaltables} disaggregates performance by network and training-sample size. For $N\leq100$, every formulation has a median final solving time below 1 second, so the stronger relaxations provide little practical benefit at these sizes. At $N=500$, \texttt{ALL} requires less MIP solution time and explores fewer branch-and-bound nodes than \texttt{RANK} on all 48 instances. The median MIP times are $5.12$ and $14.88$ seconds, respectively. When preprocessing, model construction, and root separation are included, \texttt{ALL} has lower total algorithm time on 34 of the 48 instances, with medians of $13.06$ and $15.04$ seconds, respectively. Thus, these additional costs reduce the time savings from the smaller search trees.

The $N=1{,}000$ extension in Table~\ref{tab:n1000-instance-performance} makes this difference more pronounced and shows the limits of \texttt{RANK} alone. The median root gaps are $16.95\%$ for \texttt{RANK} and $2.63\%$ for \texttt{ALL}. \texttt{ALL} reaches the prescribed relative gap in all 12 instances, with a median MIP time of $39.67$ seconds. \texttt{RANK} reaches the target gap in six instances. The remaining $6$ runs reach the time limit with feasible incumbents. The median recorded MIP time for \texttt{RANK} is $3{,}520.14$ seconds. \texttt{MIP} and \texttt{SQC} reach the time limit in 12 and 5 instances, respectively.

\begin{table}[htbp]\centering
\caption{Instance-level performance for the $(30,100)$, $N=1{,}000$ scale extension.}
\label{tab:n1000-instance-performance}
\setlength{\tabcolsep}{1.5pt}
\resizebox{\linewidth}{!}{%
\begin{tabular}{ll|rrrr|rrrr|rrrr|rrrr}\toprule[1pt]\midrule
 & & \multicolumn{4}{c}{\texttt{MIP}} & \multicolumn{4}{c}{\texttt{SQC}} & \multicolumn{4}{c}{\texttt{RANK}} & \multicolumn{4}{c}{\texttt{ALL}}\\
 \cmidrule(lr){3-6} \cmidrule(lr){7-10} \cmidrule(lr){11-14} \cmidrule(lr){15-18}
Query & Radius & $g^{\rm LP}$ & $g^{\rm final}$ & $T^{\rm MIP}$ & $B$ & $g^{\rm LP}$ & $g^{\rm final}$ & $T^{\rm MIP}$ & $B$ & $g^{\rm LP}$ & $g^{\rm final}$ & $T^{\rm MIP}$ & $B$ & $g^{\rm LP}$ & $g^{\rm final}$ & $T^{\rm MIP}$ & $B$\\\midrule
\multirow{4}{*}{Low} & NM & 99.96 & 11.969 & 3600.35 & 97{,}716 & 14.96 & 5.019 & 3600.61 & 210{,}097 & 17.10 & 1.481 & 3600.61 & 93{,}674 & 0.52 & 0.008 & 46.29 & 1{,}449\\
 & 0.1 & 96.65 & 9.168 & 3600.39 & 64{,}472 & 10.27 & 0.884 & 3600.57 & 141{,}896 & 16.61 & 2.263 & 3600.36 & 58{,}153 & 3.10 & 0.008 & 209.33 & 6{,}293\\
 & 0.5 & 84.87 & 5.280 & 3600.47 & 37{,}026 & 5.90 & 0.010 & 74.48 & 2{,}530 & 14.25 & 0.009 & 2402.15 & 30{,}559 & 3.29 & 0.007 & 33.04 & 933\\
 & 1.0 & 74.20 & 5.303 & 3600.35 & 33{,}260 & 4.92 & 0.005 & 54.58 & 2{,}597 & 13.60 & 0.010 & 2065.96 & 32{,}583 & 2.69 & 0.007 & 31.32 & 1{,}282\\
\addlinespace[2pt]\cmidrule{1-18}\addlinespace[2pt]
\multirow{4}{*}{Central} & NM & 99.96 & 6.612 & 3600.41 & 97{,}301 & 15.06 & 0.714 & 3600.85 & 180{,}303 & 18.00 & 0.009 & 436.70 & 7{,}985 & 0.14 & 0.000 & 10.67 & 1\\
 & 0.1 & 97.33 & 1.316 & 3600.32 & 51{,}747 & 10.55 & 0.010 & 353.76 & 11{,}958 & 18.72 & 0.007 & 1012.70 & 12{,}191 & 2.13 & 0.010 & 20.21 & 410\\
 & 0.5 & 88.10 & 5.084 & 3600.19 & 44{,}771 & 6.92 & 0.009 & 496.80 & 20{,}598 & 15.38 & 1.962 & 3600.22 & 56{,}296 & 3.45 & 0.007 & 220.34 & 8{,}468\\
 & 1.0 & 79.07 & 4.367 & 3600.18 & 37{,}276 & 5.49 & 0.008 & 64.15 & 3{,}234 & 14.17 & 0.010 & 1725.18 & 27{,}582 & 2.58 & 0.000 & 27.63 & 949\\
\addlinespace[2pt]\cmidrule{1-18}\addlinespace[2pt]
\multirow{4}{*}{High} & NM & 99.97 & 8.874 & 3600.41 & 106{,}634 & 16.87 & 3.865 & 3600.43 & 216{,}221 & 19.38 & 0.010 & 3440.09 & 69{,}919 & 0.27 & 0.005 & 15.68 & 140\\
 & 0.1 & 98.53 & 14.211 & 3600.47 & 72{,}054 & 15.18 & 3.293 & 3600.68 & 180{,}789 & 21.03 & 3.387 & 3600.20 & 59{,}663 & 2.45 & 0.008 & 319.15 & 14{,}254\\
 & 0.5 & 93.04 & 6.878 & 3600.42 & 46{,}980 & 10.01 & 0.010 & 1440.76 & 61{,}498 & 19.18 & 3.020 & 3600.53 & 47{,}096 & 3.70 & 0.006 & 123.27 & 4{,}363\\
 & 1.0 & 86.67 & 3.429 & 3600.27 & 44{,}682 & 7.29 & 0.010 & 184.16 & 7{,}548 & 16.80 & 1.444 & 3600.20 & 46{,}864 & 2.96 & 0.005 & 50.27 & 1{,}822\\
\midrule\bottomrule[1pt]\end{tabular}}
\end{table}

\begin{figure}[htbp]
  \centering
  \caption{Cumulative fraction of instances solved over time for \texttt{MIP}, \texttt{SQC}, \texttt{RANK}, and \texttt{ALL}.}
  \label{fig:solved-fraction-profiles}
  \begin{minipage}[t]{0.49\linewidth}
    \centering
    \includegraphics[width=\linewidth]{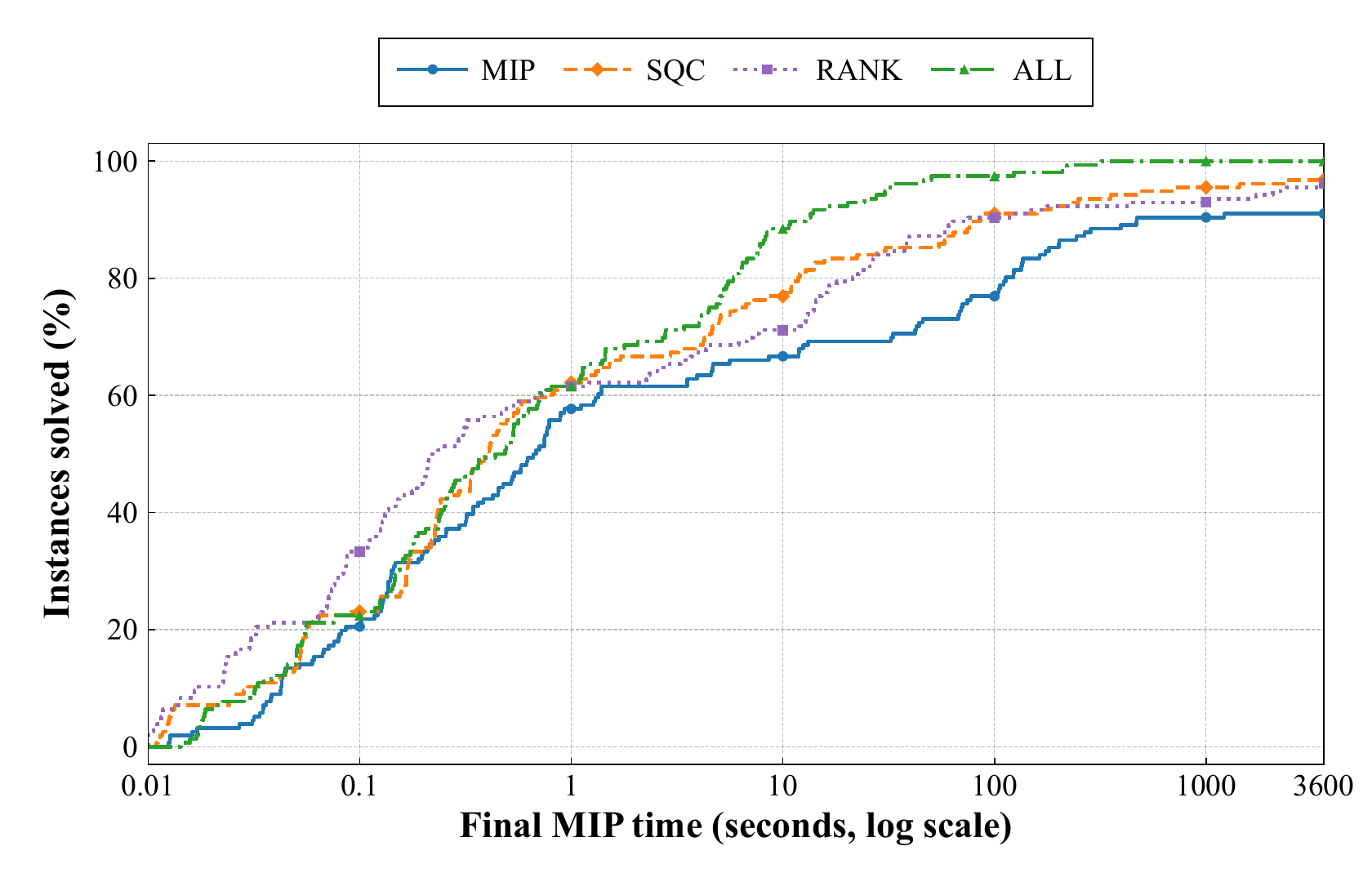}
    \par\smallskip\footnotesize\textbf{(a)} All instances.
  \end{minipage}\hfill
  \begin{minipage}[t]{0.49\linewidth}
    \centering
    \includegraphics[width=\linewidth]{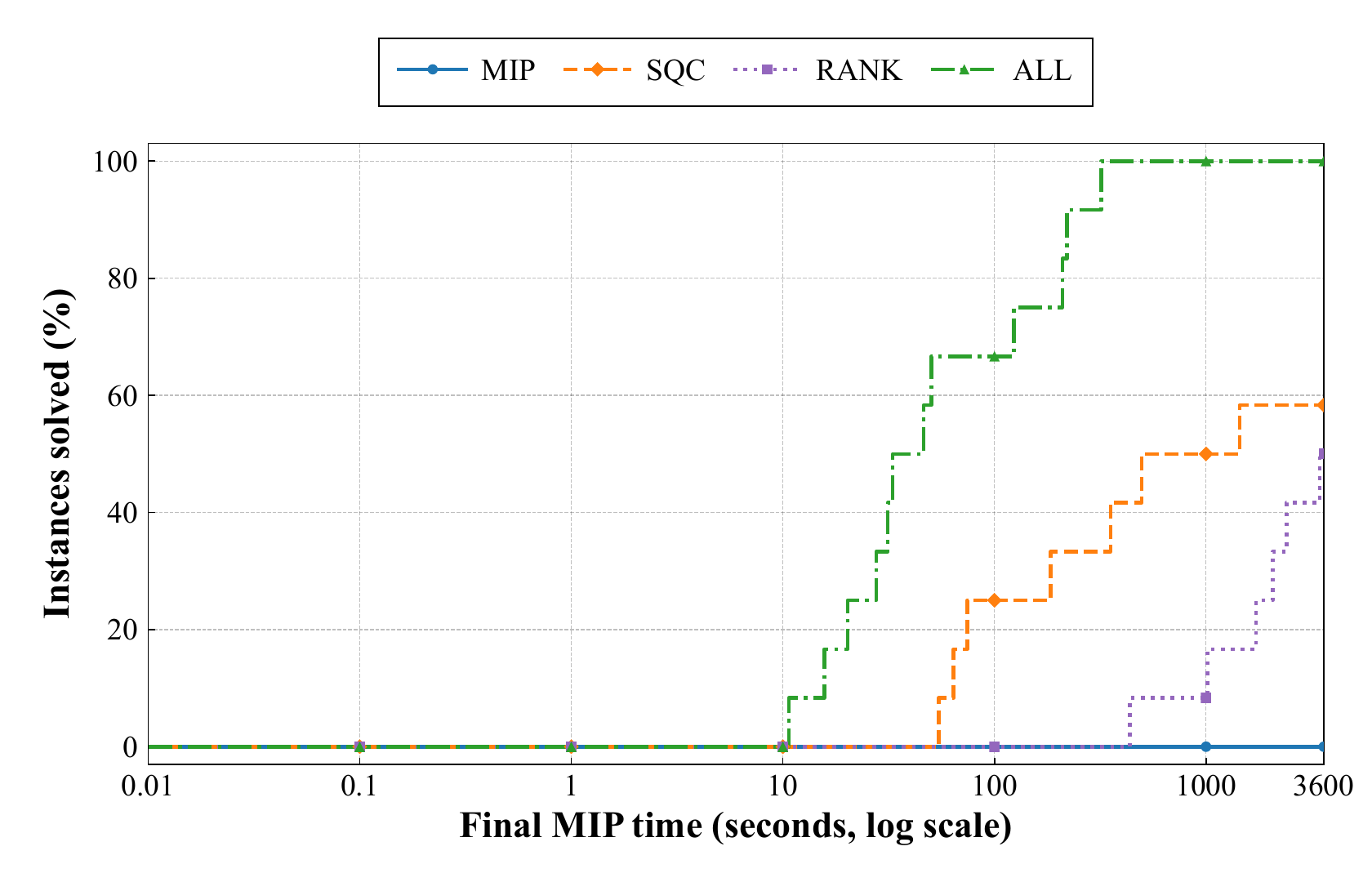}
    \par\smallskip\footnotesize\textbf{(b)} Instances with $(F,D)=(30,100)$ and $N=1{,}000$.
  \end{minipage}
\end{figure}


\subsubsection{Benchmark Comparison and Out-of-Sample Performance}
\label{sec:oos-performance}

We compare the C-DRCCP model with three benchmark models: the unconditional sample-average approximation (SAA) model, the unconditional DRCCP model, and the contextual SAA (C-SAA) model. The unconditional SAA and DRCCP models ignore the contextual information $X$. The C-SAA model uses the observations in $\mathcal{N}_\gamma(x_0)$ and permits at most $\lfloor\epsilon n_{\rm local}\rfloor$ joint violations among the resulting $n_{\rm local}$ outcomes. A C-SAA case with $n_{\rm local}=0$ is recorded as unavailable. The unconditional DRCCP model is the algebraic special case of our formulation with zero contextual transportation cost and conditioning mass fixed at one. Appendix~\ref{appendix:model-relations} gives the exact relation and distinguishes this specialization from merely assigning zero cost to a conditioning event whose mass remains variable. The C-DRCCP model is solved with the \texttt{ALL} formulation.

The out-of-sample tests use $10$ independent training replications. We study the networks $(F,D)\in\{(10,50),(15,100)\}$, the contexts corresponding to low, central, and high demand, and nested training sizes $N\in\{50,100\}$. The main design sets $\epsilon=0.10$ and, for the C-DRCCP model, sets $\varepsilon$ equal to the independently calibrated $99\%$ Hoeffding lower bound on the neighborhood probability. The contextual and unconditional radius grids are calibrated separately. Single-factor sensitivity analyses vary $\epsilon\in\{0.05,0.10,0.20\}$ and $\varepsilon/\underline p\in\{0.50,0.75,1.00\}$ on the $(10,50)$ network at the central query. Each decision is evaluated using two separate out-of-sample datasets: $5{,}000$ neighborhood observations and $5{,}000$ pointwise observations.

\paragraph{Evaluation metrics.}
For every feasible shipment plan $q$, we estimate the neighborhood risk $p_{\mathcal{N}}(q)=\mathbb P(\exists d:\ Y_d\geq\sum_fq_{fd} \mid X\in\mathcal{N}_\gamma(x_0))$ and the pointwise risk $p_{x_0}(q)$ obtained by replacing the conditioning event with $X=x_0$. Reliability (Rel.) is the percentage of training replications for which the one-sided $95\%$ upper confidence bound on the violation probability does not exceed $\epsilon$. The tables report mean cost and mean estimated risks, denoted by $\widehat p_{\mathcal N}$ and $\widehat p_{x_0}$, across the training replications specified in the experimental design. We use $\bar p_{\mathcal N}$ and $\bar p_{x_0}$ for risk averages over multiple experimental settings.

\paragraph{Nominal performance.}
Table~\ref{tab:oos-nominal-neighborhood} reports mean transportation cost, neighborhood risk, and reliability under the nominal DGP. The C-SAA model is unreliable across all reported cells, as its local training sample contains only $2$--$18$ observations. Robustification corrects this behavior: the C-DRCCP model lowers the violation probability relative to the C-SAA model in each case, but at a higher cost. The comparison with the unconditional DRCCP model is a cost--risk trade-off. The unconditional model has the lowest raw violation rate, whereas the C-DRCCP model is generally less costly. The C-DRCCP model is reliable at low query and in most central-query cells, but the high-query and larger-network cases provide important counterexamples to a uniform reliability claim. Table~\ref{tab:oos-pointwise} evaluates the same decisions at $X=x_0$, for which costs are unchanged from Table~\ref{tab:oos-nominal-neighborhood} and are omitted here. Model performance is similar at the low and central queries. 
\begin{table}[htbp]
\centering
\caption{Cost, neighborhood risk, and reliability under the nominal DGP.}
\label{tab:oos-nominal-neighborhood}
\setlength{\tabcolsep}{2.0pt}
\resizebox{.9\linewidth}{!}{%
\begin{tabular}{ccl|ccc|ccc|ccc|ccc}\toprule[1pt]\midrule
&&&\multicolumn{3}{c|}{SAA}&\multicolumn{3}{c|}{C-SAA}&\multicolumn{3}{c|}{DRCCP}&\multicolumn{3}{c}{C-DRCCP}\\
\cmidrule(lr){4-6}\cmidrule(lr){7-9}\cmidrule(lr){10-12}\cmidrule(lr){13-15}
$(F,D)$&$N$&Query&Cost&$\widehat{p}_{\mathcal{N}}$&Rel. [\%]&Cost&$\widehat{p}_{\mathcal{N}}$&Rel. [\%]&Cost&$\widehat{p}_{\mathcal{N}}$&Rel. [\%]&Cost&$\widehat{p}_{\mathcal{N}}$&Rel. [\%]\\
\midrule
\multirow{6}{*}{$(10,50)$}
&\multirow{3}{*}{50}&Low     &837.0&0.106&30&609.8&0.896&0&1002.7&0.008&100&969.2&0.009&100\\
&&Central &837.0&0.321&0 &668.2&0.926&0&1002.7&0.044&80 &975.0&0.051&100\\
&&High    &837.0&0.691&0 &735.5&0.944&0&1002.7&0.227&0  &963.5&0.332&0\\
\cmidrule{2-15}
&\multirow{3}{*}{100}&Low     &892.0&0.037&100&625.4&0.847&0&1078.9&0.001&100&960.5&0.011&100\\
&&Central &892.0&0.154&10 &695.8&0.868&0&1078.9&0.011&100&990.3&0.052&90\\
&&High    &892.0&0.487&0  &775.9&0.885&0&1078.9&0.094&70 &1008.8&0.213&10\\
\midrule
\multirow{6}{*}{$(15,100)$}
&\multirow{3}{*}{50}&Low     &2368.2&0.101&40&1508.8&0.960&0&2881.0&0.006&100&2662.9&0.028&90\\
&&Central &2368.2&0.349&0 &1783.1&0.960&0&2881.0&0.050&90 &2664.2&0.181&50\\
&&High    &2368.2&0.710&0 &2051.0&0.950&0&2881.0&0.240&10 &2659.9&0.455&0\\
\cmidrule{2-15}
&\multirow{3}{*}{100}&Low     &2506.0&0.033&100&1594.6&0.914&0&3075.1&0.002&100&2763.1&0.011&100\\
&&Central &2506.0&0.179&10 &1877.4&0.888&0&3075.1&0.014&100&2746.4&0.104&80\\
&&High    &2506.0&0.522&0  &2154.5&0.901&0&3075.1&0.112&40 &2753.6&0.347&0\\
\midrule\bottomrule[1pt]
\end{tabular}
}
\end{table}

\begin{table}[htbp]
\centering
\caption{Pointwise risk and reliability under the nominal DGP.}
\label{tab:oos-pointwise}
\setlength{\tabcolsep}{6pt}
\resizebox{0.8\linewidth}{!}{%
\begin{tabular}{ccl|cc|cc|cc|cc}\toprule[1pt]\midrule
&&&\multicolumn{2}{c|}{SAA}&\multicolumn{2}{c|}{C-SAA}&\multicolumn{2}{c|}{DRCCP}&\multicolumn{2}{c}{C-DRCCP}\\
\cmidrule(lr){4-5}\cmidrule(lr){6-7}\cmidrule(lr){8-9}\cmidrule(lr){10-11}
$(F,D)$&$N$&Query&$\widehat{p}_{x_0}$&Rel. [\%]&$\widehat{p}_{x_0}$&Rel. [\%]&$\widehat{p}_{x_0}$&Rel. [\%]&$\widehat{p}_{x_0}$&Rel. [\%]\\
\midrule
\multirow{6}{*}{$(10,50)$}
&\multirow{3}{*}{50}&Low     &0.043&100&0.853&0&0.001&100&0.001&100\\
&&Central &0.276&10 &0.921&0&0.028&100&0.033&100\\
&&High    &0.704&0  &0.959&0&0.199&20 &0.309&0\\
\cmidrule{2-11}
&\multirow{3}{*}{100}&Low     &0.009&100&0.787&0&0.000&100&0.002&100\\
&&Central &0.113&50 &0.855&0&0.005&100&0.036&90\\
&&High    &0.475&0  &0.904&0&0.066&70 &0.183&40\\
\midrule
\multirow{6}{*}{$(15,100)$}
&\multirow{3}{*}{50}&Low     &0.039&100&0.943&0&0.002&100&0.010&100\\
&&Central &0.279&10 &0.960&0&0.026&100&0.146&70\\
&&High    &0.733&0  &0.964&0&0.217&10 &0.449&0\\
\cmidrule{2-11}
&\multirow{3}{*}{100}&Low     &0.009&100&0.878&0&0.000&100&0.004&100\\
&&Central &0.119&30 &0.882&0&0.005&100&0.074&90\\
&&High    &0.521&0  &0.921&0&0.087&50 &0.326&0\\
\midrule\bottomrule[1pt]
\end{tabular}
}
\end{table}

\paragraph{Distributional robustness.}
To assess robustness, we keep every decision fixed at its nominally trained value and change one interpretable DGP mechanism at a time:
\begin{itemize}[leftmargin=1.6em,itemsep=2pt,topsep=3pt]
\item \emph{Demand shift.} We set $Y^{(m)}=mY$ for $m\in\{1.05,1.10,1.20\}$. This shift represents a systematic underestimation of demand and scales both its level and absolute dispersion.
\item \emph{Volatility shift.} We multiply the standard deviations of the residuals by $s_\sigma\in\{1.25,1.50\}$. The lognormal term is recentered to preserve $\mathbb E[Y_d\mid X]$ to isolate changes in conditional dispersion and tail risk.
\item \emph{Dependence shift.} We vary the residual correlation over $\rho_{\rm res}\in\{0.50,0.70,0.90\}$ while preserving every center's conditional mean and variance. This shift changes the extent to which extreme demands occur simultaneously across centers.
\item \emph{Covariate shift.} We shift the marginal distribution of contexts away from its nominal location at two magnitudes, $\delta_X\in\{0.5,1.0\}$, while leaving the conditional demand distribution $P(Y\mid X)$ unchanged. This changes which contexts populate the fixed query neighborhood.
\end{itemize}

Tables~\ref{tab:oos-shift-summary} and~\ref{tab:oos-shift-paths} report neighborhood and pointwise risk, respectively. Each table cell aggregates $120$ evaluations: two networks, two training sizes, three queries, and ten training replications. Table~\ref{tab:oos-shift-summary} shows how neighborhood risk changes at each shift level. Higher demand and greater residual volatility increase the mean violation probability and reduce reliability across all models. The C-DRCCP model remains substantially more reliable than the C-SAA model at every level, while the DRCCP model has the lowest mean violation probability in every reported setting. The covariate shift has a smaller effect. In contrast, increasing residual correlation lowers the probability that at least one center fails. Because the marginal distributions remain fixed, stronger positive dependence makes failures at different centers more likely to overlap, thereby reducing the probability of their union.

\begin{table}[htbp]
\centering
\caption{Neighborhood risk and reliability under distribution shifts.}
\label{tab:oos-shift-summary}
\setlength{\tabcolsep}{2.5pt}
\resizebox{.8\linewidth}{!}{%
\begin{tabular}{ll|cc|cc|cc|cc}\toprule[1pt]\midrule
&&\multicolumn{2}{c|}{SAA}&\multicolumn{2}{c|}{C-SAA}&\multicolumn{2}{c|}{DRCCP}&\multicolumn{2}{c}{C-DRCCP}\\
Shift family&Level&$\bar{p}_{\mathcal{N}}$&Rel. [\%]&$\bar{p}_{\mathcal{N}}$&Rel. [\%]&$\bar{p}_{\mathcal{N}}$&Rel. [\%]&$\bar{p}_{\mathcal{N}}$&Rel. [\%]\\
\cmidrule{1-10}
Nominal DGP&---&0.308&24.2&0.912&0.0&0.067&74.2&0.149&60.0\\
\cmidrule{1-10}
\multirow{3}{*}{Demand level}
&$m=1.05$&0.433&13.3&0.962&0.0&0.123&58.3&0.229&45.8\\
&$m=1.10$&0.560&1.7&0.985&0.0&0.200&43.3&0.327&26.7\\
&$m=1.20$&0.777&0.0&0.998&0.0&0.401&15.8&0.548&1.7\\
\cmidrule{1-10}
\multirow{2}{*}{Residual volatility}
&$s_\sigma=1.25$&0.453&2.5&0.941&0.0&0.154&50.0&0.259&28.3\\
&$s_\sigma=1.50$&0.583&0.0&0.957&0.0&0.268&20.8&0.385&4.2\\
\cmidrule{1-10}
\multirow{3}{*}{Residual dependence}
&$\rho_{\rm res}=0.50$&0.243&31.7&0.834&0.0&0.052&79.2&0.118&62.5\\
&$\rho_{\rm res}=0.70$&0.182&44.2&0.738&0.0&0.038&87.5&0.089&68.3\\
&$\rho_{\rm res}=0.90$&0.123&55.0&0.622&0.0&0.025&94.2&0.062&77.5\\
\cmidrule{1-10}
\multirow{2}{*}{Covariate}
&$\delta_X=0.50$&0.326&23.3&0.922&0.0&0.074&70.8&0.161&60.0\\
&$\delta_X=1.00$&0.344&21.7&0.931&0.0&0.083&67.5&0.174&57.5\\
\midrule\bottomrule[1pt]
\end{tabular}
}
\end{table}

Table~\ref{tab:oos-shift-paths} evaluates the same decisions at $X=x_0$. Higher demand and greater residual volatility increase the probability of violation and reduce reliability. At every reported level, the DRCCP model has the lowest mean violation probability and the highest reliability, followed by the C-DRCCP model, the SAA model, and the C-SAA model. Greater residual correlation lowers pointwise risk for the same reason that it lowers neighborhood risk: failures at different centers overlap more often. We omit the covariate shift because it leaves $P(Y\mid X=x_0)$ unchanged.

\begin{table}[htbp]
\centering
\caption{Pointwise risk and reliability under distribution shifts.}
\label{tab:oos-shift-paths}
\setlength{\tabcolsep}{2.5pt}
\resizebox{.8\linewidth}{!}{%
\begin{tabular}{ll|cc|cc|cc|cc}\toprule[1pt]\midrule
&&\multicolumn{2}{c|}{SAA}&\multicolumn{2}{c|}{C-SAA}&\multicolumn{2}{c|}{DRCCP}&\multicolumn{2}{c}{C-DRCCP}\\
Shift family&Level&$\bar{p}_{x_0}$&Rel. [\%]&$\bar{p}_{x_0}$&Rel. [\%]&$\bar{p}_{x_0}$&Rel. [\%]&$\bar{p}_{x_0}$&Rel. [\%]\\
\cmidrule{1-10}
Nominal DGP&---&0.277&41.7&0.902&0.0&0.053&79.2&0.131&65.8\\
\cmidrule{1-10}
\multirow{3}{*}{Demand level}
&$m=1.05$&0.395&23.3&0.960&0.0&0.103&65.0&0.206&55.0\\
&$m=1.10$&0.515&13.3&0.985&0.0&0.176&54.2&0.297&35.8\\
&$m=1.20$&0.732&0.0 &0.999&0.0&0.366&26.7&0.502&13.3\\
\cmidrule{1-10}
\multirow{2}{*}{Residual volatility}
&$s_\sigma=1.25$&0.421&15.8&0.936&0.0&0.136&55.8&0.237&40.0\\
&$s_\sigma=1.50$&0.552&0.0 &0.955&0.0&0.248&30.0&0.360&18.3\\
\cmidrule{1-10}
\multirow{3}{*}{Residual dependence}
&$\rho_{\rm res}=0.50$&0.213&46.7&0.812&0.0&0.040&85.0&0.100&66.7\\
&$\rho_{\rm res}=0.70$&0.154&53.3&0.702&0.0&0.028&89.2&0.073&74.2\\
&$\rho_{\rm res}=0.90$&0.099&61.7&0.571&0.0&0.017&96.7&0.048&83.3\\
\midrule\bottomrule[1pt]
\end{tabular}
}
\end{table}

\paragraph{Value of contextual information.}
Tables~\ref{tab:fixed-raw-radius-neighborhood} and~\ref{tab:fixed-raw-radius-pointwise} in Appendix~\ref{sec:additionaltables} compare the DRCCP and C-DRCCP models over a sequence of Wasserstein radii. Both risk measures indicate that larger training samples allow the C-DRCCP model to align protection more closely with the risk at the target query. Although the DRCCP model also benefits from additional data, it is independent of the contextual information (query) and, therefore, cannot simultaneously reduce conservativeness at the low query and increase protection at the high query. In contrast, as $N$ grows, the risk disadvantage of the C-DRCCP model at the low query narrows while its cost generally remains lower, and its risk advantage at the high query becomes more consistent. These patterns indicate a more effective allocation of protection across contexts rather than uniform dominance of one model over the other.

\subsection{Anesthesiologist Deployment Case Study}\label{sec:staff-deployment}

\citet{rath2026staff} study anesthesiologist planning across a healthcare system with multiple hospital locations. Six weeks before surgery, each available physician is assigned to a regular shift, placed in a systemwide on-call pool, or left unassigned. Three days before surgery, location $s$ reports a forecast $F_{s,t}$ of the number of anesthesiologists required on date $t$. Physicians in the pool can then be deployed to specific locations. We study this three-day decision while treating the regular assignments and pool as fixed parameters.

In this case study, we investigate whether forecasts should determine not only the predicted workload, but also the amount of protection placed around that prediction. The DRCCP model protects against uncertainty without conditioning on the forecasts. In contrast, the C-DRCCP model can direct protection toward dates whose forecasts identify an adverse workload regime.

\subsubsection{C-DRCCP Deployment Model}

Let $S$ be the set of hospital locations and $I_t$ the set of physicians available on date $t$. Let $B_s$ denote the uncertain workload at location $s$. Parameter $v_{is}$ indicates whether physician $i$ is qualified at location $s$. The upstream model in~\citet{rath2026staff} determines the staffing plan $P_t=(\bar{x}_t,\bar{p}_t)$, where $\bar{x}_{is,t}$ gives the regular assignment and $\bar{p}_{i,t}$ indicates membership in the on-call pool. 

The binary decision $q_{is,t}$ equals one if on-call physician $i$ is activated and deployed to location $s$. Qualification and availability require $q_{is,t}\leq \bar{p}_{i,t}v_{is}$ and $\sum_{s\in S} q_{is,t}\leq\bar{p}_{i,t}$, respectively. Define the numbers of regularly assigned and activated on-call physicians at location $s$ by $R_{s,t}=\sum_{i\in I_t}\bar{x}_{is,t}$ and $Q_{s,t}(q)=\sum_{i\in I_t}q_{is,t}$. A regular physician supplies $t_r=480$ minutes, whereas an activated on-call physician supplies $t_q=360$ minutes. The physical capacity at location $s$ is $C_{s,t}(q)=t_rR_{s,t}+t_qQ_{s,t}(q)$. We choose an allowance $u_{s,t}\geq0$ for workload beyond physical capacity, which is measured in workload minutes and selected jointly with the deployment decision three days before surgery. Operationally, it may represent extended shifts, temporary or float staff, rescheduling, or the cost of delayed service. We value the allowance at an overtime rate of \$375 per hour used by~\citet{rath2026staff}. A pool physician costs \$500 if activated and \$300 otherwise. After dropping the fixed pool cost, the C-DRCCP deployment model is
\begin{subequations}\label{eq:staff-deployment-model}
\begin{align}
\min_{q,u}\quad
&200\sum_{i\in I_t}\sum_{s\in S}q_{is,t}+6.25\sum_{s\in S}u_{s,t}\\
\mathrm{s.t.}\quad
&\sup_{\mathbb P\in\mathcal P_{\theta_t,\varepsilon_t}(\mathcal N_t)}
\mathbb P\left(
\bigcup_{s\in S}\{B_s>C_{s,t}(q)+u_{s,t}\}
\,\middle|\,X\in\mathcal N_t
\right)\leq0.10\\
&q_{is,t}\leq\bar p_{i,t}v_{is} &&\forall i\in I_t,\ s\in S\\
&\sum_{s\in S}q_{is,t}\leq\bar p_{i,t} &&\forall i\in I_t\\
&q_{is,t}\in\{0,1\}\quad\forall i\in I_t,\ s\in S,
\qquad 0\leq u_{s,t}\leq U_{s,t}\quad\forall s\in S,
\end{align}
\end{subequations}
where the upper bounds $U_{s,t}$ are computed from the training workloads and are large enough to preserve an optimal solution. Model~\eqref{eq:staff-deployment-model} is solved with the exact MIP formulation in Section~\ref{sec:mip}.

\subsubsection{Benchmark Models and Evaluation Design}

The benchmark models address two questions. Comparing the SAA and DRCCP models, and the C-SAA and C-DRCCP models, assesses whether distributional robustness improves reliability. Comparing the DRCCP and C-DRCCP models assesses whether contextual information helps place that protection. The Point model provides a deterministic benchmark.

To compare uncertainty models on the same decision problem, we adapt the quantile-regression (QR) model of~\citet{rath2026staff} to the common $(q,u)$ formulation. Following their equations~(24)--(26), site and weekday fixed effects and a common forecast slope are fitted using training dates only. Let $B^U_{s,t}(F_{s,t};\rho)$ denote the upper endpoint of the fitted central $\rho = 0.99$ interval. The adapted QR model solves the problem
\begin{align*}
\min_{q,u}\left\{200\sum_{i\in I_t} \sum_{s\in S}q_{is,t}+6.25\sum_{s\in S}u_{s,t} \;\middle|\; B^U_{s,t}(F_{s,t};\rho)\leq C_{s,t}(q)+u_{s,t} \quad \forall s\in S\right\}
\end{align*}
together with the deployment constraints in~\eqref{eq:staff-deployment-model}. The adapted QR benchmark uses conditional prediction intervals to set workload allowances under the same deployment objective. The adapted QR, DRCCP, and C-DRCCP models differ in their treatment of uncertainty.

We use the real-world data provided by~\citet{rath2026staff}: 10 hospital locations, 159 training days, and 52 test days. The forecast vector $F_t=(F_{s,t})_{s\in S}$ provides the contextual information, and the outcome vector contains concurrency-adjusted billed minutes. The distance between two observations sums the absolute differences in their forecast contexts and workloads across the 10 locations. We use 127 of the training dates for estimation and 32 for validation, and freeze the same first-stage staffing plan across all development candidates. For the C-DRCCP model, the development experiments vary the context neighborhood, the lower bound on its probability, and the Wasserstein radius, and record validation cost and failure frequency. We use the C-DRCCP and DRCCP parameter settings from these development experiments. Each model is fitted using all 159 training dates and evaluated on the same 52 test dates.

We construct a stress test in which high forecasts correspond to higher workloads. We select the two locations with the largest average forecasts in the training data, where a date is classified as high-demand when the forecasts at both locations exceed their training medians. On these dates, we increase the workload at every location by $0.75$ of its training standard deviation. Forecasts and staffing inputs remain unchanged. For this experiment, we repeat estimation and validation using the modified training workloads, without using the test outcomes, and keep the baseline staffing plans fixed.

\paragraph{Evaluation metrics.}
A joint allowance failure occurs on a test date if workload at any location exceeds physical capacity plus the allowance chosen before the workload is observed. We report the number of failure dates and the one-sided $95\%$ upper confidence bound on the test failure probability. Excess is the total workload above physical capacity plus allowance, summed over all locations and the 52 test dates. Planning cost is the objective value of model~\eqref{eq:staff-deployment-model}, whereas operating cost uses the called, uncalled, overtime, and idle rates in the reference application. Both costs are reported as daily averages.

\begin{samepage}
\subsubsection{Operational Value of Context} We examine the value of robustness and contextual information.

\paragraph{The value of robustness.}
Table~\ref{tab:staff-main-results} compares test reliability and cost.
On the original released data, the Point, C-SAA, and SAA models incur 52, 18, and 9 joint allowance failures, respectively. Relative to these nonrobust benchmarks, the C-DRCCP and DRCCP models reduce the corresponding failure counts from 18 to 3 and from 9 to 2. These results indicate that localization alone does not provide reliable protection. Because the number of observations relevant to a particular forecast context is limited, the local empirical distribution remains vulnerable to sampling error. Distributional robustness compensates for this uncertainty, although it requires a larger ex ante planning margin. The adapted QR model has four failure dates at a lower planning cost than either robust model. The increase in realized operating cost is much smaller than the increase in planning cost. On the original data, the operating costs of C-SAA, SAA, C-DRCCP, DRCCP, and adapted QR lie within a range of less than $1.4\%$.
\end{samepage}

\begin{table}[htbp]
\centering
\caption{Reliability and cost on the test data under the common deployment objective.}
\label{tab:staff-main-results}
\setlength{\tabcolsep}{4.5pt}
\resizebox{.95\linewidth}{!}{%
\begin{tabular}{l|rrrrr|rrrrr}
\toprule[1pt]\midrule
& \multicolumn{5}{c|}{\emph{Panel A. Original released data}}
& \multicolumn{5}{c}{\emph{Panel B. Controlled stress}}\\
\cmidrule(lr){2-6}\cmidrule(lr){7-11}
Model
& \shortstack{Fail.\\/52} & \shortstack{Upper bound\\(\%)} & \shortstack{Excess\\(min)} & \shortstack{Planning\\(\$/day)} & \shortstack{Operating\\(\$/day)}
& \shortstack{Fail.\\/52} & \shortstack{Upper bound\\(\%)} & \shortstack{Excess\\(min)} & \shortstack{Planning\\(\$/day)} & \shortstack{Operating\\(\$/day)}\\
\midrule
Point      & 52 & 100.0 & 57,342 &  1,496 & 20,323 & 52 & 100.0 & 75,528 &   2,734 & 21,722\\
C-SAA      & 18 &  46.9 &  8,698 & 39,533 & 21,106 & 18 &  46.9 &  8,578 &  49,297 & 22,188\\
SAA        &  9 &  28.3 &  4,170 & 51,285 & 20,974 &  9 &  28.3 &  5,182 &  64,460 & 22,242\\
C-DRCCP    &  3 &  14.2 &    818 & 91,517 & 21,128 &  1 &   8.8 &    438 & 113,361 & 22,211\\
DRCCP      &  2 &  11.6 &    603 & 87,730 & 20,941 &  2 &  11.6 &    787 & 101,077 & 22,299\\
Adapted QR & 4  &  16.7 &  1,434 & 86,769 & 21,232 &  5 &  19.2 &  1,908 &  94,347 & 22,400\\
\midrule\bottomrule[1pt]
\end{tabular}%
}
\end{table}

\paragraph{The value of contextual information.}
Under the controlled stress, the C-DRCCP model has one joint allowance failure, compared with two under the DRCCP model and five under the adapted QR model. The single failure under the C-DRCCP model is contained in both comparison failure sets. Its cumulative allowance excess across sites and dates is 438 minutes, $44.3\%$ below that of the DRCCP model and $77.0\%$ below that of the adapted QR model. Its upper confidence bound is $8.8\%$. The C-DRCCP model uses a larger ex ante planning margin, with planning costs $12.2\%$ and $20.2\%$ above those of the DRCCP and adapted QR models, respectively. Yet it has the lowest realized operating cost of the three policies.

\begin{figure}[htbp]
  \centering
  \caption{Contextual reallocation of allowance and test failures by forecast regime.}
  \label{fig:staff-context-value}
  \includegraphics[width=0.8\linewidth]{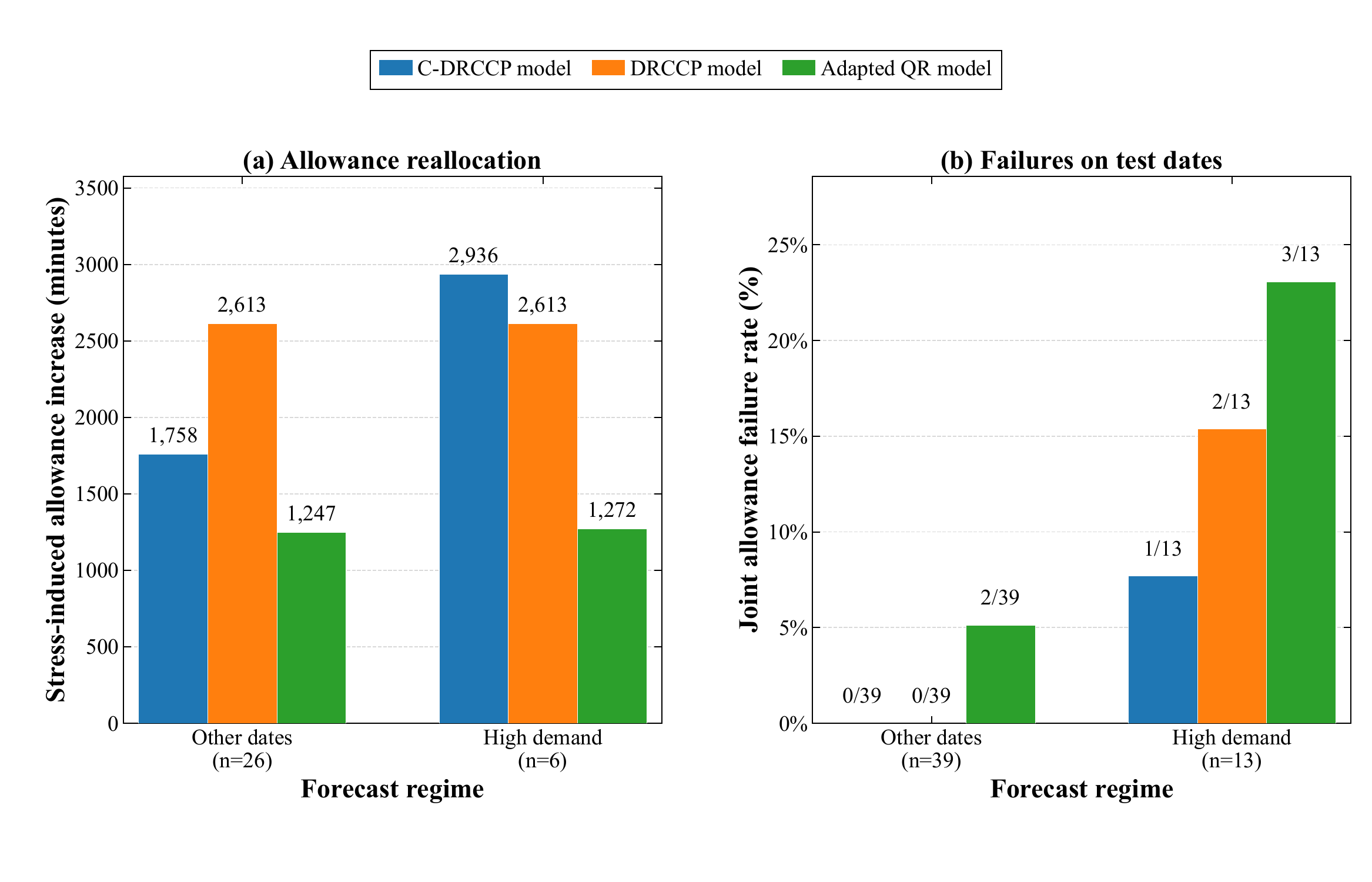}
\end{figure}

Figure~\ref{fig:staff-context-value} illustrates how contextual information changes the allocation of protection. The C-DRCCP model increases protection more on high-demand dates, increasing allowance by 2{,}936 minutes on the 6 high-demand validation dates and by 1{,}758 minutes on the remaining 26 dates. In contrast, the DRCCP model increases allowance by 2{,}613 minutes in both regimes, while the adapted QR model changes by similar amounts across them. The test failures follow this allocation pattern: among the high-demand dates, the C-DRCCP, DRCCP, and adapted QR models fail on one, two, and three dates, respectively, whereas the C-DRCCP model has no failures on the other 39 dates. 

\section{Conclusion}\label{sec:conclusion}
We develop exact reformulations for contextual distributionally robust chance constraints over joint Wasserstein balls. By tracking the conditioning mass and the portion assigned to failure, we can obtain an exact MIP reformulation for affine polyhedral safety systems. We further show that robust feasibility depends on whether the minimum transportation cost to reach the conditional risk boundary meets or exceeds the transportation budget. This characterization yields an exact inequality reformulation and explains why normalizing contextual mass allocations loses information about the remaining transportation budget. Preserving this information gives strengthened quantile cuts that strictly dominate their standard counterparts and yields tighter valid coefficients in the MIP formulation. The ordering of contextual costs also yields rank inequalities that combine binary restrictions across allocations. Under suitable coverage conditions, feasible solutions control the true conditional risk, and local regularity extends this guarantee to the target context.

The numerical results show that these valid inequalities substantially tighten the root relaxation and reduce MIP solution times. In the capacitated transportation problem, the \texttt{ALL} formulation reduces the mean root gap from $89.48\%$ to $0.69\%$. Its computational benefit is clearest at $N=1{,}000$: it solves all 12 instances with a median final MIP time of $39.67$ seconds, whereas every run of the basic formulation reaches the time limit of $3{,}600$ seconds. The out-of-sample experiments also show that the value of context depends on the amount of data. When local observations are sparse, the resulting estimation error can make the DRCCP model a safer baseline. With more data, the C-DRCCP model can estimate contextual differences more reliably and choose a level of protection that better matches the risk associated with the observed context, thereby improving the balance between cost and reliability.


\bibliographystyle{informs2014}
\bibliography{ref}
\newpage
\ECSwitch
\begin{APPENDICES}

\renewcommand*{\theHsection}{appendix.\arabic{section}}
\renewcommand*{\theHequation}{appendix.\arabic{equation}}
\renewcommand*{\theHtable}{appendix.\arabic{table}}
\renewcommand*{\theHfigure}{appendix.\arabic{figure}}
\renewcommand*{\theHtheorem}{appendix.\arabic{theorem}}
\def\theHproposition{appendix.\arabic{proposition}}
\def\theHcorollary{appendix.\arabic{corollary}}
\def\theHlemma{appendix.\arabic{lemma}}
\def\theHremark{appendix.\arabic{remark}}
\def\theHexample{appendix.\arabic{example}}

\section{Relations with the Unconditional DRCCP Model}\label{appendix:model-relations}
We maintain Assumption~\ref{ass:transport} and condition~\eqref{eq:strict-radius}, and abbreviate $\theta_{\min}(x_0,\gamma,\varepsilon)$ as $\theta_{\min}$. We use the polytope of mass allocations at the conditional risk boundary $\mathcal{C}_{\varepsilon,\epsilon}$ and the value function $\vartheta$ defined in~\eqref{eq:critical-allocation-polytope}--\eqref{eq:critical-radius-value}. We also use the allocation polytope $\mathcal{W}_\theta$ and residual budget $B(w)=\theta-K_0-\bar\kappa^\top w$ from~\eqref{eq:context-allocation-polytope}. Thus, $B(w)>0$ means that allocation $w$ can be formed, with a transportation budget remaining for changing the outcomes.

We first define an unconditional benchmark on the outcome space. Let $W_{\mathcal{Y}}$ denote the type-1 Wasserstein distance induced by $d_{\mathcal{Y}}$. For a probability vector $p\in\R_+^N$ with $\mathbf{1}^\top p=1$, let $\widehat{\mathbb{Q}}_p\coloneqq\sum_{i\in[N]}p_i\delta_{\hat{y}_i}$. For a Wasserstein radius $\rho>0$ and risk level $\alpha\in(0,1)$, define
\begin{equation*}
\mathcal{Z}_{\rm DR}(p,\rho,\alpha)
\coloneqq
\left\{
z\in\mathcal{Z}\;\middle|\;
\sup_{\substack{\mathbb{Q}: W_{\mathcal{Y}}(\mathbb{Q},\widehat{\mathbb{Q}}_p)\leq\rho}}
\mathbb{Q}\left(\mathcal{S}^c(z)\right)\leq\alpha
\right\},
\end{equation*}
Its critical transportation cost is
\begin{equation}\label{eq:weighted-critical-cost}
\varphi_{\alpha,p}(d)
\coloneqq
\min\left\{
d^\top u\;\middle|\;
0\leq u\leq p,\ \mathbf{1}^\top u=\alpha
\right\}.
\end{equation}
For the uniform vector $p=\mathbf{1}/N$, we write $\varphi_\alpha$ and $\mathcal{Z}_{\rm DR}(\rho,\alpha)$.

Lemma~\ref{lem:unconditional-critical-cost} is the weighted finite-support version of the standard critical transportation characterization for Wasserstein DRCCP models~\citep{xie2021drcc,chen2022data}. We state it explicitly because both the empirical weights and the radius will depend on the context allocation.

\begin{lemma}\label{lem:unconditional-critical-cost}
For every probability vector $p$, risk level $\alpha\in(0,1)$, radius $\rho>0$, and decision $z\in\mathcal{Z}$,
\begin{equation}\label{eq:unconditional-critical-equivalence}
z\in\mathcal{Z}_{\rm DR}(p,\rho,\alpha)
\quad\Longleftrightarrow\quad
\varphi_{\alpha,p}(\Delta(z))\geq\rho.
\end{equation}
Moreover, let $\widehat{\mathbb{P}}_p=\sum_{i\in[N]}p_i\delta_{(\hat{x}_i,\hat{y}_i)}$. Under the additive product metric, the outcome marginals of the distributions in $\{\mathbb{P}\;|\;W(\mathbb{P},\widehat{\mathbb{P}}_p)\leq\rho\}$ form exactly the ball $\{\mathbb{Q}\;|\;W_{\mathcal{Y}}(\mathbb{Q},\widehat{\mathbb{Q}}_p)\leq\rho\}$.
\end{lemma}

The outcome benchmark in Lemma~\ref{lem:unconditional-critical-cost} is defined separately because the formal choice $\mathcal{N}=\mathcal{X}$ violates the finite distance to $\mathcal{N}^c$ required by Assumption~\ref{ass:transport}. 

Let $\mathcal{Z}_{\rm C}(\theta,\varepsilon,\epsilon)$ denote the feasible region of the C-DRCCP model~\eqref{opt:drc3p}. For every $w\in\mathcal{W}_\theta$ with $B(w)>0$, define $p_i(w)\coloneqq\frac{w_i}{\mathbf{1}^\top w}$ and $\rho(w)\coloneqq\frac{B(w)}{\mathbf{1}^\top w}$. The vector $p(w)$ gives the empirical outcome weights induced by the mass allocation $w$, while $\rho(w)$ is the remaining transportation budget per unit of conditioning mass.

\begin{theorem}\label{thm:weighted-drccp-intersection}
The C-DRCCP feasible region satisfies
\begin{equation}\label{eq:weighted-drccp-intersection}
\mathcal{Z}_{\rm C}(\theta,\varepsilon,\epsilon)
=
\bigcap_{\substack{w\in\mathcal{W}_\theta\\ B(w)>0}}
\mathcal{Z}_{\rm DR}\left(p(w),\rho(w),\epsilon\right).
\end{equation}
\end{theorem}

Theorem~\ref{thm:weighted-drccp-intersection} gives the exact rule for comparing the two models. In general, the C-DRCCP model cannot be obtained from the unconditional benchmark by rescaling its radius alone. It imposes a family of weighted DRCCP constraints in which the empirical probabilities $p(w)$ and radius $\rho(w)$ are determined by the same allocation. Assigning more mass to empirical sources whose contexts are inexpensive to place in the neighborhood changes both quantities at once.

Define the cost of transporting all empirical mass into the conditioning neighborhood as
\[K_{\rm in}\coloneqq K_0+\frac{1}{N}\bar\kappa^\top\mathbf{1}=\frac{1}{N}\sum_{i\in[N]}\kappa_i^{\rm in}.\]

\begin{theorem}\label{thm:unconditional-drccp-bounds}
For every finite $d\in\R_+^N$,
\begin{equation}\label{eq:unconditional-critical-bounds}
\theta_{\min}+\varphi_{\epsilon\varepsilon}(d)
\leq
\vartheta(d)
\leq
K_{\rm in}+\varphi_\epsilon(d).
\end{equation}
Consequently,
\begin{equation}\label{eq:unconditional-sufficient-region}
\mathcal{Z}_{\rm DR}
\left(\theta-\theta_{\min},\epsilon\varepsilon\right)
\subseteq
\mathcal{Z}_{\rm C}(\theta,\varepsilon,\epsilon).
\end{equation}
If $\theta>K_{\rm in}$, then
\begin{equation}\label{eq:unconditional-necessary-region}
\mathcal{Z}_{\rm C}(\theta,\varepsilon,\epsilon)
\subseteq
\mathcal{Z}_{\rm DR}
\left(\theta-K_{\rm in},\epsilon\right).
\end{equation}
Both inclusions can be strict.
\end{theorem}

The two bounds serve different purposes. The first inclusion gives a sufficient unconditional benchmark: it reserves the minimum context cost and protects the smallest possible failure mass $\epsilon\varepsilon$. The second gives a necessary unconditional benchmark obtained by transporting all mass into the neighborhood. When $\varepsilon<1$, these bounds therefore use the distinct risk levels $\epsilon\varepsilon$ and $\epsilon$.

\begin{corollary}\label{cor:exact-unconditional-cases}
The following exact relations hold.
\begin{enumerate}[label=(\roman*),leftmargin=2em]
\item If $\varepsilon=1$, then $w=\mathbf{1}/N$ is forced, $\theta_{\min}=K_{\rm in}$, and $\mathcal{Z}_{\rm C}(\theta,1,\epsilon)=\mathcal{Z}_{\rm DR}(\theta-K_{\rm in},\epsilon)$.
\item If $\bar\kappa_i=\bar\kappa\geq0$ for every $i\in[N]$, then $\theta_{\min}=K_0+\bar\kappa\varepsilon$ and $\mathcal{Z}_{\rm C}(\theta,\varepsilon,\epsilon)=\mathcal{Z}_{\rm DR}\left(\theta-\theta_{\min},\epsilon\varepsilon\right)$.
\end{enumerate}
\end{corollary}

The sign condition in part~(ii) cannot be omitted. If the common coefficient $\bar\kappa$ is negative, assigning more mass to the neighborhood reduces contextual transportation cost. The optimal conditioning mass can then depend on the outcome distances, so the transformation of the risk level in part~(ii) need not remain valid.

Part~(i) also clarifies two distinct special cases. If the conditioning event is defined separately as the global event $\mathcal{N}=\mathcal{X}$, then its probability is one under every distribution. Lemma~\ref{lem:unconditional-critical-cost} implies that the resulting joint Wasserstein model is exactly $\mathcal{Z}_{\rm DR}(\theta,\epsilon)$, independently of $\varepsilon$. In contrast, suppose that a proper conditioning neighborhood has $K_0=0$ and $\bar\kappa=0$, but its mass is only required to be at least $\varepsilon$. Part~(ii) gives $\mathcal{Z}_{\rm C}(\theta,\varepsilon,\epsilon)=\mathcal{Z}_{\rm DR}(\theta,\epsilon\varepsilon)$. Thus, zero contextual transportation cost alone does not establish equivalence with the unconditional model. For this equivalence to hold in general, the conditioning event must also have probability one under every distribution in the ambiguity set.

\section{Additional Numerical Results}\label{sec:additionaltables}

\paragraph{Performance across scales.}
Table~\ref{tab:formal-scale-performance} reports the median root gap, maximum final gap, median solving time, and median search-tree size for each network, training size, and formulation. For $N\leq100$, the median final solving times remain below 1 second despite the large differences in root gaps. At $N=500$, the baseline \texttt{MIP} formulation has median root gaps of $92.89\%$--$96.28\%$, median final solving times of $5.14$--$136.90$ seconds, and median trees of $1{,}591$--$6{,}019$ nodes. \texttt{ALL} has corresponding ranges of $1.20\%$--$2.51\%$, $1.32$--$9.60$ seconds, and $60$--$128.5$ nodes, respectively.

The median root gaps under \texttt{SQC} are $5.58\%$--$7.79\%$, compared with $8.16\%$--$8.68\%$ under \texttt{QC}, $11.78\%$--$13.43\%$ under \texttt{FMC}, and $18.40\%$--$26.32\%$ under the formulations that use only probability cuts. For $(F,D)=(30,100)$ and $N=500$, \texttt{SQC} reduces the baseline median root gap from $92.89\%$ to $6.05\%$, the median final solving time from $134.96$ to $7.79$ seconds, and the median tree from $4{,}886$ to $504.5$ nodes. \texttt{ALL} has corresponding values of $1.43\%$, $5.67$ seconds, and $113$ nodes. On these same 12 instances, \texttt{RANK} has medians of $14.36\%$, $35.58$ seconds, and $1{,}277.5$ nodes. Thus, the strengthened quantile cuts provide most of the initial tightening, while the complementary inequalities in \texttt{ALL} further reduce the search effort.

\setlength{\rotFPtop}{0pt plus 1fil}\setlength{\rotFPbot}{0pt plus 1fil}
\begin{sidewaystable}[htbp]
\makeatletter\let\@makecaption\@maketablecaption\makeatother
\centering
\caption{Root-relaxation and final MIP performance by network and training-sample size.}
\label{tab:formal-scale-performance}
\setlength{\tabcolsep}{2.5pt}
\resizebox{.9\linewidth}{!}{%
\begin{tabular}{ll|rrrr|rrrr|rrrr|rrrr|rrrr}\toprule[1pt]\midrule
 & & \multicolumn{4}{c}{\texttt{MIP}} & \multicolumn{4}{c}{\texttt{PC}} & \multicolumn{4}{c}{\texttt{QC}} & \multicolumn{4}{c}{\texttt{FMC}} & \multicolumn{4}{c}{\texttt{SQC}}\\
$(F,D)$ & $N$ & $g^{\rm LP}$ & $g^{\rm final}$ & $T^{\rm MIP}$ & $B$ & $g^{\rm LP}$ & $g^{\rm final}$ & $T^{\rm MIP}$ & $B$ & $g^{\rm LP}$ & $g^{\rm final}$ & $T^{\rm MIP}$ & $B$ & $g^{\rm LP}$ & $g^{\rm final}$ & $T^{\rm MIP}$ & $B$ & $g^{\rm LP}$ & $g^{\rm final}$ & $T^{\rm MIP}$ & $B$\\\midrule
$(5,20)$ & 50 & 88.88 & 0.00 & 0.03 & 5 & 18.26 & 0.00 & 0.02 & 1 & 0.00 & 0.00 & 0.01 & 1 & 7.94 & 0.00 & 0.03 & 1 & 0.00 & 0.00 & 0.01 & 1\\
 & 100 & 96.10 & 0.00 & 0.08 & 10 & 23.69 & 0.00 & 0.11 & 18 & 0.11 & 0.00 & 0.04 & 1 & 4.23 & 0.00 & 0.06 & 6 & 0.09 & 0.00 & 0.04 & 1\\
 & 500 & 96.28 & 0.00 & 5.14 & 1{,}591 & 23.39 & 0.00 & 3.57 & 809 & 8.68 & 0.00 & 1.99 & 314 & 13.43 & 0.00 & 3.64 & 612.5 & 7.79 & 0.00 & 1.52 & 212.5\\
\addlinespace[3pt]
$(10,50)$ & 50 & 90.68 & 0.00 & 0.06 & 3 & 15.62 & 0.00 & 0.06 & 1 & 0.00 & 0.00 & 0.05 & 1 & 6.65 & 0.00 & 0.09 & 1 & 0.00 & 0.00 & 0.05 & 1\\
 & 100 & 95.12 & 0.00 & 0.29 & 47.5 & 15.04 & 0.00 & 0.30 & 29 & 2.34 & 0.00 & 0.24 & 3.5 & 6.64 & 0.00 & 0.32 & 11.5 & 1.72 & 0.00 & 0.25 & 3\\
 & 500 & 94.21 & 0.00 & 58.16 & 4{,}929 & 18.40 & 0.00 & 27.70 & 2{,}144 & 8.27 & 0.00 & 21.50 & 1{,}407.5 & 11.78 & 0.00 & 26.77 & 1{,}667.5 & 6.30 & 0.00 & 11.51 & 715.5\\
\addlinespace[3pt]
$(15,100)$ & 50 & 93.27 & 0.00 & 0.12 & 1.5 & 19.90 & 0.00 & 0.09 & 1 & 0.00 & 0.00 & 0.17 & 1 & 1.55 & 0.00 & 0.20 & 1 & 0.00 & 0.00 & 0.17 & 1\\
 & 100 & 93.51 & 0.00 & 0.77 & 44 & 17.78 & 0.00 & 0.70 & 26 & 0.01 & 0.00 & 0.48 & 1.5 & 6.08 & 0.00 & 0.68 & 1.5 & 0.01 & 0.00 & 0.45 & 1\\
 & 500 & 93.20 & 2.93 & 136.90 & 6{,}019 & 20.46 & 0.00 & 62.07 & 1{,}957 & 8.16 & 0.63 & 20.62 & 448 & 13.06 & 0.00 & 79.70 & 2{,}100 & 5.58 & 0.00 & 14.12 & 270.5\\
\addlinespace[3pt]
$(30,100)$ & 50 & 92.93 & 0.00 & 0.14 & 1 & 18.41 & 0.00 & 0.14 & 1 & 0.00 & 0.00 & 0.23 & 1 & 1.62 & 0.00 & 0.24 & 1 & 0.00 & 0.00 & 0.23 & 1\\
 & 100 & 93.23 & 0.00 & 0.93 & 52.5 & 17.10 & 0.01 & 0.97 & 33 & 0.01 & 0.00 & 0.45 & 1 & 6.99 & 0.00 & 0.58 & 2.5 & 0.01 & 0.00 & 0.45 & 1\\
 & 500 & 92.89 & 2.98 & 134.96 & 4{,}886 & 20.78 & 0.01 & 57.57 & 1{,}307 & 8.42 & 0.01 & 10.77 & 753 & 12.84 & 0.01 & 27.33 & 2{,}229 & 6.05 & 0.01 & 7.79 & 504.5\\
\midrule\end{tabular}}
\par\medskip
\resizebox{.9\linewidth}{!}{%
\begin{tabular}{ll|rrrr|rrrr|rrrr|rrrr|rrrr}\midrule
 & & \multicolumn{4}{c}{\texttt{SQC+MIX}} & \multicolumn{4}{c}{\texttt{SP}} & \multicolumn{4}{c}{\texttt{FAH}} & \multicolumn{4}{c}{\texttt{RANK}} & \multicolumn{4}{c}{\texttt{ALL}}\\
$(F,D)$ & $N$ & $g^{\rm LP}$ & $g^{\rm final}$ & $T^{\rm MIP}$ & $B$ & $g^{\rm LP}$ & $g^{\rm final}$ & $T^{\rm MIP}$ & $B$ & $g^{\rm LP}$ & $g^{\rm final}$ & $T^{\rm MIP}$ & $B$ & $g^{\rm LP}$ & $g^{\rm final}$ & $T^{\rm MIP}$ & $B$ & $g^{\rm LP}$ & $g^{\rm final}$ & $T^{\rm MIP}$ & $B$\\\midrule
$(5,20)$ & 50 & 0.00 & 0.00 & 0.01 & 1 & 12.65 & 0.00 & 0.03 & 1 & 12.65 & 0.00 & 0.03 & 1 & 0.00 & 0.00 & 0.01 & 1 & 0.00 & 0.00 & 0.02 & 1\\
 & 100 & 0.09 & 0.00 & 0.04 & 1 & 20.24 & 0.00 & 0.07 & 7.5 & 15.49 & 0.00 & 0.07 & 7.5 & 0.11 & 0.00 & 0.02 & 1 & 0.02 & 0.00 & 0.04 & 1\\
 & 500 & 7.79 & 0.00 & 1.77 & 245 & 26.32 & 0.00 & 3.76 & 880 & 26.32 & 0.00 & 3.68 & 880 & 18.44 & 0.00 & 3.11 & 447 & 2.51 & 0.00 & 1.32 & 60\\
\addlinespace[3pt]
$(10,50)$ & 50 & 0.00 & 0.00 & 0.06 & 1 & 11.23 & 0.00 & 0.04 & 1 & 11.23 & 0.00 & 0.04 & 1 & 0.00 & 0.00 & 0.03 & 1 & 0.00 & 0.00 & 0.05 & 1\\
 & 100 & 1.72 & 0.00 & 0.25 & 3 & 17.40 & 0.00 & 0.23 & 22 & 17.40 & 0.00 & 0.23 & 22 & 4.09 & 0.00 & 0.13 & 1 & 0.04 & 0.00 & 0.17 & 1\\
 & 500 & 6.30 & 0.00 & 12.25 & 696.5 & 20.50 & 0.00 & 18.37 & 2{,}175.5 & 20.50 & 0.00 & 17.84 & 2{,}175.5 & 13.42 & 0.00 & 13.58 & 1{,}567 & 1.92 & 0.00 & 4.47 & 128.5\\
\addlinespace[3pt]
$(15,100)$ & 50 & 0.00 & 0.00 & 0.17 & 1 & 1.88 & 0.00 & 0.09 & 1 & 1.88 & 0.00 & 0.09 & 1 & 0.00 & 0.00 & 0.08 & 1 & 0.00 & 0.00 & 0.15 & 1\\
 & 100 & 0.01 & 0.00 & 0.49 & 1 & 17.34 & 0.00 & 0.92 & 15 & 17.34 & 0.00 & 0.92 & 15 & 0.01 & 0.00 & 0.23 & 1 & 0.01 & 0.00 & 0.46 & 1\\
 & 500 & 4.73 & 0.00 & 15.27 & 273.5 & 22.67 & 0.00 & 49.58 & 1{,}831 & 22.66 & 0.00 & 40.67 & 1{,}760 & 13.95 & 0.00 & 24.60 & 1{,}294 & 1.20 & 0.00 & 9.60 & 93\\
\addlinespace[3pt]
$(30,100)$ & 50 & 0.00 & 0.00 & 0.23 & 1 & 1.96 & 0.01 & 0.13 & 1 & 1.96 & 0.01 & 0.13 & 1 & 0.00 & 0.00 & 0.13 & 1 & 0.00 & 0.00 & 0.25 & 1\\
 & 100 & 0.01 & 0.01 & 0.59 & 1 & 16.92 & 0.00 & 0.76 & 17 & 16.92 & 0.00 & 0.75 & 17 & 0.01 & 0.00 & 0.32 & 1 & 0.01 & 0.01 & 0.59 & 1\\
 & 500 & 5.48 & 0.01 & 8.53 & 449.5 & 22.54 & 0.01 & 63.85 & 1{,}937.5 & 22.53 & 0.01 & 63.47 & 2{,}052 & 14.36 & 0.01 & 35.58 & 1{,}277.5 & 1.43 & 0.01 & 5.67 & 113\\
\midrule\bottomrule[1pt]\end{tabular}}
\par\smallskip\footnotesize Note. Each cell summarizes 12 runs (three queries and four radius levels): median $g^{\rm LP}$, maximum $g^{\rm final}$, median $T^{\rm MIP}$, and median $B$. 
\end{sidewaystable}

\paragraph{Fixed-radius comparison.}
Tables~\ref{tab:fixed-raw-radius-neighborhood} and~\ref{tab:fixed-raw-radius-pointwise} provide the detailed sensitivity analysis summarized in Section~\ref{sec:oos-performance}. The models use different ground spaces, so equal numerical radii need not define ambiguity sets of equal strength. The comparison instead shows how each model responds to the same radius grid. In both tables, entries report mean cost followed by mean violation probability for the indicated risk measure. The $N=200$ and $N=300$ rows use one nested replication, and a dash denotes an unavailable result.

\begingroup
\setlength{\rotFPtop}{0pt plus 1fil}
\setlength{\rotFPbot}{0pt plus 1fil}
\begin{sidewaystable}[htbp]
\makeatletter
\let\@makecaption\@maketablecaption
\makeatother
\centering
\caption{Cost and neighborhood risk at fixed numerical Wasserstein radii.}
\label{tab:fixed-raw-radius-neighborhood}
\scriptsize
\renewcommand{\arraystretch}{1.4}
\setlength{\tabcolsep}{2.5pt}
\begin{tabular}{ccl|cc|cc|cc|cc|cc}\toprule[1pt]\midrule
&&&\multicolumn{2}{c|}{$0.010$}&\multicolumn{2}{c|}{$0.0125$}&\multicolumn{2}{c|}{$0.015$}&\multicolumn{2}{c|}{$0.0175$}&\multicolumn{2}{c}{$0.020$}\\
$(F,D)$ & $N$ & Query & DRCCP & C-DRCCP & DRCCP & C-DRCCP & DRCCP & C-DRCCP & DRCCP & C-DRCCP & DRCCP & C-DRCCP\\
\midrule
\multirow{12}{*}{$(10,50)$} & \multirow{3}{*}{50} & Low & $882.5\,/\;0.050$ & $777.7\,/\;0.257$ & $888.4\,/\;0.046$ & $812.1\,/\;0.165$ & $894.0\,/\;0.044$ & $845.3\,/\;0.108$ & $899.0\,/\;0.040$ & $872.9\,/\;0.078$ & $903.9\,/\;0.037$ & $912.9\,/\;0.045$\\
 & & Central & $882.5\,/\;0.190$ & $857.0\,/\;0.246$ & $888.4\,/\;0.177$ & $901.6\,/\;0.154$ & $894.0\,/\;0.172$ & $943.0\,/\;0.105$ & $899.0\,/\;0.161$ & $979.8\,/\;0.064$ & $903.9\,/\;0.150$ & $1013.6\,/\;0.043$\\
 & & High & $882.5\,/\;0.537$ & $947.2\,/\;0.370$ & $888.4\,/\;0.519$ & $988.9\,/\;0.267$ & $894.0\,/\;0.509$ & $1015.8\,/\;0.214$ & $899.0\,/\;0.491$ & $1049.9\,/\;0.162$ & $903.9\,/\;0.474$ & $1080.9\,/\;0.125$\\
\cmidrule(lr){2-13}
 & \multirow{3}{*}{100} & Low & $947.2\,/\;0.014$ & $817.9\,/\;0.124$ & $953.9\,/\;0.013$ & $851.9\,/\;0.066$ & $960.4\,/\;0.011$ & $887.4\,/\;0.038$ & $966.5\,/\;0.010$ & $926.1\,/\;0.021$ & $972.4\,/\;0.009$ & $970.1\,/\;0.009$\\
 & & Central & $947.2\,/\;0.074$ & $896.6\,/\;0.146$ & $953.9\,/\;0.066$ & $942.5\,/\;0.082$ & $960.4\,/\;0.059$ & $997.3\,/\;0.040$ & $966.5\,/\;0.054$ & $1044.6\,/\;0.023$ & $972.4\,/\;0.050$ & $1084.2\,/\;0.013$\\
 & & High & $947.2\,/\;0.323$ & $1020.5\,/\;0.182$ & $953.9\,/\;0.303$ & $1066.4\,/\;0.119$ & $960.4\,/\;0.286$ & $1103.4\,/\;0.080$ & $966.5\,/\;0.271$ & $1142.4\,/\;0.055$ & $972.4\,/\;0.256$ & $1184.0\,/\;0.039$\\
\cmidrule(lr){2-13}
 & \multirow{3}{*}{200} & Low & $966.9\,/\;0.008$ & $841.6\,/\;0.054$ & $973.8\,/\;0.007$ & $886.2\,/\;0.023$ & $980.1\,/\;0.006$ & $946.6\,/\;0.009$ & $985.9\,/\;0.005$ & $992.7\,/\;0.004$ & $991.8\,/\;0.005$ & $1043.8\,/\;0.002$\\
 & & Central & $966.9\,/\;0.038$ & $925.6\,/\;0.077$ & $973.8\,/\;0.037$ & $992.1\,/\;0.035$ & $980.1\,/\;0.034$ & $1043.2\,/\;0.015$ & $985.9\,/\;0.031$ & $1079.8\,/\;0.008$ & $991.8\,/\;0.028$ & $1119.7\,/\;0.004$\\
 & & High & $966.9\,/\;0.236$ & $1064.3\,/\;0.090$ & $973.8\,/\;0.229$ & $1109.4\,/\;0.058$ & $980.1\,/\;0.219$ & $1148.2\,/\;0.038$ & $985.9\,/\;0.204$ & $1195.3\,/\;0.024$ & $991.8\,/\;0.192$ & $1279.0\,/\;0.015$\\
\cmidrule(lr){2-13}
 & \multirow{3}{*}{300} & Low & $989.4\,/\;0.003$ & $857.7\,/\;0.041$ & $997.4\,/\;0.003$ & $927.4\,/\;0.010$ & $1004.7\,/\;0.002$ & $979.6\,/\;0.004$ & $1011.6\,/\;0.002$ & $1024.7\,/\;0.002$ & $1018.0\,/\;0.002$ & $1080.3\,/\;0.001$\\
 & & Central & $989.4\,/\;0.023$ & $961.4\,/\;0.043$ & $997.4\,/\;0.020$ & $1017.7\,/\;0.016$ & $1004.7\,/\;0.018$ & $1073.9\,/\;0.007$ & $1011.6\,/\;0.016$ & $1112.2\,/\;0.003$ & $1018.0\,/\;0.015$ & $1152.9\,/\;0.002$\\
 & & High & $989.4\,/\;0.173$ & $1086.5\,/\;0.061$ & $997.4\,/\;0.165$ & $1142.9\,/\;0.038$ & $1004.7\,/\;0.150$ & $1181.3\,/\;0.026$ & $1011.6\,/\;0.145$ & $1223.4\,/\;0.017$ & $1018.0\,/\;0.137$ & --\\
\midrule
\multirow{12}{*}{$(15,100)$} & \multirow{3}{*}{50} & Low & $2522.3\,/\;0.051$ & $2093.3\,/\;0.330$ & $2539.6\,/\;0.046$ & $2234.8\,/\;0.202$ & $2556.7\,/\;0.040$ & $2351.6\,/\;0.128$ & $2573.4\,/\;0.035$ & $2475.0\,/\;0.085$ & $2589.2\,/\;0.033$ & $2577.7\,/\;0.055$\\
 & & Central & $2522.3\,/\;0.217$ & $2382.4\,/\;0.339$ & $2539.6\,/\;0.201$ & $2521.0\,/\;0.207$ & $2556.7\,/\;0.183$ & $2639.6\,/\;0.127$ & $2573.4\,/\;0.171$ & $2749.7\,/\;0.078$ & $2589.2\,/\;0.162$ & $2853.7\,/\;0.047$\\
 & & High & $2522.3\,/\;0.555$ & $2499.9\,/\;0.605$ & $2539.6\,/\;0.535$ & $2599.3\,/\;0.498$ & $2556.7\,/\;0.511$ & $2702.5\,/\;0.397$ & $2573.4\,/\;0.493$ & $2805.1\,/\;0.309$ & $2589.2\,/\;0.478$ & $2897.5\,/\;0.241$\\
\cmidrule(lr){2-13}
 & \multirow{3}{*}{100} & Low & $2670.3\,/\;0.013$ & $2247.6\,/\;0.165$ & $2691.3\,/\;0.012$ & $2382.3\,/\;0.093$ & $2710.1\,/\;0.011$ & $2502.4\,/\;0.048$ & $2727.1\,/\;0.010$ & $2628.3\,/\;0.025$ & $2743.4\,/\;0.009$ & $2748.6\,/\;0.013$\\
 & & Central & $2670.3\,/\;0.089$ & $2573.5\,/\;0.154$ & $2691.3\,/\;0.081$ & $2733.1\,/\;0.080$ & $2710.1\,/\;0.076$ & $2835.8\,/\;0.047$ & $2727.1\,/\;0.070$ & $2961.5\,/\;0.026$ & $2743.4\,/\;0.065$ & $3084.3\,/\;0.014$\\
 & & High & $2670.3\,/\;0.358$ & $2774.6\,/\;0.301$ & $2691.3\,/\;0.339$ & $2873.0\,/\;0.230$ & $2710.1\,/\;0.327$ & $2989.8\,/\;0.168$ & $2727.1\,/\;0.312$ & $3090.9\,/\;0.118$ & $2743.4\,/\;0.298$ & $3201.8\,/\;0.081$\\
\cmidrule(lr){2-13}
 & \multirow{3}{*}{200} & Low & $2891.1\,/\;0.002$ & $2530.2\,/\;0.032$ & $2912.7\,/\;0.002$ & $2697.5\,/\;0.012$ & $2934.2\,/\;0.001$ & $2812.7\,/\;0.004$ & $2953.5\,/\;0.001$ & $2952.9\,/\;0.002$ & $2972.7\,/\;0.001$ & $3110.5\,/\;0.001$\\
 & & Central & $2891.1\,/\;0.022$ & $2747.8\,/\;0.062$ & $2912.7\,/\;0.019$ & $2971.7\,/\;0.018$ & $2934.2\,/\;0.017$ & $3093.9\,/\;0.010$ & $2953.5\,/\;0.015$ & $3206.4\,/\;0.005$ & $2972.7\,/\;0.013$ & $3309.2\,/\;0.003$\\
 & & High & $2891.1\,/\;0.150$ & $3062.0\,/\;0.097$ & $2912.7\,/\;0.141$ & $3183.0\,/\;0.067$ & $2934.2\,/\;0.133$ & $3285.5\,/\;0.045$ & $2953.5\,/\;0.125$ & $3441.2\,/\;0.025$ & $2972.7\,/\;0.118$ & $3597.8\,/\;0.019$\\
\cmidrule(lr){2-13}
 & \multirow{3}{*}{300} & Low & $2923.6\,/\;0.002$ & $2545.2\,/\;0.023$ & $2946.6\,/\;0.002$ & $2698.3\,/\;0.010$ & $2968.8\,/\;0.002$ & $2805.4\,/\;0.005$ & $2988.1\,/\;0.002$ & $2928.4\,/\;0.002$ & $3007.4\,/\;0.001$ & $3112.7\,/\;0.001$\\
 & & Central & $2923.6\,/\;0.017$ & $2809.5\,/\;0.037$ & $2946.6\,/\;0.015$ & $3020.5\,/\;0.011$ & $2968.8\,/\;0.014$ & $3157.5\,/\;0.006$ & $2988.1\,/\;0.013$ & $3252.8\,/\;0.003$ & $3007.4\,/\;0.012$ & $3358.8\,/\;0.002$\\
 & & High & $2923.6\,/\;0.135$ & $3123.3\,/\;0.075$ & $2946.6\,/\;0.125$ & $3230.3\,/\;0.052$ & $2968.8\,/\;0.120$ & $3334.2\,/\;0.036$ & $2988.1\,/\;0.113$ & $3459.5\,/\;0.021$ & $3007.4\,/\;0.106$ & $3657.1\,/\;0.012$\\
\midrule\bottomrule[1pt]
\end{tabular}
\end{sidewaystable}
\endgroup

Table~\ref{tab:fixed-raw-radius-neighborhood} shows that increasing the radius raises the objective and lowers neighborhood risk for both robust models, with a larger response under the C-DRCCP model. At radii of at least $0.015$, the C-DRCCP model has lower risk than the DRCCP model in every comparable cell for the central and high queries. At the low query, the ordering is generally reversed because the DRCCP decision does not depend on the query and is already conservative there. For example, for $(F,D)=(10,50)$, $N=100$, and radius $0.015$, the DRCCP model has objective value $960.4$ at every query, with risks $0.011$, $0.059$, and $0.286$ from the low to the high query. The corresponding objective values for the C-DRCCP model are $887.4$, $997.3$, and $1{,}103.4$, and its risks are $0.038$, $0.040$, and $0.080$. Thus, conditioning reduces protection at the low query and increases it where demand is higher.

\begingroup
\setlength{\rotFPtop}{0pt plus 1fil}
\setlength{\rotFPbot}{0pt plus 1fil}
\begin{sidewaystable}[htbp]
\makeatletter
\let\@makecaption\@maketablecaption
\makeatother
\centering
\caption{Cost and pointwise risk at fixed numerical Wasserstein radii.}
\label{tab:fixed-raw-radius-pointwise}
\scriptsize
\renewcommand{\arraystretch}{1.4}
\setlength{\tabcolsep}{1.8pt}
\begin{tabular}{ccl|cc|cc|cc|cc|cc}\toprule[1pt]\midrule
&&&\multicolumn{2}{c|}{$0.010$}&\multicolumn{2}{c|}{$0.0125$}&\multicolumn{2}{c|}{$0.015$}&\multicolumn{2}{c|}{$0.0175$}&\multicolumn{2}{c}{$0.020$}\\
$(F,D)$ & $N$ & Query & DRCCP & C-DRCCP & DRCCP & C-DRCCP & DRCCP & C-DRCCP & DRCCP & C-DRCCP & DRCCP & C-DRCCP\\\midrule
\multirow{12}{*}{$(10,50)$} & \multirow{3}{*}{50} & Low & $882.5\,/\;0.015$ & $777.7\,/\;0.163$ & $888.4\,/\;0.012$ & $812.1\,/\;0.094$ & $894.0\,/\;0.012$ & $845.3\,/\;0.053$ & $899.0\,/\;0.011$ & $872.9\,/\;0.034$ & $903.9\,/\;0.010$ & $912.9\,/\;0.016$\\
 & & Central & $882.5\,/\;0.148$ & $857.0\,/\;0.199$ & $888.4\,/\;0.136$ & $901.6\,/\;0.117$ & $894.0\,/\;0.133$ & $943.0\,/\;0.079$ & $899.0\,/\;0.124$ & $979.8\,/\;0.044$ & $903.9\,/\;0.114$ & $1013.6\,/\;0.028$\\
 & & High & $882.5\,/\;0.537$ & $947.2\,/\;0.348$ & $888.4\,/\;0.517$ & $988.9\,/\;0.236$ & $894.0\,/\;0.506$ & $1015.8\,/\;0.181$ & $899.0\,/\;0.488$ & $1049.9\,/\;0.130$ & $903.9\,/\;0.469$ & $1080.9\,/\;0.094$\\
\cmidrule(lr){2-13}
 & \multirow{3}{*}{100} & Low & $947.2\,/\;0.002$ & $817.9\,/\;0.056$ & $953.9\,/\;0.002$ & $851.9\,/\;0.024$ & $960.4\,/\;0.001$ & $887.4\,/\;0.011$ & $966.5\,/\;0.001$ & $926.1\,/\;0.005$ & $972.4\,/\;0.001$ & $970.1\,/\;0.002$\\
 & & Central & $947.2\,/\;0.046$ & $896.6\,/\;0.106$ & $953.9\,/\;0.041$ & $942.5\,/\;0.055$ & $960.4\,/\;0.036$ & $997.3\,/\;0.024$ & $966.5\,/\;0.032$ & $1044.6\,/\;0.013$ & $972.4\,/\;0.029$ & $1084.2\,/\;0.007$\\
 & & High & $947.2\,/\;0.291$ & $1020.5\,/\;0.149$ & $953.9\,/\;0.270$ & $1066.4\,/\;0.088$ & $960.4\,/\;0.251$ & $1103.4\,/\;0.055$ & $966.5\,/\;0.235$ & $1142.4\,/\;0.035$ & $972.4\,/\;0.220$ & $1184.0\,/\;0.023$\\
\cmidrule(lr){2-13}
 & \multirow{3}{*}{200} & Low & $966.9\,/\;0.000$ & $841.6\,/\;0.013$ & $973.8\,/\;0.000$ & $886.2\,/\;0.004$ & $980.1\,/\;0.000$ & $946.6\,/\;0.001$ & $985.9\,/\;0.000$ & $992.7\,/\;0.000$ & $991.8\,/\;0.000$ & $1043.8\,/\;0.000$\\
 & & Central & $966.9\,/\;0.024$ & $925.6\,/\;0.052$ & $973.8\,/\;0.023$ & $992.1\,/\;0.019$ & $980.1\,/\;0.020$ & $1043.2\,/\;0.009$ & $985.9\,/\;0.018$ & $1079.8\,/\;0.005$ & $991.8\,/\;0.016$ & $1119.7\,/\;0.002$\\
 & & High & $966.9\,/\;0.197$ & $1064.3\,/\;0.056$ & $973.8\,/\;0.192$ & $1109.4\,/\;0.031$ & $980.1\,/\;0.182$ & $1148.2\,/\;0.019$ & $985.9\,/\;0.169$ & $1195.3\,/\;0.012$ & $991.8\,/\;0.156$ & $1279.0\,/\;0.007$\\
\cmidrule(lr){2-13}
 & \multirow{3}{*}{300} & Low & $989.4\,/\;0.000$ & $857.7\,/\;0.009$ & $997.4\,/\;0.000$ & $927.4\,/\;0.001$ & $1004.7\,/\;0.000$ & $979.6\,/\;0.000$ & $1011.6\,/\;0.000$ & $1024.7\,/\;0.000$ & $1018.0\,/\;0.000$ & $1080.3\,/\;0.000$\\
 & & Central & $989.4\,/\;0.012$ & $961.4\,/\;0.024$ & $997.4\,/\;0.011$ & $1017.7\,/\;0.007$ & $1004.7\,/\;0.008$ & $1073.9\,/\;0.003$ & $1011.6\,/\;0.007$ & $1112.2\,/\;0.001$ & $1018.0\,/\;0.006$ & $1152.9\,/\;0.001$\\
 & & High & $989.4\,/\;0.134$ & $1086.5\,/\;0.035$ & $997.4\,/\;0.126$ & $1142.9\,/\;0.021$ & $1004.7\,/\;0.110$ & $1181.3\,/\;0.012$ & $1011.6\,/\;0.104$ & $1223.4\,/\;0.007$ & $1018.0\,/\;0.095$ & --\\
\midrule
\multirow{12}{*}{$(15,100)$} & \multirow{3}{*}{50} & Low & $2522.3\,/\;0.015$ & $2093.3\,/\;0.211$ & $2539.6\,/\;0.013$ & $2234.8\,/\;0.111$ & $2556.7\,/\;0.011$ & $2351.6\,/\;0.062$ & $2573.4\,/\;0.010$ & $2475.0\,/\;0.037$ & $2589.2\,/\;0.010$ & $2577.7\,/\;0.024$\\
 & & Central & $2522.3\,/\;0.157$ & $2382.4\,/\;0.268$ & $2539.6\,/\;0.143$ & $2521.0\,/\;0.148$ & $2556.7\,/\;0.127$ & $2639.6\,/\;0.081$ & $2573.4\,/\;0.117$ & $2749.7\,/\;0.045$ & $2589.2\,/\;0.110$ & $2853.7\,/\;0.024$\\
 & & High & $2522.3\,/\;0.568$ & $2499.9\,/\;0.610$ & $2539.6\,/\;0.545$ & $2599.3\,/\;0.496$ & $2556.7\,/\;0.517$ & $2702.5\,/\;0.387$ & $2573.4\,/\;0.498$ & $2805.1\,/\;0.290$ & $2589.2\,/\;0.480$ & $2897.5\,/\;0.218$\\
\cmidrule(lr){2-13}
 & \multirow{3}{*}{100} & Low & $2670.3\,/\;0.003$ & $2247.6\,/\;0.081$ & $2691.3\,/\;0.002$ & $2382.3\,/\;0.040$ & $2710.1\,/\;0.002$ & $2502.4\,/\;0.017$ & $2727.1\,/\;0.002$ & $2628.3\,/\;0.008$ & $2743.4\,/\;0.002$ & $2748.6\,/\;0.004$\\
 & & Central & $2670.3\,/\;0.051$ & $2573.5\,/\;0.102$ & $2691.3\,/\;0.045$ & $2733.1\,/\;0.044$ & $2710.1\,/\;0.042$ & $2835.8\,/\;0.024$ & $2727.1\,/\;0.039$ & $2961.5\,/\;0.012$ & $2743.4\,/\;0.035$ & $3084.3\,/\;0.006$\\
 & & High & $2670.3\,/\;0.335$ & $2774.6\,/\;0.275$ & $2691.3\,/\;0.315$ & $2873.0\,/\;0.201$ & $2710.1\,/\;0.301$ & $2989.8\,/\;0.139$ & $2727.1\,/\;0.285$ & $3090.9\,/\;0.092$ & $2743.4\,/\;0.269$ & $3201.8\,/\;0.060$\\
\cmidrule(lr){2-13}
 & \multirow{3}{*}{200} & Low & $2891.1\,/\;0.001$ & $2530.2\,/\;0.010$ & $2912.7\,/\;0.001$ & $2697.5\,/\;0.002$ & $2934.2\,/\;0.001$ & $2812.7\,/\;0.001$ & $2953.5\,/\;0.001$ & $2952.9\,/\;0.000$ & $2972.7\,/\;0.001$ & $3110.5\,/\;0.000$\\
 & & Central & $2891.1\,/\;0.007$ & $2747.8\,/\;0.026$ & $2912.7\,/\;0.007$ & $2971.7\,/\;0.007$ & $2934.2\,/\;0.006$ & $3093.9\,/\;0.004$ & $2953.5\,/\;0.006$ & $3206.4\,/\;0.002$ & $2972.7\,/\;0.005$ & $3309.2\,/\;0.001$\\
 & & High & $2891.1\,/\;0.120$ & $3062.0\,/\;0.072$ & $2912.7\,/\;0.108$ & $3183.0\,/\;0.046$ & $2934.2\,/\;0.098$ & $3285.5\,/\;0.028$ & $2953.5\,/\;0.091$ & $3441.2\,/\;0.014$ & $2972.7\,/\;0.084$ & $3597.8\,/\;0.011$\\
\cmidrule(lr){2-13}
 & \multirow{3}{*}{300} & Low & $2923.6\,/\;0.001$ & $2545.2\,/\;0.007$ & $2946.6\,/\;0.001$ & $2698.3\,/\;0.002$ & $2968.8\,/\;0.001$ & $2805.4\,/\;0.001$ & $2988.1\,/\;0.001$ & $2928.4\,/\;0.000$ & $3007.4\,/\;0.001$ & $3112.7\,/\;0.000$\\
 & & Central & $2923.6\,/\;0.006$ & $2809.5\,/\;0.016$ & $2946.6\,/\;0.006$ & $3020.5\,/\;0.003$ & $2968.8\,/\;0.005$ & $3157.5\,/\;0.002$ & $2988.1\,/\;0.004$ & $3252.8\,/\;0.001$ & $3007.4\,/\;0.004$ & $3358.8\,/\;0.001$\\
 & & High & $2923.6\,/\;0.102$ & $3123.3\,/\;0.051$ & $2946.6\,/\;0.092$ & $3230.3\,/\;0.033$ & $2968.8\,/\;0.089$ & $3334.2\,/\;0.021$ & $2988.1\,/\;0.082$ & $3459.5\,/\;0.011$ & $3007.4\,/\;0.077$ & $3657.1\,/\;0.006$\\
\midrule\bottomrule[1pt]
\end{tabular}
\end{sidewaystable}
\endgroup

Table~\ref{tab:fixed-raw-radius-pointwise} evaluates the same decisions at $X=x_0$ and reproduces the neighborhood pattern. In the representative cell above, the pointwise risks from the low to the high query are $0.001$, $0.036$, and $0.251$ for the DRCCP model and $0.011$, $0.024$, and $0.055$ for the C-DRCCP model. The agreement between the two risk measures shows that the reallocation is not an artifact of averaging over the conditioning neighborhood. The $N=200$ and $N=300$ rows present the same qualitative behavior, but each uses one nested replication and is therefore diagnostic rather than replicated evidence about sample size.

\subsection{Computational Comparison}\label{appendix:inequality-implementation}
We compare three exact implementations of the C-DRCCP model with the strengthening constraints in \texttt{ALL}. The \emph{compact} implementation is formulation~\eqref{opt:exact-mip} with its dual block. The \emph{external} and \emph{lazy} implementations replace that block with inequalities from~\eqref{eq:critical-radius-cuts}. The external implementation adds violated inequalities between successive MIP solves, whereas the lazy implementation adds them through a callback within one branch-and-bound tree.

We include additional network and sample sizes used only for this implementation comparison. All runs terminate with zero objective-bound gap and pass independent feasibility checks. We report total algorithm time $T^{\rm alg}$ using the timing convention defined in Section~\ref{sec:computational}. The first four rows of Table~\ref{tab:constraint-generation-performance} report medians over six paired replications, the next four report medians over two, and the final row is a single stress instance. The ratio column reports the median of the paired ratios.

\begin{table}[htbp]
\centering
\caption{Computational comparison of three exact implementations.}
\label{tab:constraint-generation-performance}
\scriptsize
\setlength{\tabcolsep}{2.5pt}
\resizebox{0.6\linewidth}{!}{%
\begin{tabular}{ccrrrrr}\toprule[1pt]\midrule
& & \multicolumn{3}{c}{Median $T^{\rm alg}$ [s]} & & \\
\cmidrule(lr){3-5}
 $(F,D,N)$ & $a$ & Compact & External & Lazy & $T_L^{\rm alg}/T_C^{\rm alg}$ & Lazy faster\\
\midrule
 $(15,100,140)$ & 0.10 & 3.060 & \textbf{1.813} & 2.352 & 0.771 & 6/6\\
 $(15,300,100)$ & 0.10 & 4.424 & 8.054 & \textbf{4.106} & 0.928 & 6/6\\
 $(15,300,150)$ & 0.05 & 10.131 & 17.085 & \textbf{7.784} & 0.769 & 6/6\\
 $(15,800,100)$ & 0.10 & 12.264 & \textbf{10.053} & 10.159 & 0.829 & 6/6\\
 $(15,300,125)$ & 0.10 & \textbf{5.818} & -- & 6.465 & 1.111 & 0/2\\
 $(15,300,130)$ & 0.10 & \textbf{6.334} & -- & 7.268 & 1.148 & 0/2\\
 $(15,300,150)$ & 0.02 & 6.899 & -- & \textbf{6.861} & 0.994 & 2/2\\
 $(15,300,200)$ & 0.02 & \textbf{16.584} & -- & 17.446 & 1.052 & 0/2\\
 $(15,100,500)$ & 0.50 & \textbf{19.867} & -- & 262.121 & 13.194 & 0/1\\
\midrule\bottomrule[1pt]
\end{tabular}%
}
\end{table}

The lazy implementation has a lower median algorithm time than the compact implementation in five of the nine settings. In the first four settings, it is faster in all 24 paired replications and reduces the median time by $7.2\%$--$23.1\%$. Among the four settings for which all three implementations are reported, the external and lazy implementations each attain the lowest median time twice. In the other five settings, lazy enforcement is essentially tied with compact at $(F,D,N)=(15,300,150)$ and $a=0.02$, is $5.2\%$--$14.8\%$ slower in three settings, and is $13.2$ times slower in the $N=500$ stress instance. Thus, lazy enforcement can avoid the cost of the compact dual block and the repeated solves of the external implementation, but its advantage depends on the number of inequalities needed during the search.
\end{APPENDICES}

\end{document}